\documentclass[acmtog]{acmart}

\usepackage[shortlabels]{enumitem}
\usepackage{algorithm}
\usepackage{algorithmic}
\usepackage{multirow}

\newtheorem{theorem}{Theorem}[section]

\newtheorem{lemma}[theorem]{Lemma}

\newtheorem{problem}[theorem]{Problem}

\theoremstyle{definition}

\theoremstyle{definition}
\newtheorem{remark}[theorem]{Remark}

\newcommand{\RN}[1]{%
	\textup{\uppercase\expandafter{\romannumeral#1}}%
}
\def\*#1{\mathbf{#1}}
\newcommand{\half}{{k + \frac{1}{2}}}
\newcommand{\divg}{\mathrm{div}_g}
\newcommand{\divgh}{\mathrm{div}_{g_h}}
\newcommand{\IP}[2]{\left\langle #1, #2 \right\rangle_g}
\newcommand{\IPA}[2]{\left\langle #1, #2 \right\rangle}
\newcommand{\IPH}[2]{\left\langle #1, #2 \right\rangle_{g_h}}

\definecolor{dollarbill}{rgb}{0.52, 0.73, 0.4}
\definecolor{asparagus}{rgb}{0.53, 0.66, 0.42}
\definecolor{applegreen}{rgb}{0.55, 0.71, 0.0}
\newcommand{\greener}[1]{\textcolor{black}{#1}}
\newcommand{\green}[1]{\textcolor{black}{#1}}

\AtBeginDocument{%
  \providecommand\BibTeX{{%
    \normalfont B\kern-0.5em{\scshape i\kern-0.25em b}\kern-0.8em\TeX}}}

\setcopyright{acmcopyright}
\acmJournal{TOG}
\acmYear{2020} \acmVolume{1} \acmNumber{1} \acmArticle{1} \acmMonth{1} \acmPrice{15.00}\acmDOI{10.1145/3369387}

\begin{document}

\title{\green{Computational} p-Willmore Flow with \green{Conformal Penalty}}

\author{Anthony Gruber}
\email{anthony.gruber@ttu.edu}
\orcid{0000-0001-7107-5307}

\author{Eugenio Aulisa}
\email{eugenio.aulisa@ttu.edu}

\affiliation{%
  \institution{Texas Tech University}
  \streetaddress{P.O. Box 41042}
  \city{Lubbock}
  \state{Texas}
  \postcode{79409}
}

\begin{abstract}
    The unsigned p-Willmore functional introduced in \cite{mondino2011} generalizes important geometric functionals which measure the area and Willmore energy of immersed surfaces.  Presently, techniques from \cite{dziuk2008} are adapted to compute the first variation of this functional as a weak-form system of equations, which are subsequently used to develop a model for the p-Willmore flow of closed surfaces in $\mathbb{R}^3$.  This model is amenable to constraints on surface area and enclosed volume, and is shown to decrease the p-Willmore energy monotonically.  In addition, a penalty-based regularization procedure is formulated to prevent artificial mesh degeneration along the flow; inspired by a conformality condition derived in \cite{kamberov1996}, \green{this procedure encourages angle-preservation in a closed and oriented surface immersion as it evolves.  Following this, a finite-element discretization of both procedures is discussed, an algorithm for running the flow is given, and an application to mesh editing is presented.}
\end{abstract}




\begin{teaserfigure}
\begin{center}
    \begin{minipage}[c]{0.18\textwidth}
    \includegraphics[width=\textwidth]{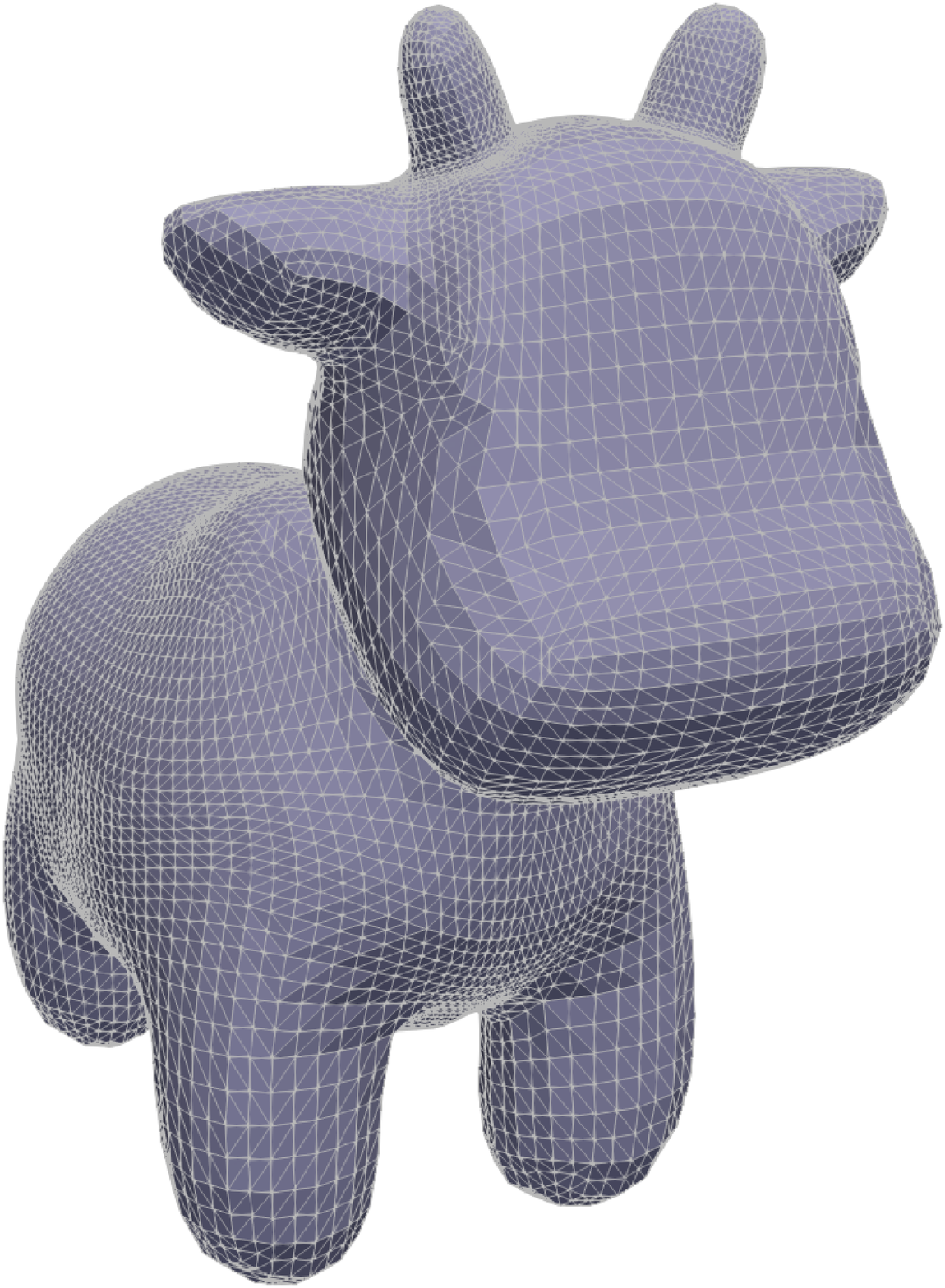}
    \end{minipage}
    \begin{minipage}[c]{0.18\textwidth}
    \includegraphics[width=\textwidth]{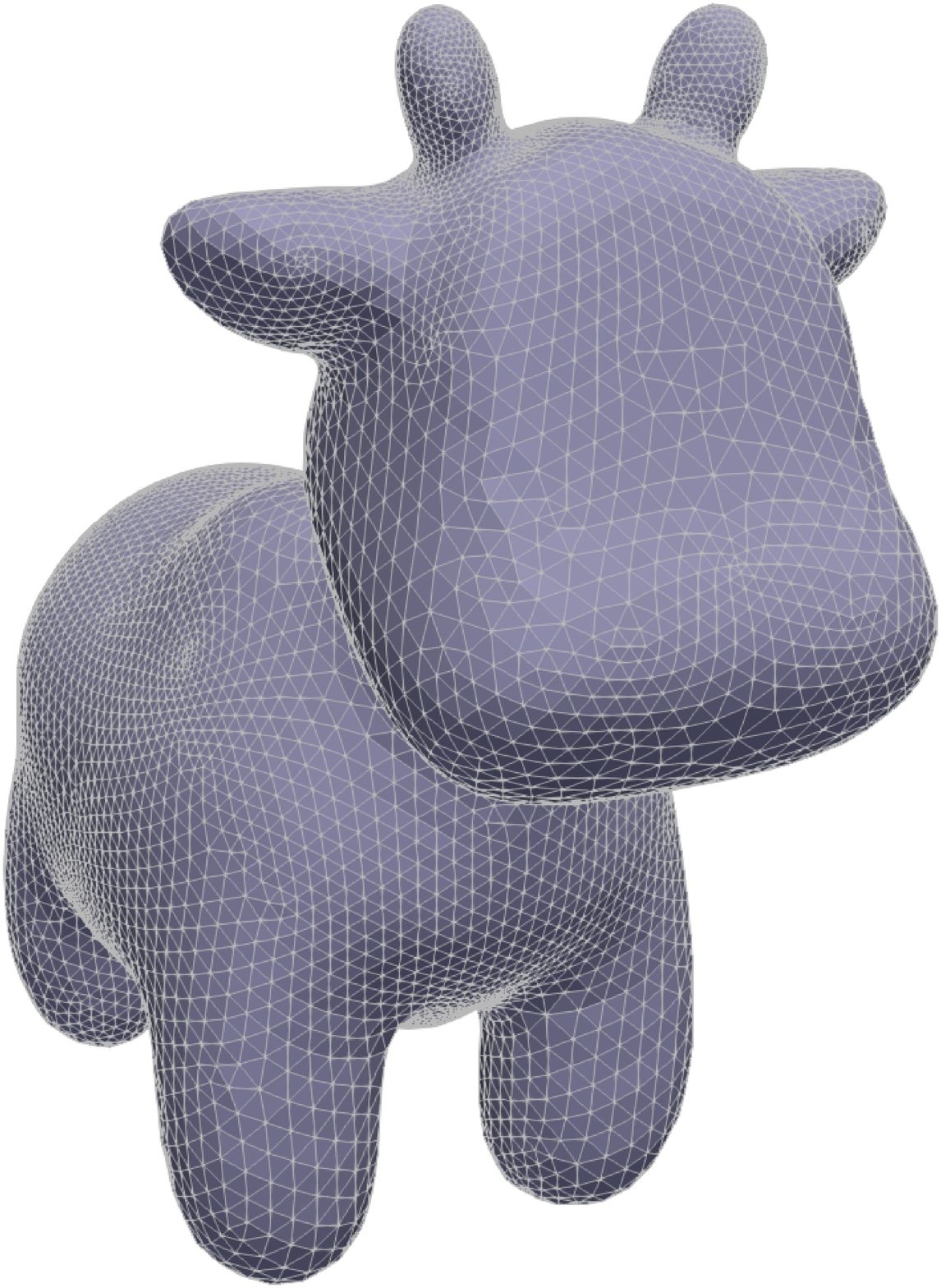}
    \end{minipage}
    \begin{minipage}[c]{0.18\textwidth}
    \includegraphics[width=\textwidth]{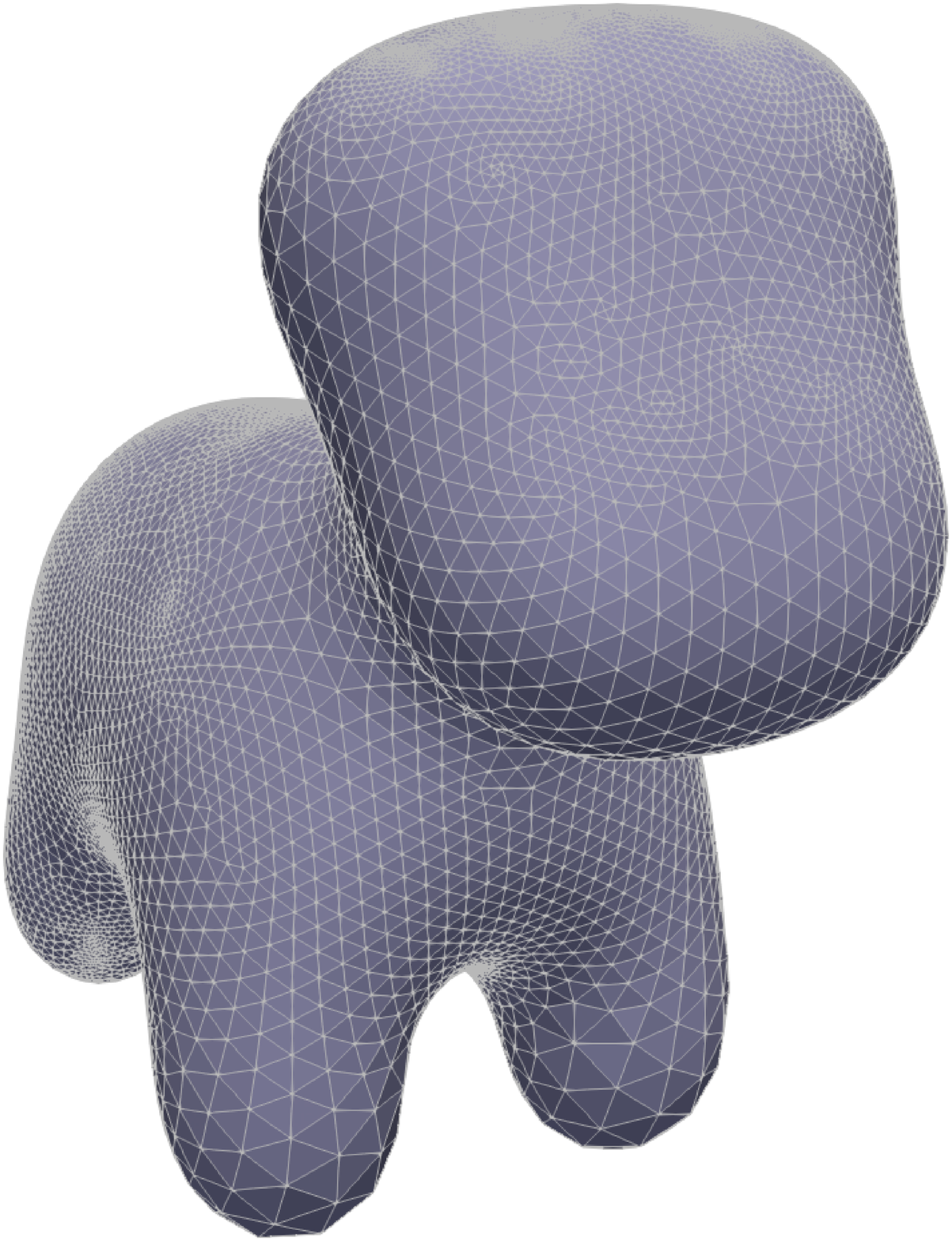}
    \end{minipage}
    \begin{minipage}[c]{0.18\textwidth}
    \includegraphics[width=\textwidth]{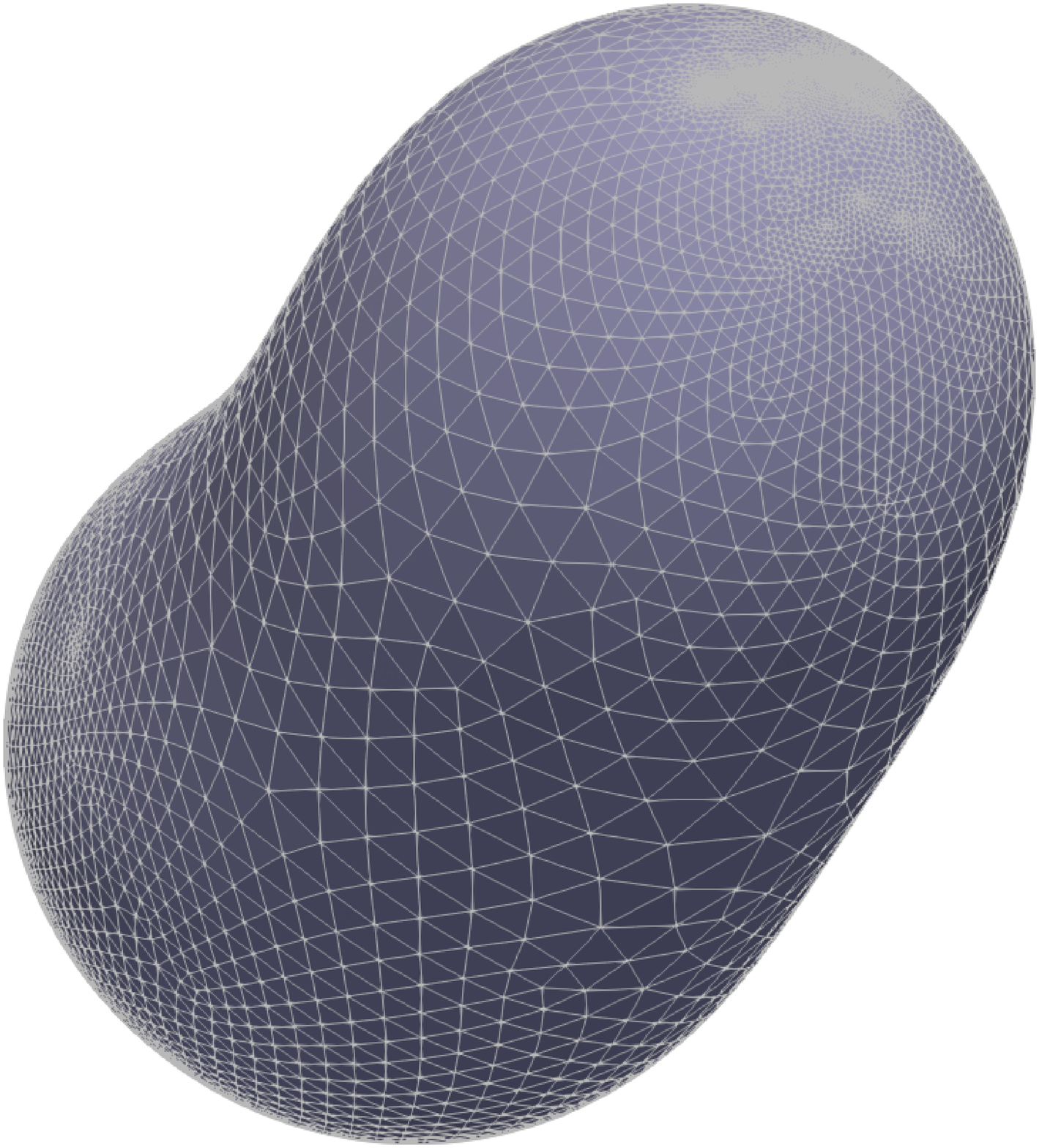}
    \end{minipage}
    \begin{minipage}[c]{0.20\textwidth}
    \includegraphics[width=\textwidth]{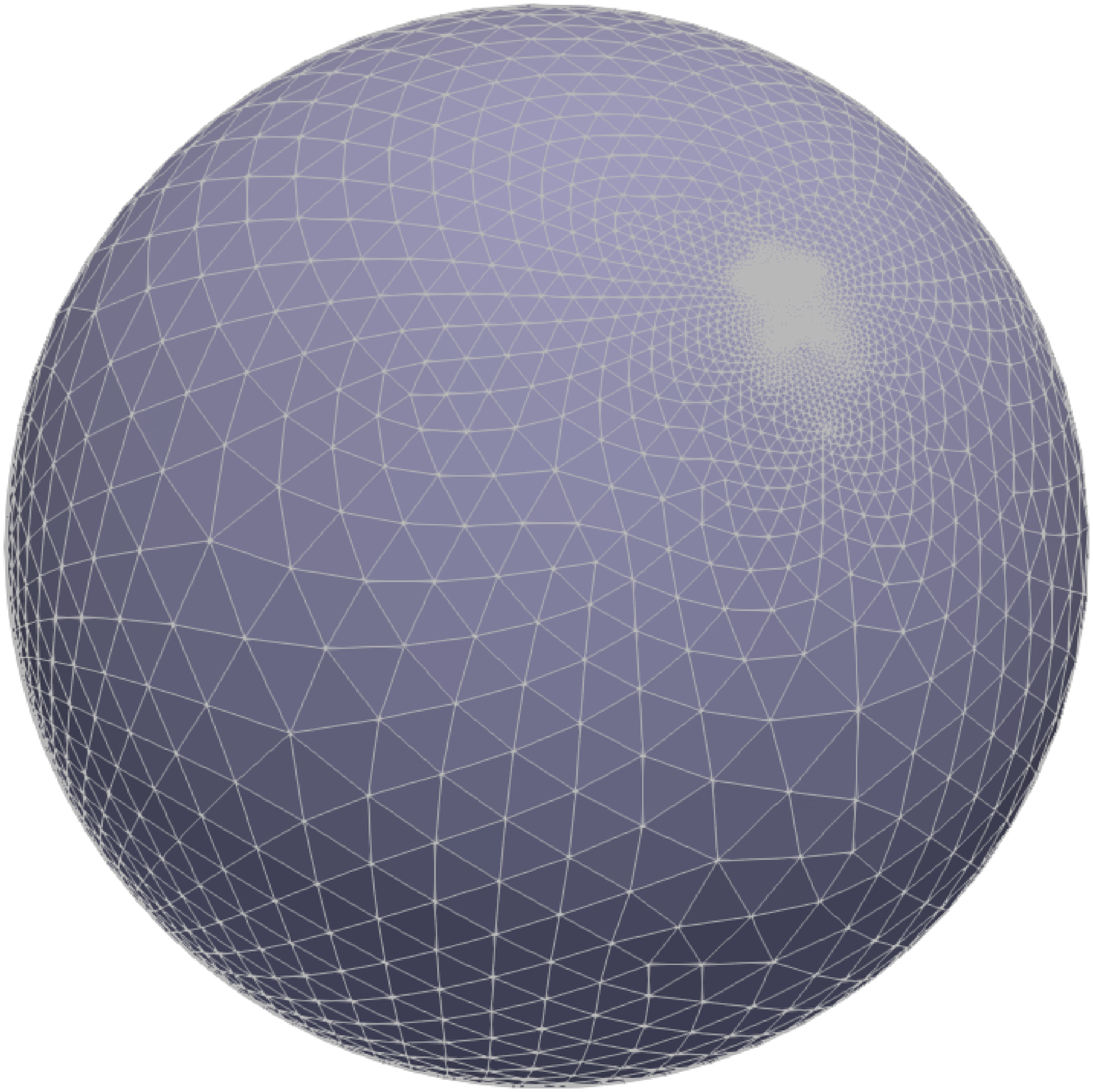}
    \end{minipage}
\end{center}
\caption[2-Willmore flow with conformal penalty]{Area-preserving 2-Willmore flow with conformal penalty applied to a cow mesh of 23.4k triangles.  Time steps pictured: 0,1,50,90,130. Area change <0.3\%.}
\label{fig:funcow}
\end{teaserfigure}

\maketitle

\section{Introduction}

As another example, the reliable \green{Helfrich-Canham} model for biomembranes (see \cite{helfrich1973}) is based on the well-studied Willmore energy (see \cite{willmore1965,white1973,weiner1978,bohle2008,mondino2011,marques2014,athukorallage2015} and references therein)
\begin{equation*}
    \mathcal{W}^2(u) = \int_{M} H^2\,\green{d\mu_g},
\end{equation*}
whose $L^2$-gradient flow has been proven to converge smoothly to a global minimum when the surface genus and initial energy are sufficiently low \cite{kuwert2001,mondino2014} \green{(c.f. Figure~\ref{fig:torus}).}  Due to its pleasing aesthetic character, the Willmore flow has further attracted the interest of computational mathematicians and scientists, and has been studied numerically in a  variety of contexts including conformal geometry, geometric partial differential equations, and computer graphics.  See e.g. \cite{crane2013,dziuk2013,joshi2007} and the references therein.

\subsection{Related work}
Besides the inherent mathematical challenges present in geometric flows (involving e.g. convergence, changes in global topology, and singularity formation), \greener{their governing equations} introduce a number of computational difficulties as well.  In particular, discrete surfaces are often stored as piecewise-linear data, such as meshes of simplices, and it is taxing to find a satisfactory method of expression for second-order geometric phenomena such as curvature.  There have been two broad approaches to this problem in the current literature, which can be thought of colloquially as arising from \textit{discrete} versus \textit{discretized} perspectives on the issue.

In discrete geometry, the aim is to use global characterizations from geometry and topology to develop fully-discrete analogues of classical geometric quantities, which are in some sense independent from their original (continuous) definitions.  Tools such as exterior calculus, the Gauss-Bonnet and Stokes' Theorems are employed to define length, area, curvature, etc. on a simplicial surface, which is accomplished through enforcing global geometric relationships rather than considering local values at specific places (nodes) on a mesh.  The main advantages of this approach are relative independence from mesh quality, and sparse  linear formulations which are fast to solve.  Some notable disadvantages present here are the restriction of such methods (so far) to triangular meshes, and the fact that several equivalent definitions of geometric quantities in the smooth setting become inequivalent when treated in this way (see \cite{crane2017} for details).  Further information on this area can be found in \cite{droske2004,deckelnick2006,meyer2003,bobenko2008,gu2009} and the references therein.

Conversely, discretized geometry involves approximating continuous geometric quantities as well as possible by using a good choice of nodal mesh points, so that the difference between the continuous and discrete objects vanishes in the limit of mesh refinement.  Traditional finite element mathematics is based on this idea, whereby the necessary calculations are done locally and element-wise without any particular adherence to global phenomena except in the limit.  The primary advantage of this approach is its flexibility with respect to applications, problem formulations, and mesh data.  Its main disadvantages are its inherent sensitivity to mesh quality, and its agnosticism with respect to the global aspects of surface geometry.  See \cite{dziuk2013} for a compendium of knowledge and techniques in this area.
\begin{remark}
In fact, the failure of the finite element method to capture global relationships was a primary motivation for the development of a discrete geometric theory, as mentioned in \cite{bobenko2008,gu2009}. 
\end{remark}

\green{Due in part to their useful application to problems such as mesh editing (see \cite{bobenko2005}),} the computational details of geometric flows have been examined previously from both of the above perspectives.  In \cite{dziuk2008}, the author studies parametric Willmore flow using finite-element methods.  In particular, \green{the author} develops and discretizes a model for the Willmore flow of surfaces, detailing some examples and proving \green{stability} of this discretization.  On the other hand, the authors in \cite{crane2017} use ideas from discrete conformal geometry to develop a conformally-constrained model for the Willmore flow.  More precisely, they develop results which enable the direct manipulation of surface curvature, allowing for angle-preserving mesh positions to be recovered using a natural integrability condition. Beyond the Willmore flow, many computational studies have also been done which focus on the mean and Gauss curvature flows, Ricci flow, and Yamabe flow of surfaces; see \cite{deckelnick2005,joshi2007} and their enclosed references for more details.

This work adopts a discretized perspective similar to \cite{dziuk2008,dziuk2013} and aims to extend the computational study of curvature flows \green{that arise from functionals which depend on some power of the mean curvature of an immersed surface}.  To that end, the main object of study is the $L^2$-gradient flow of the \green{(unsigned)} \textit{p-Willmore functional} introduced in \green{\cite{mondino2011},}
\begin{equation*}
    \mathcal{W}^p(u) = \frac{1}{2^p}\int_{M} |H|^p\,\green{d\mu_g, \qquad p\geq 1.}
\end{equation*}
\green{As mentioned in \cite{gruber2019}, this definition can be extended to include the case $p=0$, so that the surface area, (unsigned) total mean curvature, and Willmore functionals are encompassed here as $\mathcal{W}^0$, $\mathcal{W}^1$, and $\mathcal{W}^2$, respectively.  It follows that the 0-Willmore flow is simply MCF, and usual Willmore flow occurs when $p=2$.}

\begin{figure}
\begin{center}
    \begin{minipage}[c]{0.12\textwidth}
    MCF \\ (0-Willmore)
    \end{minipage}
    \begin{minipage}[c]{0.12\textwidth}
    Willmore flow \\ (2-Willmore)
    \end{minipage}
    \begin{minipage}[c]{0.12\textwidth}
    4-Willmore flow \\
    \end{minipage}
    \\
    \vspace{0.2pc}
    \begin{minipage}[c]{0.15\textwidth}
    \includegraphics[width=\textwidth]{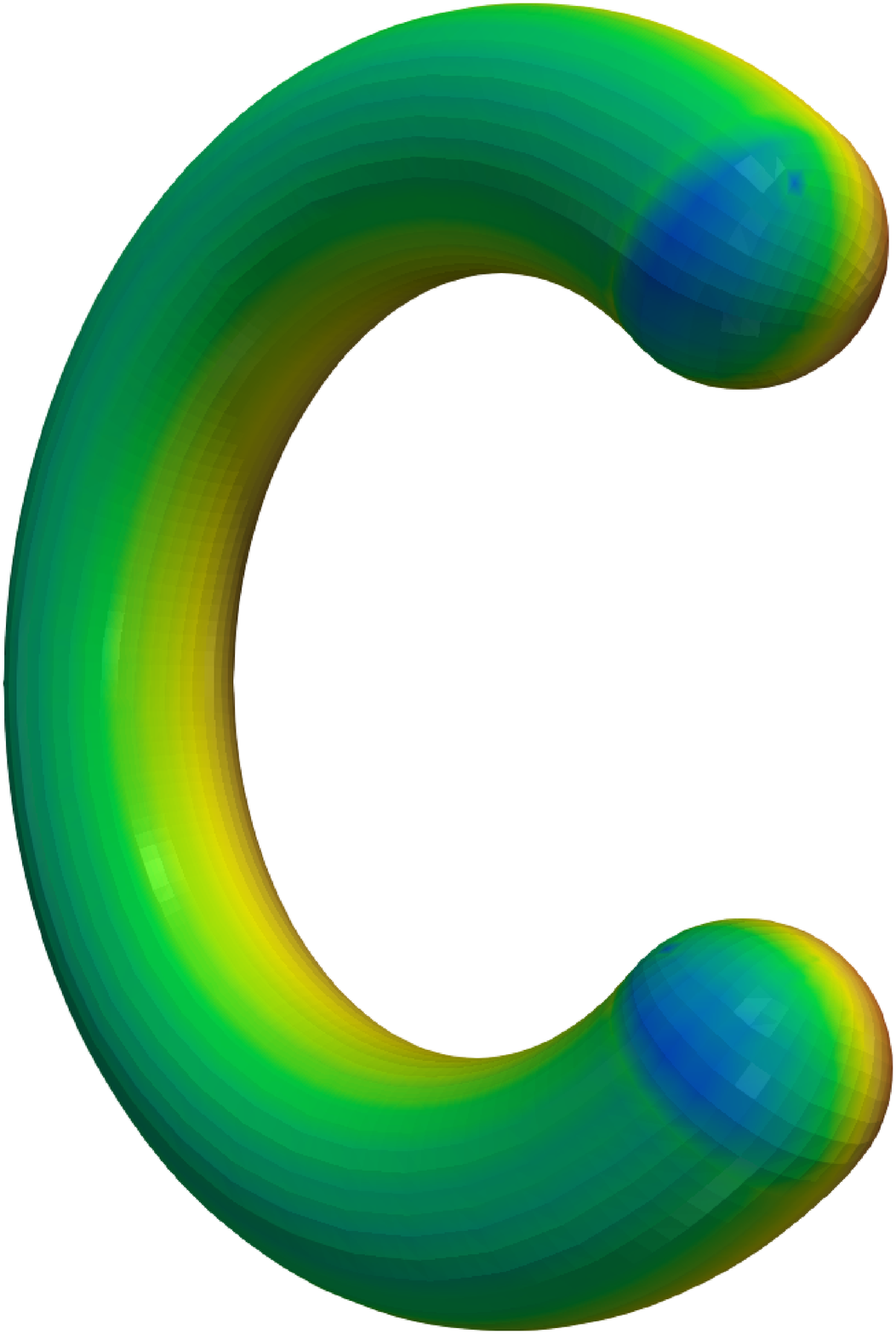}
    \end{minipage}
    \begin{minipage}[c]{0.15\textwidth}
    \includegraphics[width=\textwidth]{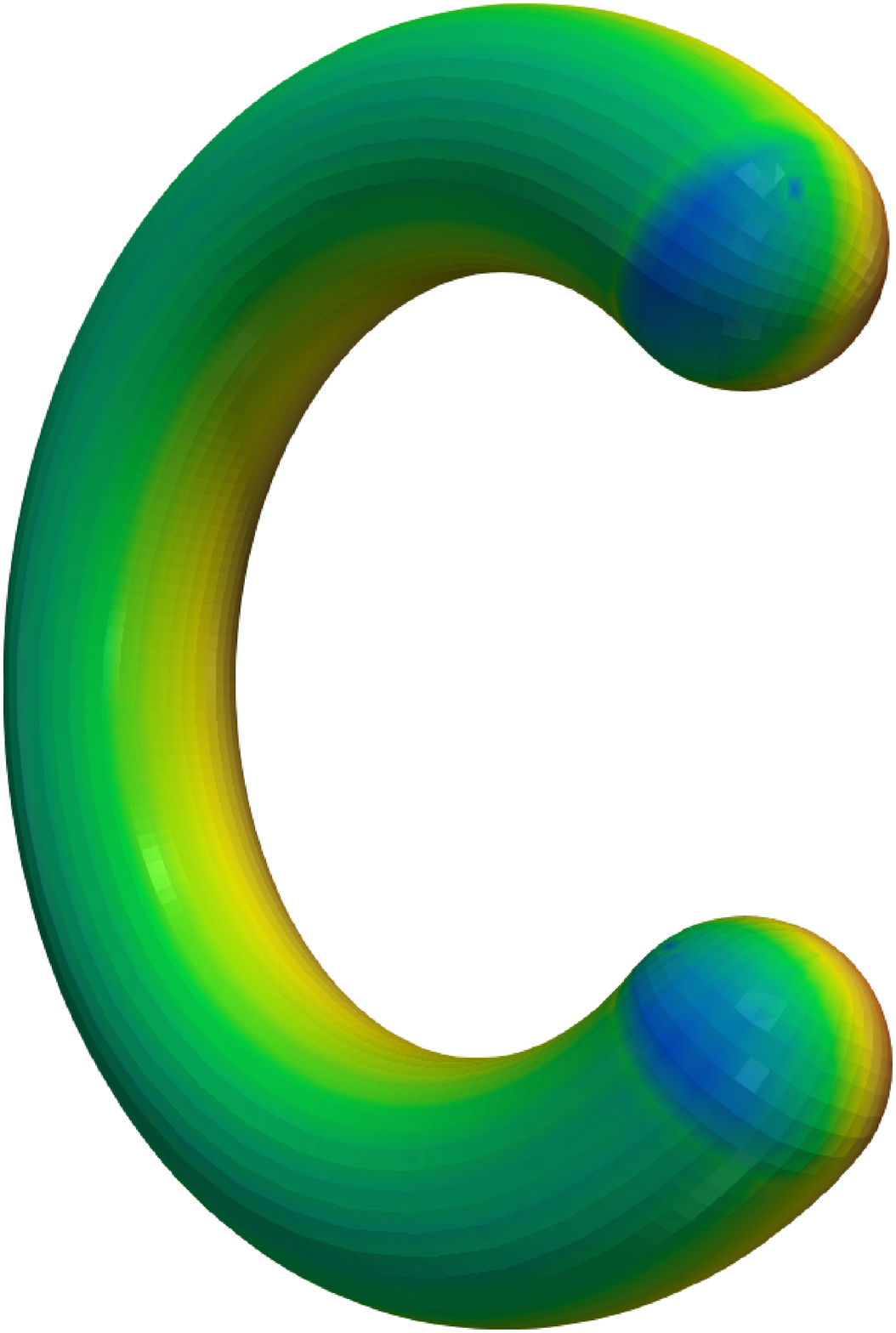}
    \end{minipage}
    \begin{minipage}[c]{0.15\textwidth}
    \includegraphics[width=\textwidth]{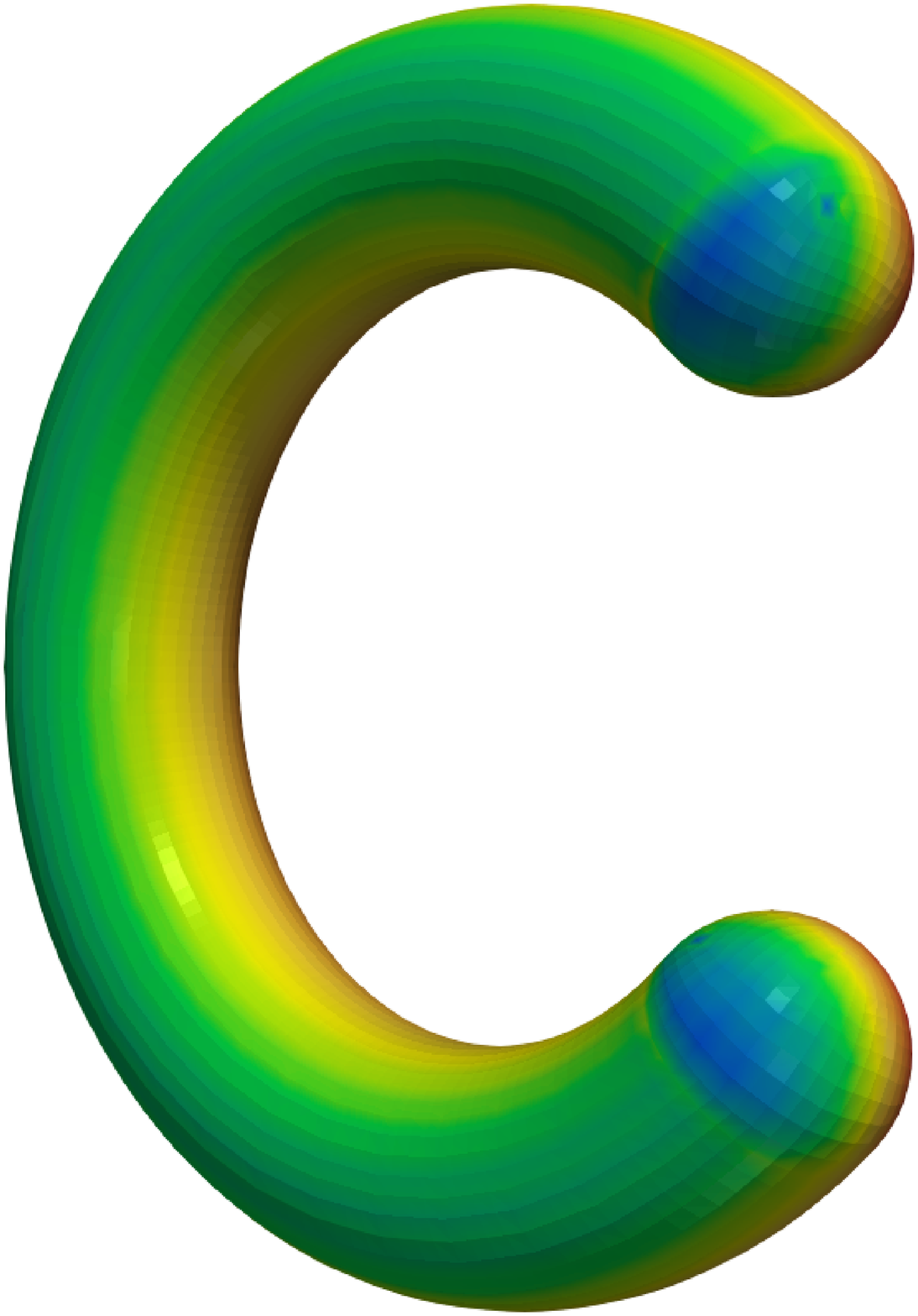}
    \end{minipage}
    \\
    \vspace{0.4pc}
    \begin{minipage}[c]{0.15\textwidth}
    \includegraphics[width=\textwidth]{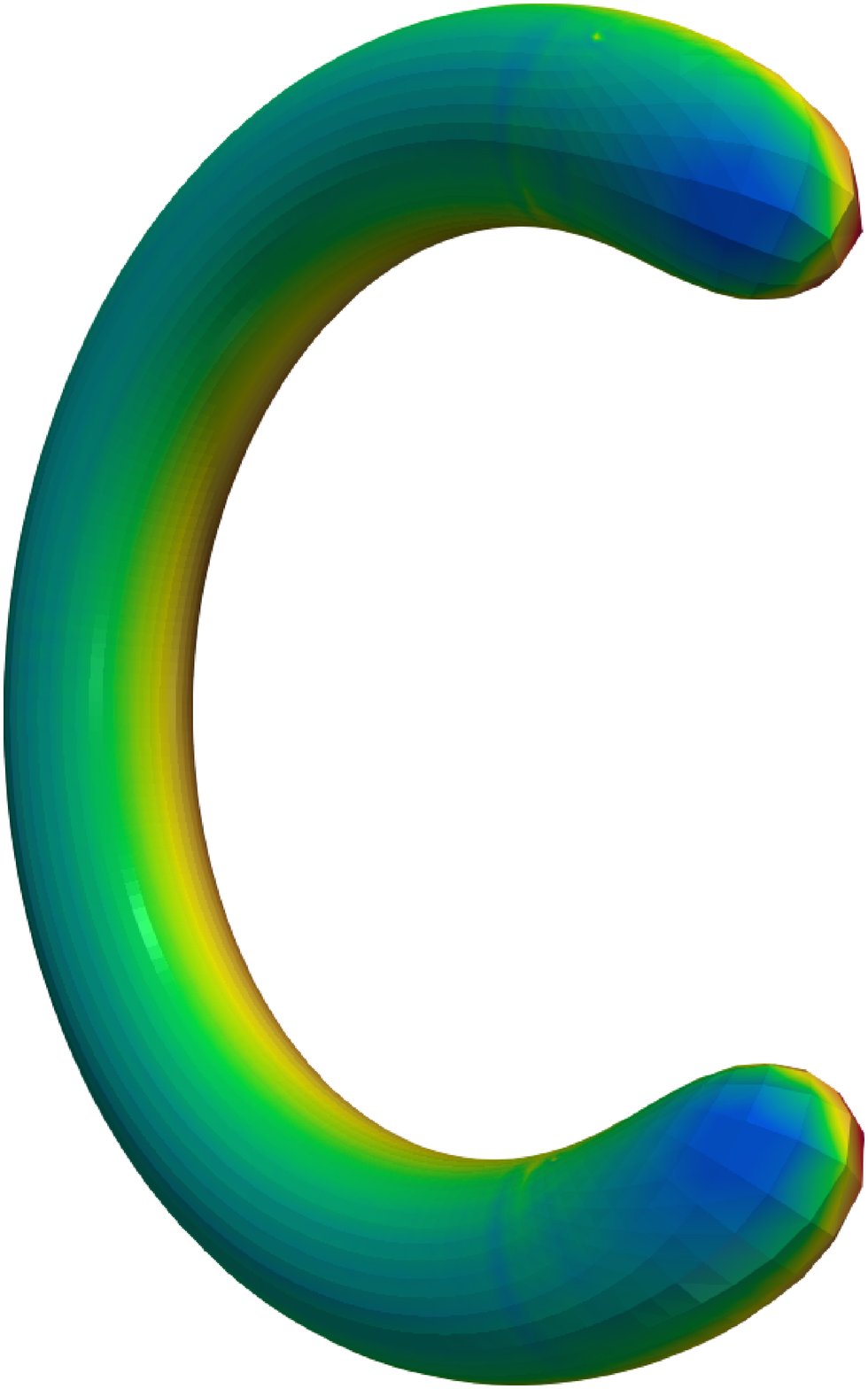}
    \end{minipage}
    \begin{minipage}[c]{0.15\textwidth}
    \includegraphics[width=\textwidth]{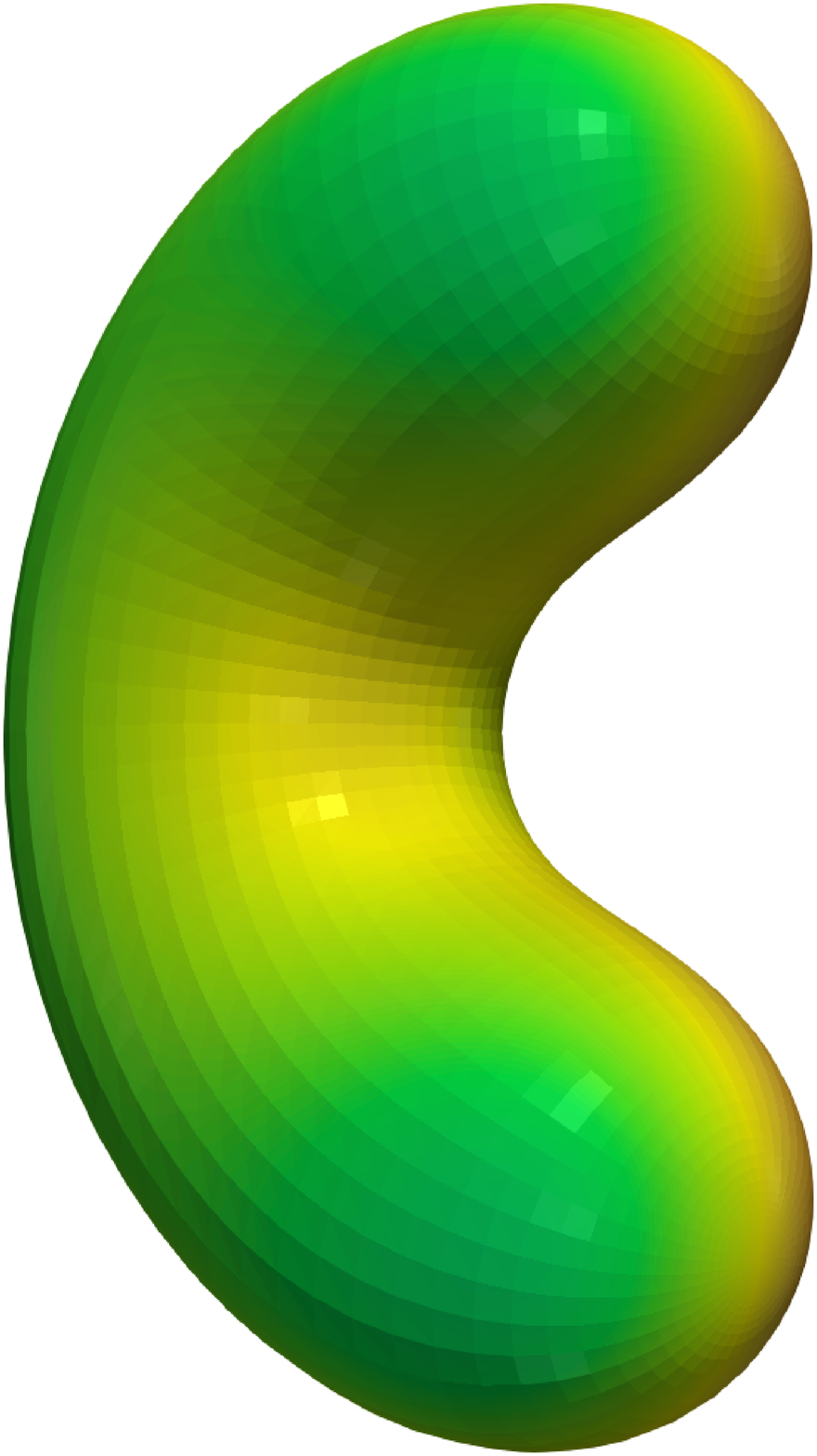}
    \end{minipage}
    \begin{minipage}[c]{0.16\textwidth}
    \includegraphics[width=\textwidth]{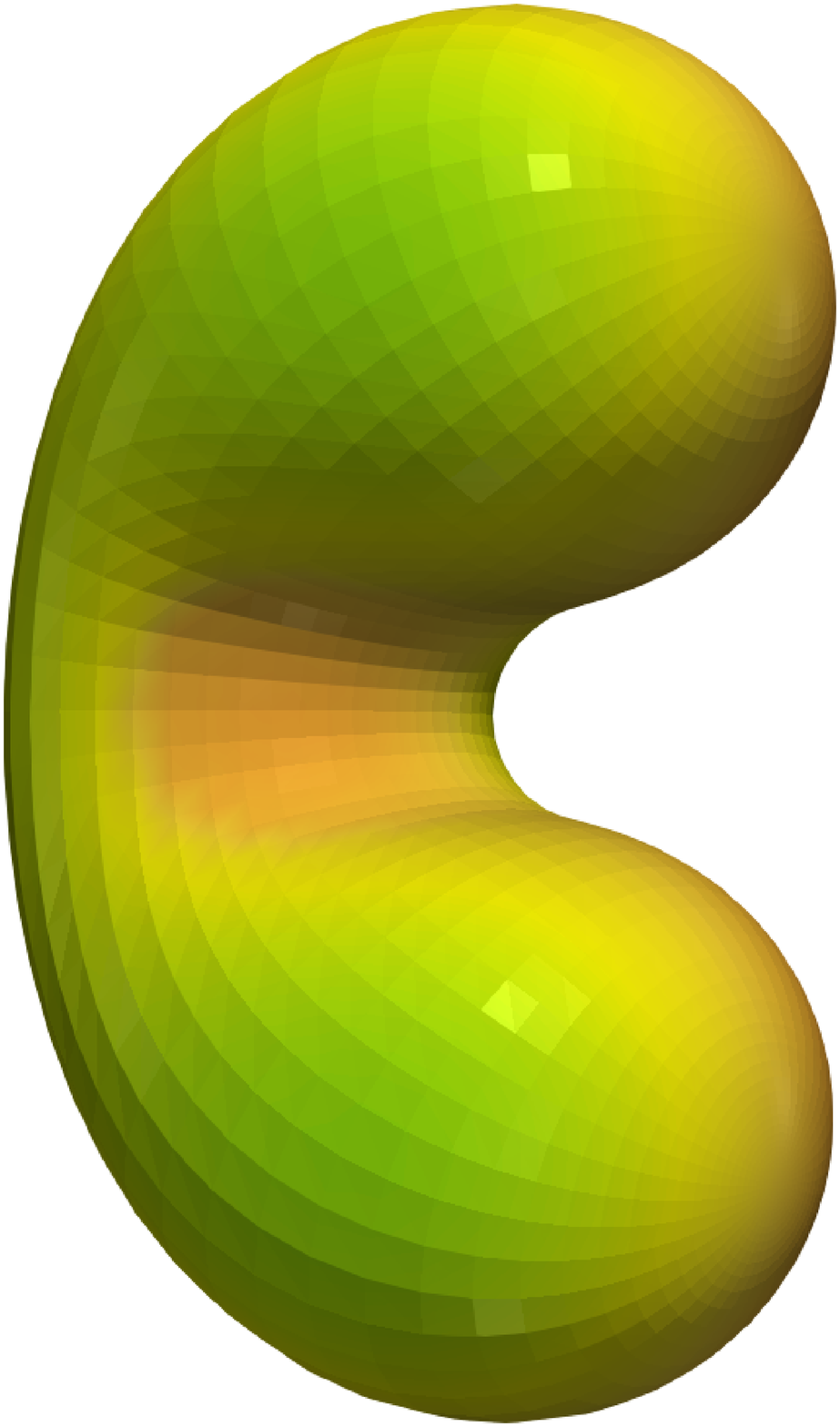}
    \end{minipage}
    \\
    \vspace{0.2pc}
    \begin{minipage}[c]{0.15\textwidth}
    \includegraphics[width=\textwidth]{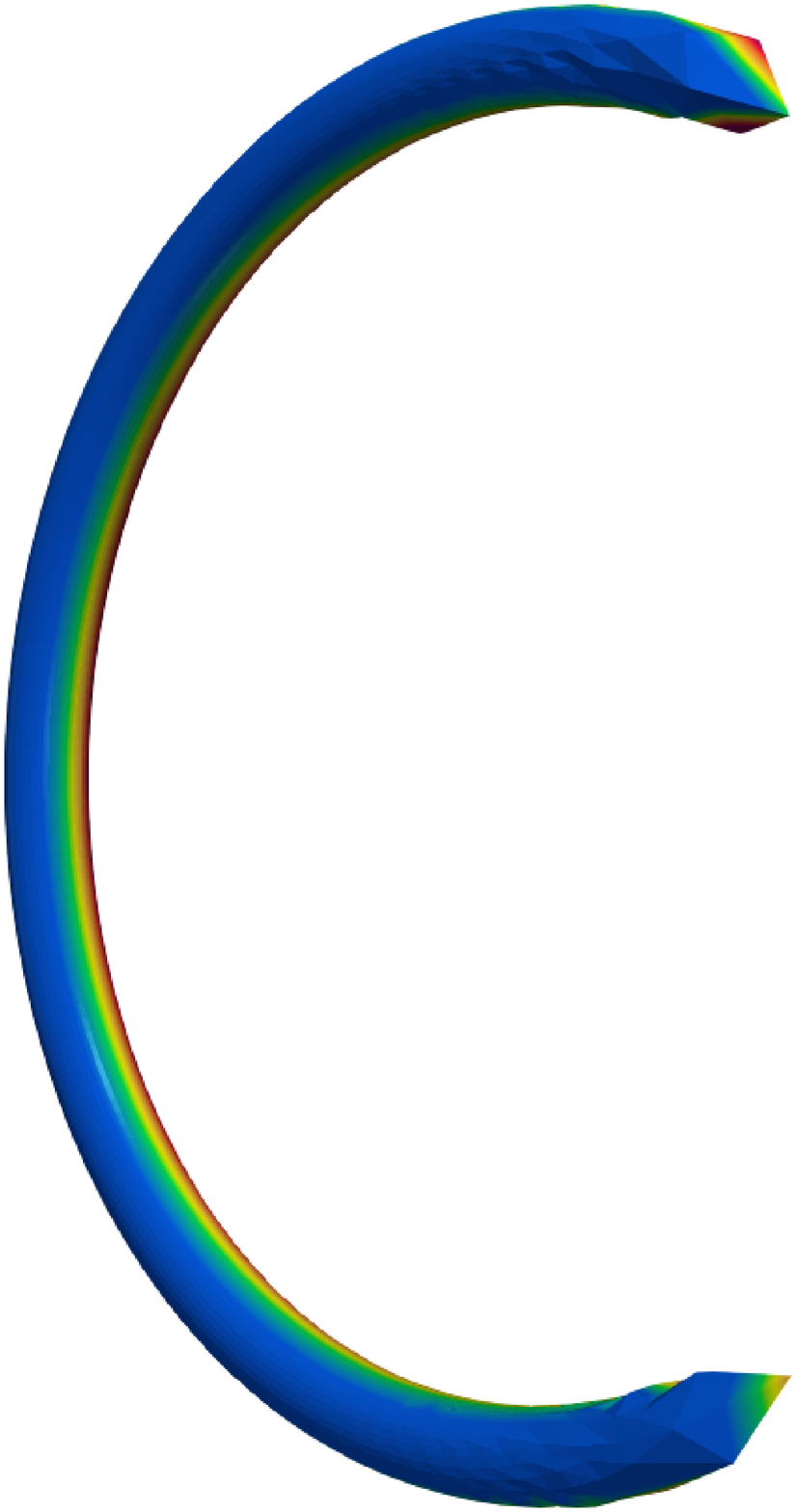}
    \end{minipage}
    \begin{minipage}[c]{0.15\textwidth}
    \includegraphics[width=\textwidth]{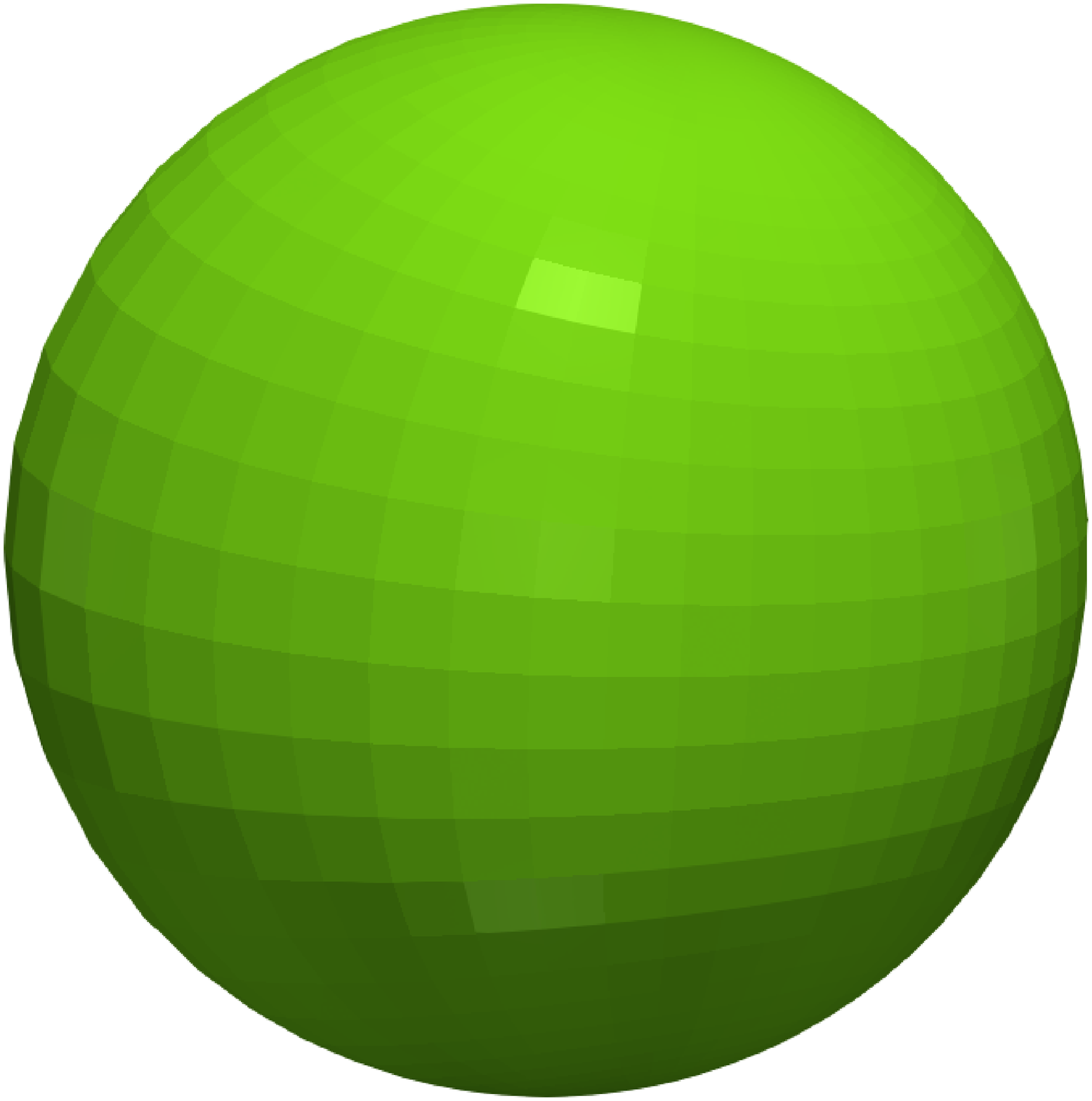}
    \end{minipage}
    \begin{minipage}[c]{0.15\textwidth}
    \includegraphics[width=\textwidth]{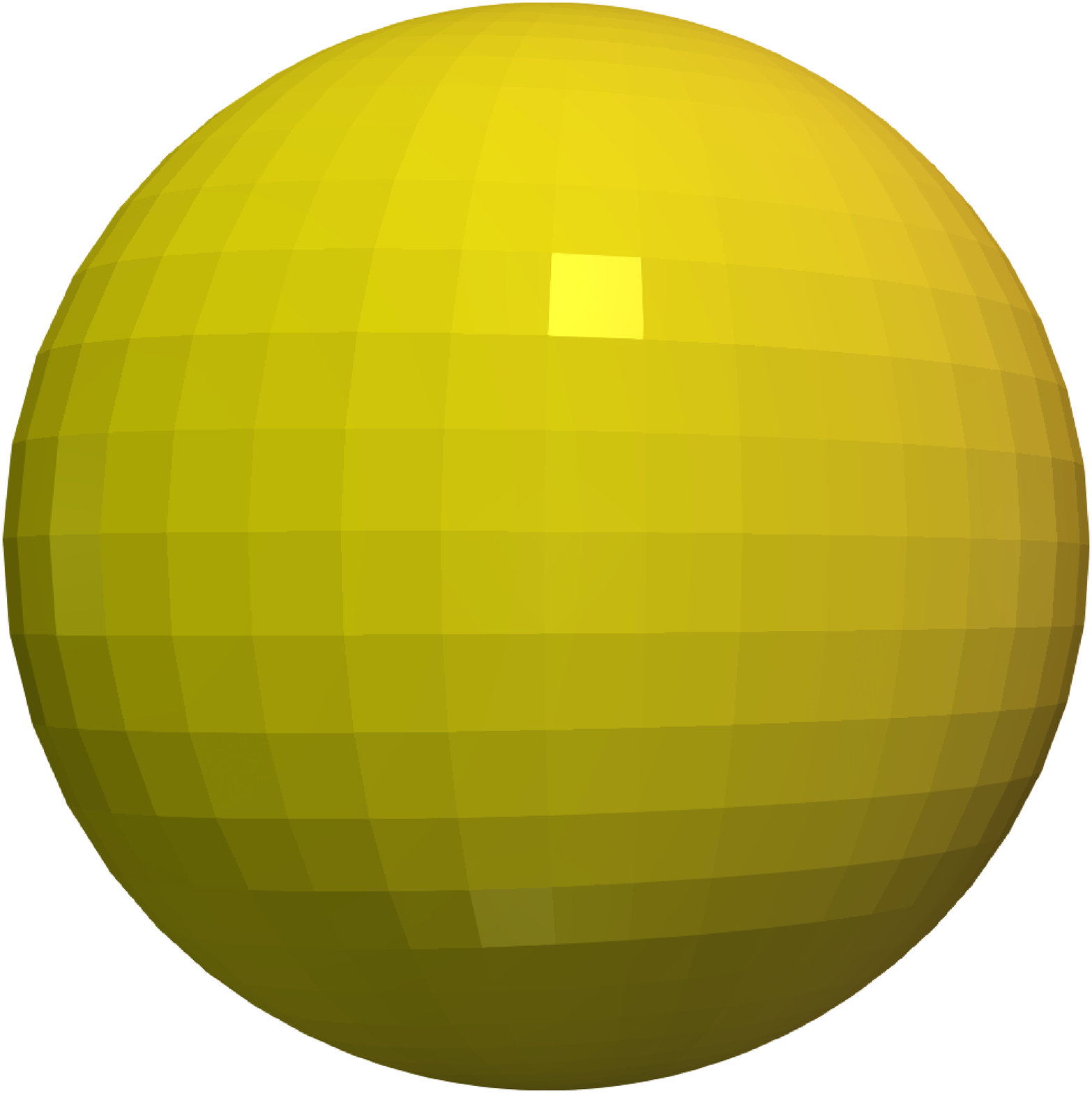}
    \end{minipage}
\end{center}
\caption[Comparison 0,2,4-Willmore flow]{p-Willmore evolution (with conformal penalty) of a letter ``C'' when $p=0,2,4$, respectively.  Colored by one component of mean curvature and oriented top to bottom.}
\label{fig:C}
\end{figure}

\green{It is well-known that the analytic properties of these flows are quite different from one another.  For example, convex surfaces evolving under MCF become extinct in finite time (see \cite{huisken1984}), while the Willmore flow can terminate in a round sphere of finite (positive) radius \cite{kuwert2001}.  In light of these differences, it is reasonable to wonder how the behavior of a geometric flow depends on the exponential weight of the mean curvature being measured, and the p-Willmore functional provides a natural way to investigate this idea.  In particular, it is apparent from simulation that when $p>2$, (at least some) surfaces which become spherical under the p-Willmore flow will instead grow indefinitely.  This is not surprising, as the p-Willmore functional is only invariant under changes of scale when $p=2$ (c.f. \cite{gruber2019thesis}).  Therefore, an immersed surface can easily decrease its p-Willmore energy by growing uniformly, so that its mean curvature decreases pointwise.  This phenomenon is displayed in Figure~\ref{fig:C}, where the p-Willmore evolution of a C-shaped surface is compared when $p=0,2,4$.  Moreover, Figure~\ref{fig:doggo} shows that even when the various p-Willmore flows terminate at a common immersion, their intermediate surfaces may be quite different depending on the value of $p$.}

\subsection{Contributions}
In the following sections, techniques from \cite{dziuk2008} will be adapted to express the $L^2$-gradient of $\mathcal{W}^p$ in a computationally-accessible way, resulting in an appropriate weak formulation of the p-Willmore flow problem.  Once the relevant system of PDE has been established, geometric constraints on surface area and enclosed volume will be considered and introduced into the flow model as Lagrange multipliers, leading to new and different behavior.  Moreover, the problem of mesh degradation along the flow will be discussed, and a minimization procedure will be given which dramatically improves mesh quality throughout the p-Willmore flow \green{at the expense of solving another nonlinear system at each time step.} This procedure is inspired by a conformality criterion of Kamberov, Pedit, and Pinkall derived in \cite{kamberov1996} and is similar in spirit to the least squares conformal mapping (LSCM) technique introduced in \cite{levy2002}.  Consequently, the p-Willmore flow and mesh regularization systems will be discretized and implemented on manifold meshes of triangles and quadrilaterals using the Finite Element Multiphysics Solver FEMuS \cite{aulisa2014}, and a fully-automated algorithm given for running the p-Willmore flow with \green{conformal penalty.  Finally, some specifics of this implementation will be discussed, as well as an application to mesh editing}.

\begin{figure}
\begin{center}
\begin{minipage}[c]{0.16\textwidth}
\includegraphics[width=\textwidth]{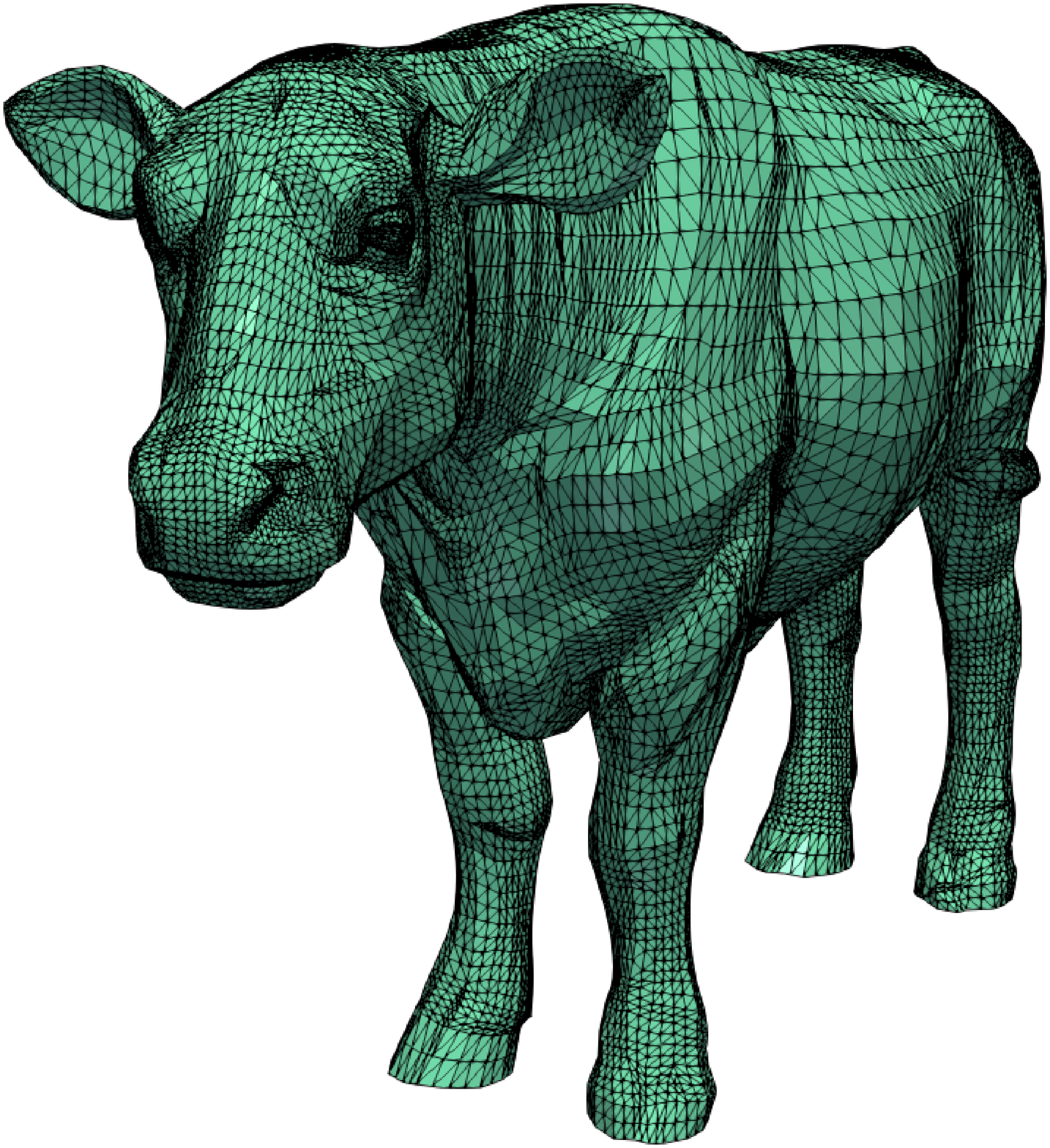}
\end{minipage}
\hspace{-0.5pc}
\begin{minipage}[c]{0.16\textwidth}
\includegraphics[width=\textwidth]{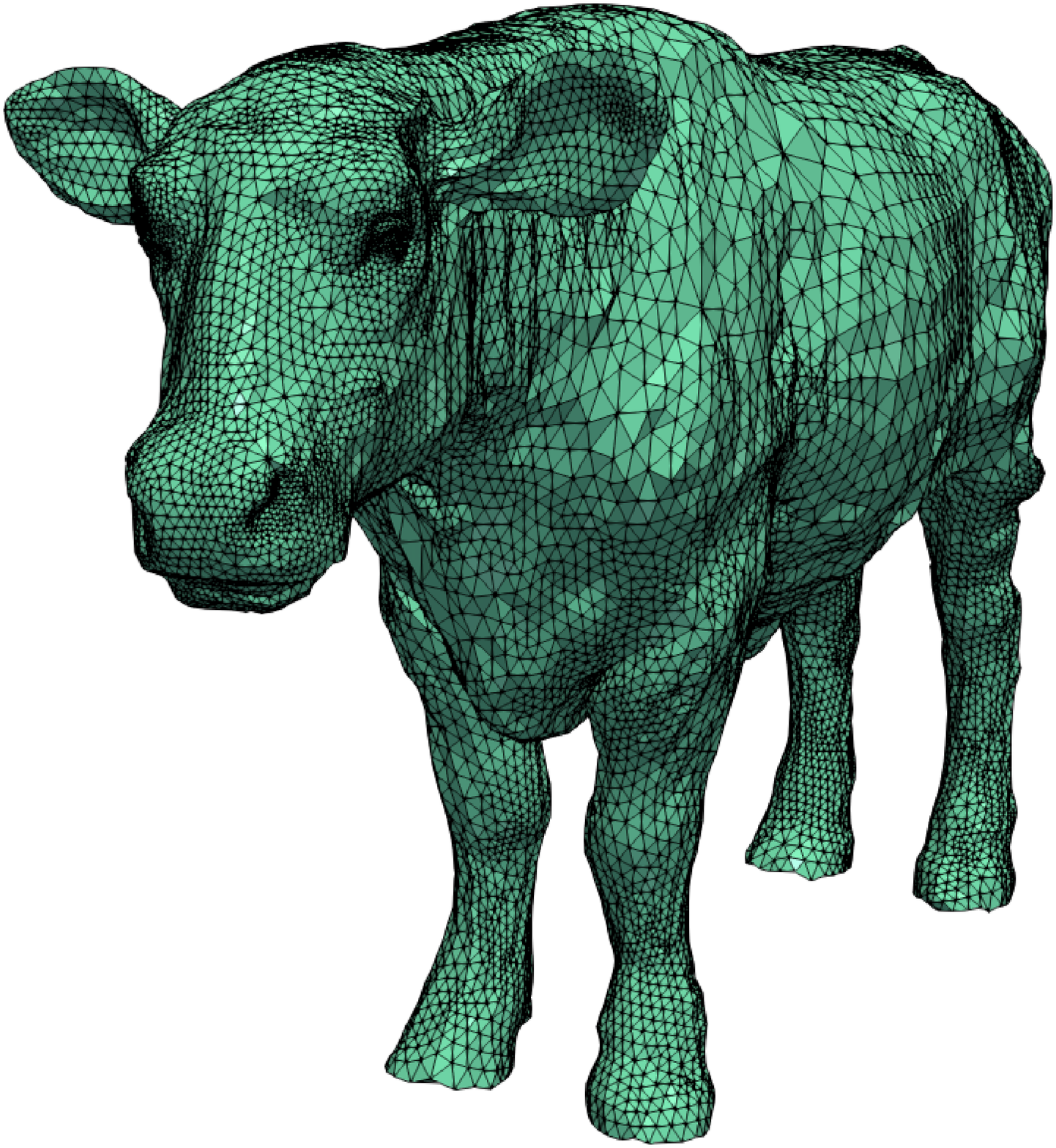}
\end{minipage}
\hspace{-0.4pc}
\begin{minipage}[c]{0.16\textwidth}
\includegraphics[width=\textwidth]{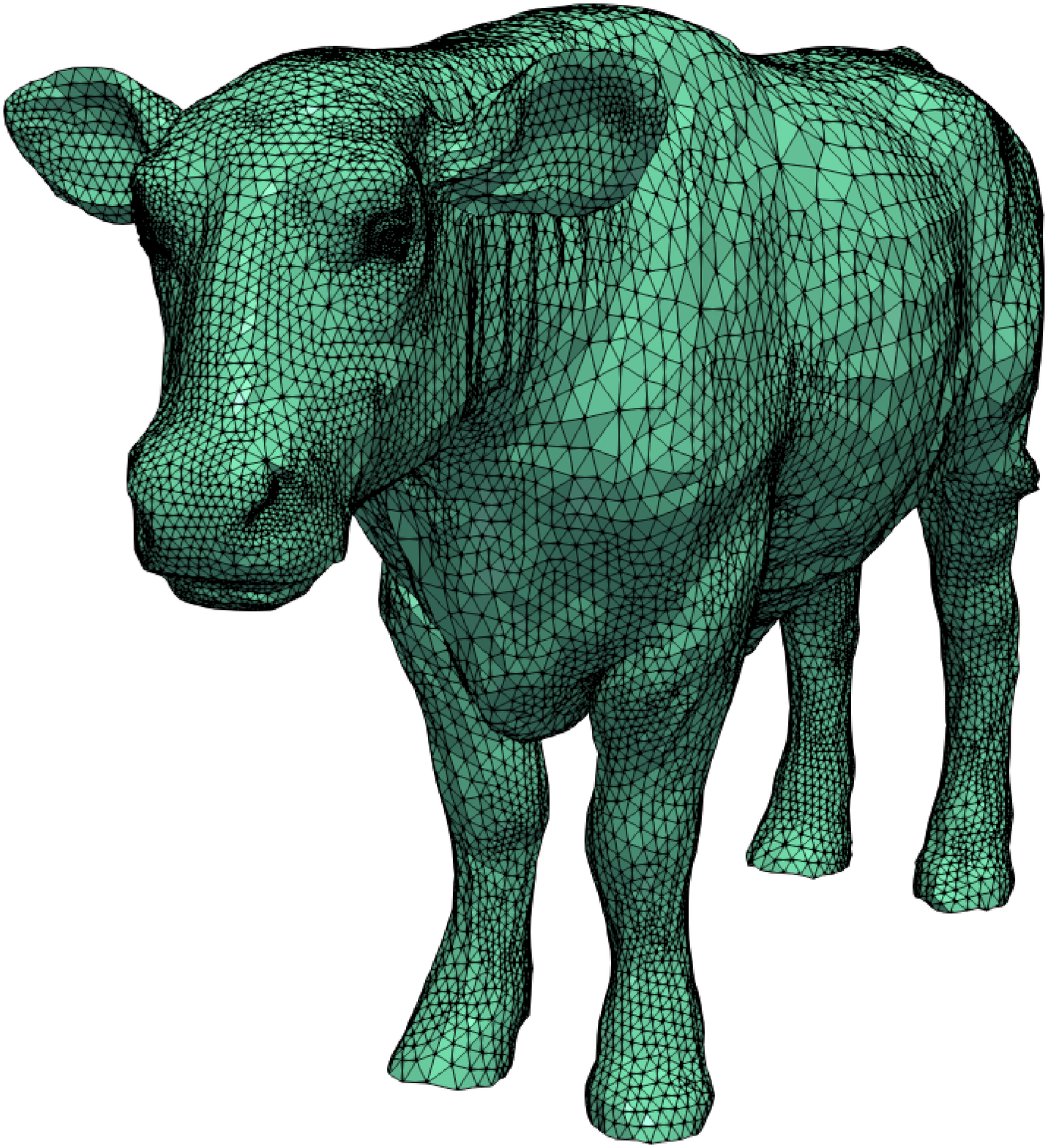}
\end{minipage}
\end{center}
\caption[Linear versus nonlinear regularization]{A performance comparison of linear versus nonlinear conformal penalty regularization on a cow with 34.5k triangles. Original mesh (left), linear algorithm (middle), nonlinear algorithm (right).  Simulations take roughly 2s, 4s, resp. on a 2.7 GHz Intel Core i5 with 8GB of RAM.}
\label{fig:moocompare}
\end{figure}

The p-Willmore flow algorithm introduced here has the following benefits:
\begin{itemize}
    \item \green{It provides a unified computational treatment of geometric flows which arise from functionals whose integrand is a power of the unsigned mean curvature, including MCF and the Willmore flow.}
    \item It is flexible with respect to geometric constraints on area and volume, as well as mesh geometry data (tri or quad) and surface genera.
    \item It affords the ability to near-conformally regularize the surface mesh along the flow, preventing mesh degeneration at the expense of \green{an additional nonlinear solve at each time step.}
    \item It is entirely minimization-based and therefore amenable to a large library of developed theory and techniques, including those in \cite{dziuk2013}.
\end{itemize}
\green{
\begin{remark}\label{rem:linvsnonlin}
The regularization procedure mentioned above can be easily modified to require only a linear solve, at the expense of more roughness in the mesh (c.f. Section 5).  See Figure~\ref{fig:moocompare} for a comparison on a realistic cow surface.  In addition, note that the conformal penalty regularization in this work is not a true constraint on the conformality class of the evolving surface. Therefore, the approach here differs from the work done on conformally-constrained Willmore surfaces in \cite{bohle2008,crane2011,schatzle2013} and others.
\end{remark}}

\green{
Though the p-Willmore flow with conformal penalty is useful, it is prudent to mention some challenges that have yet to be overcome.  In particular, the formulation considered here can be sensitive to initial data due to the high degree of nonlinearity present in the p-Willmore equation, especially when large values of $p$ are considered; typically, the flow cannot be run on rough meshes with a high degree of noise, and can be relatively unstable when $p>2$.  Moreover, the nonlinear systems involved in the p-Willmore flow algorithm are computationally demanding, requiring significant effort on fine meshes  which may be prohibitively expensive for ``real time'' use cases; specifics related to the figures in this work, including the solver time required are recorded in Table~\ref{tab:imp}.  Finally, the p-Willmore flow with conformal penalty is not yet well-understood with respect to theoretical results on consistency, stability, or convergence.  Such questions provide ample opportunity for future work in this area.
}

\section{Preliminaries}
It is beneficial to recall how to manipulate evolving surfaces mathematically.  Let \greener{ $M$ be a compact, connected $C^2$ surface without boundary. For $\varepsilon>0$, consider the family of surface immersions $u: M \times (-\varepsilon,\varepsilon) \to \mathbb{R}^3$ with images $M(t) \coloneqq u(M,t)$,} and let $\delta \coloneqq d/dt|_{t=0}$ be the variational derivative operator.  Then, if $\dot u$ denotes differentiation with respect to $t$, the initial surface $M(0)$ is said to undergo p-Willmore flow provided the equation
\begin{equation}\label{eq:pwillflow}
    \dot u = -\delta \mathcal{W}^p(u),
\end{equation}
is satisfied for all $t$ in some interval $(0,T]$. Using standard techniques from the calculus of variations, it can be shown (see \cite{gruber2019}) that for closed surfaces $M$ this condition implies the scalar equation
\begin{equation}\label{eq:gruberflow}
\green{
    \IPA{\dot u}{N} = -\frac{p}{2}\Delta_g \left( H|H|^{p-2}\right) - pH|H|^{p-2}\left(2H^2-K\right) + 2H|H|^{p}, }
\end{equation}
where $N: M \times (-\varepsilon,\varepsilon) \to S^2$ is the \green{outward-directed} unit normal vector to $M(t)$ for each $t$, $\Delta_g$ is the Laplace operator  associated to the metric $g$ on the surface, and $K$ is its Gauss curvature.  
\begin{remark}
Note that \green{from} here on the Einstein summation convention will be employed, so that any index appearing twice in an expression (once up and once down) will be implicity traced over.
\end{remark}
While equation (\ref{eq:pwillflow}) can be discretized by itself and used to define a normally-directed p-Willmore flow, it is advantageous to work directly with position instead of the mean curvature $H$.  Besides being more straightforward to implement, this allows for the consideration of tangential motion during the flow which can help regularize the surface mesh as it evolves (see \cite{dziuk2013}). 

\begin{remark}
Though position-based flow techniques are more standard in the literature, researchers in \cite{crane2013} have had success working directly with curvature.  Using a natural integrability condition, they are able to recover surface positions that maintain full conformality with respect to the reference immersion.  A major advantage of this approach is that such conformality is built directly into the flow, completely eradicating mesh degradation along the evolution.
\end{remark}
To develop a suitable model for the p-Willmore flow of surfaces, it is helpful to adopt the formalism of G. Dziuk found in \cite{dziuk2008}.  To that end, let $X: U \subset \mathbb{R}^2 \to M$ be a parametrization of (a portion of) the closed surface $M$, with \green{outward}-directed unit normal field $N$.  Then, the identity map $u: M \to \mathbb{R}^3$ defined through $u \circ X = X$ provides an isometric surface immersion, and the components of the induced metric on $M\subset \mathbb{R}^3$ are given by
\green{
\begin{equation*}
    g_{ij} = \IPA{\partial_i X}{\partial_j X} \coloneqq \IPA{X_i}{X_j}.
\end{equation*}
where $\langle\cdot,\cdot\rangle$ denotes the standard Euclidean inner product.
}
With this, the metric gradient $\nabla_g$ of a function $f:M \to \mathbb{R}^3$ can be expressed componentwise as \green{(Einstein summation assumed)}
\begin{equation*}
    \left(\nabla_g f \right) \circ X =\greener{ g^{ij}F_i \otimes X_j,}
\end{equation*}
where $F = f \circ X$ is the pullback of $f$ through the parametrization $X$, \greener{ $F_i = dF\left(X_i\right)$}, and $g^{ik}g_{kj} = \delta^i_j$.
It follows that the Laplace-Beltrami operator $\Delta_g$ on $M$ is then expressed as
\green{
\begin{equation*}
    \left(\Delta_g f\right) \circ X = \left(\divg\nabla_g f\right) \circ X = \frac{1}{\sqrt{\text{det}\,g}} \partial_j \left( \sqrt{\text{det}\,g}\, g^{ij} F_i\right),
\end{equation*}
and a simple calculation verifies that for two functions $f,h : M \to \mathbb{R}^3$, the metric inner product extends linearly to yield
\[\IP{df}{dh} = \IPA{\nabla_g f}{\nabla_g h} = \greener{ g^{ij}\langle F_i, H_j \rangle.} \]
}
Moreover, in view of the geometric identity \green{$Y \coloneqq \Delta_g u = 2HN$}, the p-Willmore functional is expressed succinctly in this framework as \green{
\begin{equation*}
    \mathcal{W}^p(u) = \frac{1}{2^p}\int_M |Y|^p\,d\mu_g.
\end{equation*}
In particular, introduction of the mean curvature vector $Y$ ensures that $\mathcal{W}^p$ is free of explicit \greener{second} derivatives of the position vector field.
}
\begin{remark}\label{rem:scaled}
Since the constant factor in front of the p-Willmore integrand merely scales the value of the functional and does not affect its geometric behavior, it will be omitted in subsequent passages with the understanding that $\mathcal{W}^p$ truly indicates $2^p\,\mathcal{W}^p$.  Note that this will manifest itself in the flow only as a uniform scaling of the temporal domain.
\end{remark}
\begin{figure}
\begin{center}
\begin{minipage}[c]{0.16\textwidth}
\includegraphics[width=\textwidth]{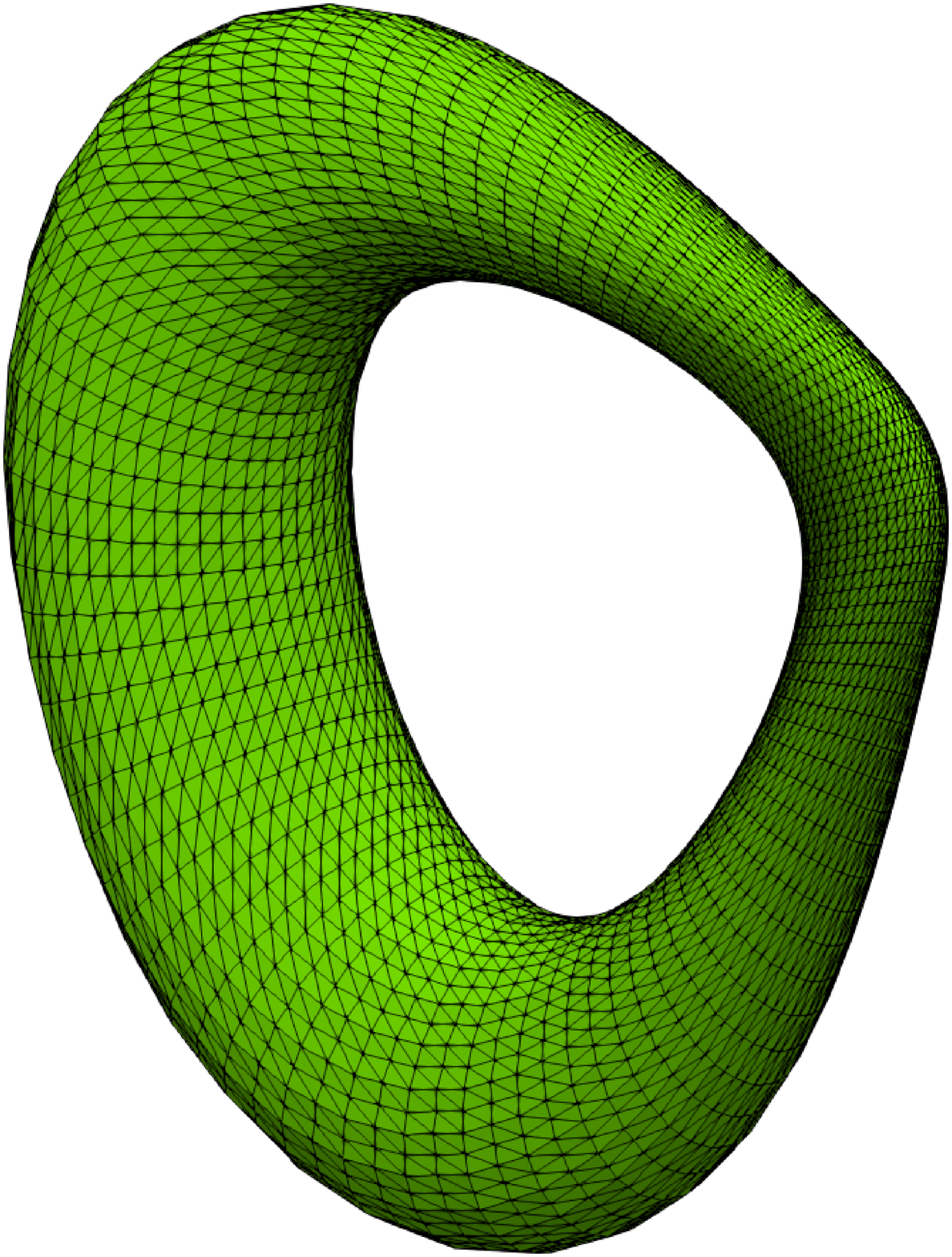}
\end{minipage}
\hspace{-0.6pc}
\begin{minipage}[c]{0.16\textwidth}
\includegraphics[width=\textwidth]{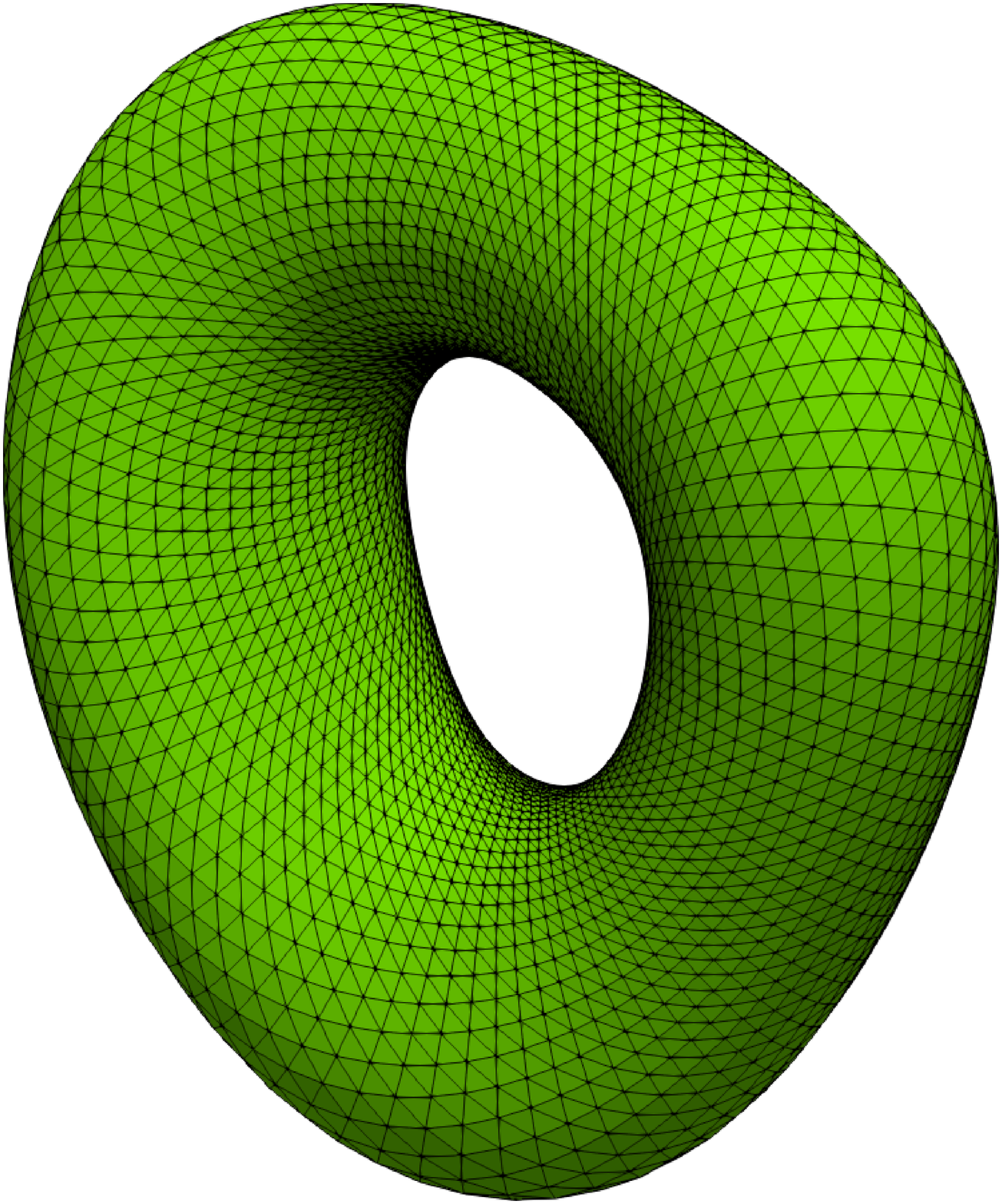}
\end{minipage}
\hspace{-0.3pc}
\begin{minipage}[c]{0.16\textwidth}
\includegraphics[width=\textwidth]{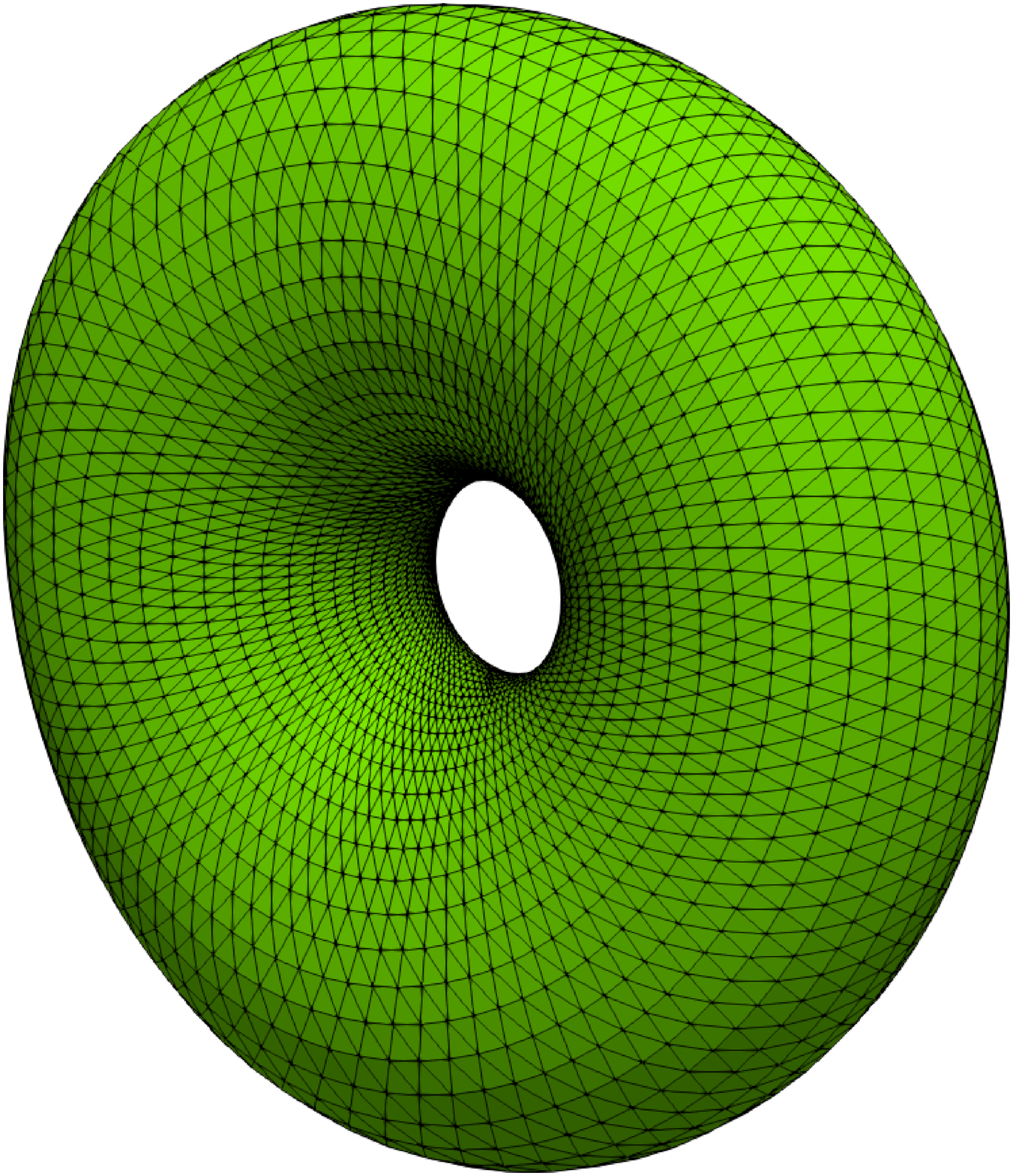}
\end{minipage}
\end{center}
\caption[2-Willmore evolution of a  torus]{\greener{Unconstrained and regularized} 2-Willmore evolution of a deformed torus to a known global minimum.  The minimizing surface is the stereographic projection of a Clifford torus in $S^3$.  \greener{Note that the conformal structure on the initial surface cannot be preserved as the pinched neck is thickened.}}
\label{fig:torus}
\end{figure}

\section{Building the p-Willmore flow model}
It is now possible to calculate the variational derivative ($L^2$-gradient) of the functional $\mathcal{W}^p$ in a way that is respectful towards computer implementation.  More precisely, the calculation presented here involves no adapted coordinate system or explicit second-order derivatives, and the variations considered are assumed to have tangential as well as normal components.  This will make it possible to accomplish the finite element discretization seen later.

Recall that when given a smooth function $\varphi: M \to \mathbb{R}^3$ and a parameter $t \in (-\epsilon, \epsilon)$, a variation of $u$ is given by
\begin{equation*}
    u(x,t) = u(x) + t\varphi(x),
\end{equation*}
where $x$ denotes a local coordinate on $M$.  This in turn induces a variation in the area functional, which can be calculated as in \cite{dziuk2008}.  In particular, there is the following lemma from that work.
\green{
\begin{lemma} \label{lem:dziuklem}
Let Greek letters indicate tensor components with respect to the standard basis for $\mathbb{R}^3$, and define $D(\varphi) = \nabla_g\varphi + (\nabla_g\varphi)^T$ through
\[D(\varphi)^{\alpha \beta} \coloneqq  g^{ij}\left(\varphi^\alpha_i X^\beta_j + \varphi_i^\beta X_j^\alpha \right).\] 
Then, in the notation above and denoting the area functional on $M$ by
\begin{equation*}
    \mathcal{A}(u) = \int_M 1\,d\mu_g,
\end{equation*}
the first and second variations of $\mathcal{A}$ may be expressed as
\begin{align*}
    \delta\mathcal{A}(u)\varphi &= \int_M \divg\varphi\,\,d\mu_g = \int_M \IP{du}{d\varphi}d\mu_g, \\
    \delta^2\mathcal{A}(u)(\varphi,\psi) &= -\int_M \left( \IP{D(\varphi)du}{d\psi} - \IP{d\varphi}{d\psi}\right)d\mu_g \\
    &+\int_M \left(\divg\varphi\right) \left(\divg\psi\right)d\mu_g.
\end{align*}
\end{lemma}
}


\begin{proof}
The proof is a direct calculation and can be found in \cite{dziuk2008}.
\end{proof}
With this in place, it is helpful also to recall an operator-splitting technique employed in \cite{dziuk2008}, which is used to reduce the order of the flow problem.  In particular, let $H^1(M;\mathbb{R}^3)$ denote the space of weakly first differentiable functions on $M$, and recall the equation $Y = \Delta_g u$.  Integrating this by parts against $\psi \in H^1(M;\mathbb{R}^3)$ then yields the relationship
\begin{equation}\label{eq:foru}\green{
    \int_M \left( \greener{ \langle Y,\psi\rangle } + \IP{du}{d\psi}\right) d\mu_g = 0, }
\end{equation}
which can be considered as a weak-form expression of the mean curvature vector $Y$. Note that due to the definition of $Y$, (\ref{eq:foru}) has the useful function of effectively reducing the order of the p-Willmore flow equation (\ref{eq:pwillflow}) by two at the expense of solving an additional PDE.
\begin{figure*}
\begin{center}
    \begin{minipage}[c]{0.22\textwidth}
    \includegraphics[width=\textwidth]{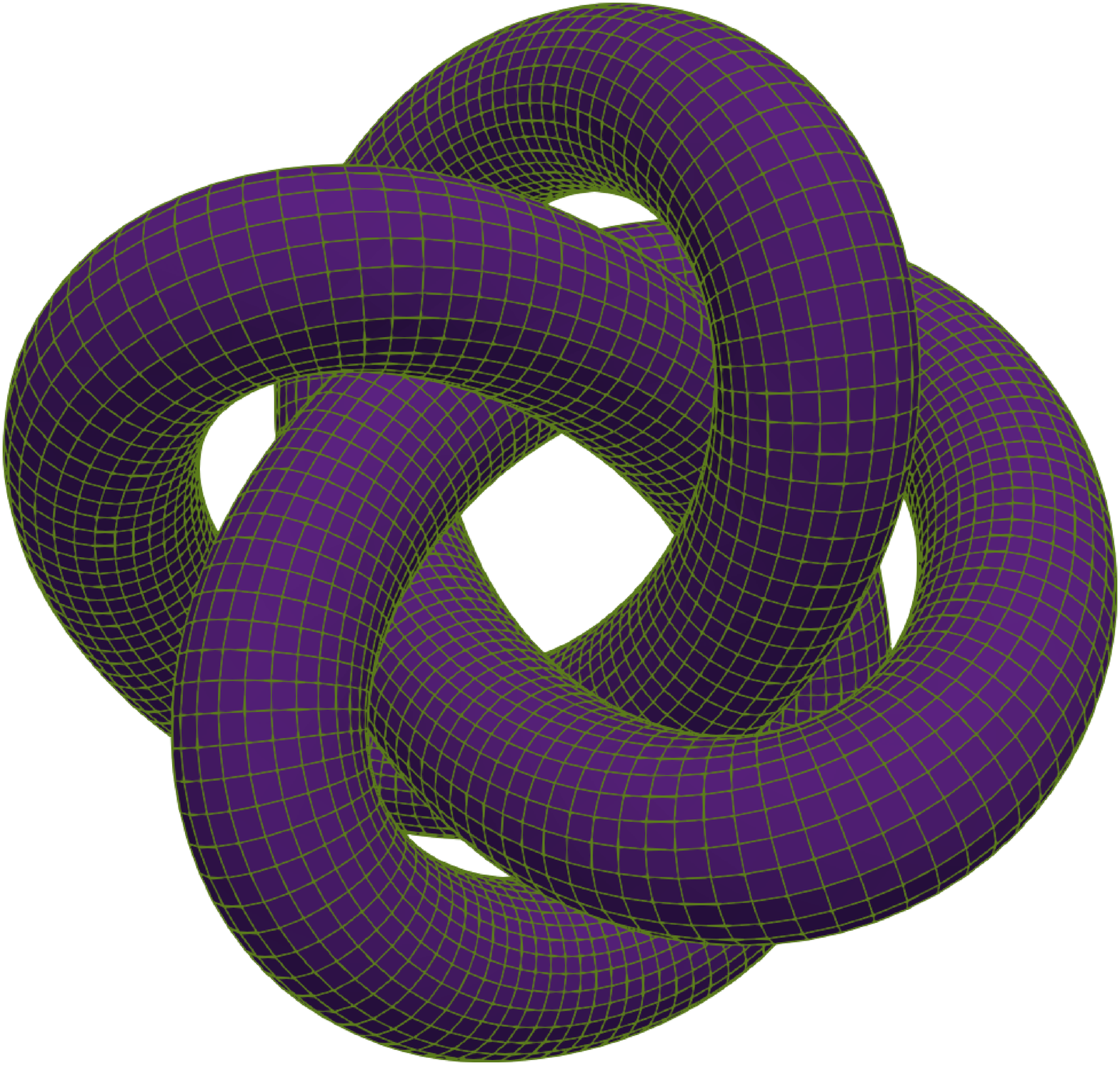}
    \end{minipage}
    \begin{minipage}[c]{0.22\textwidth}
    \includegraphics[width=\textwidth]{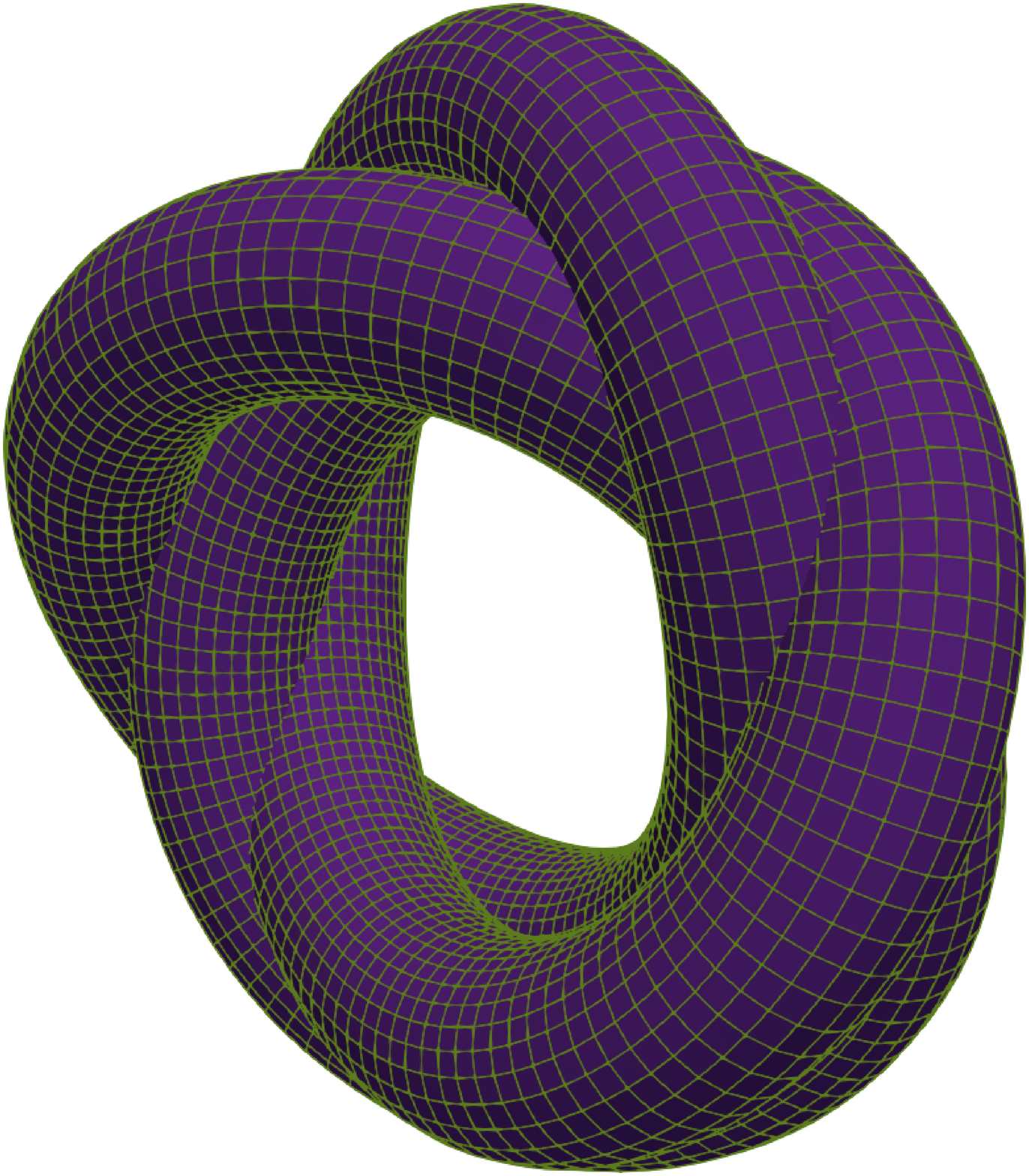}
    \end{minipage}
    \begin{minipage}[c]{0.22\textwidth}
    \includegraphics[width=\textwidth]{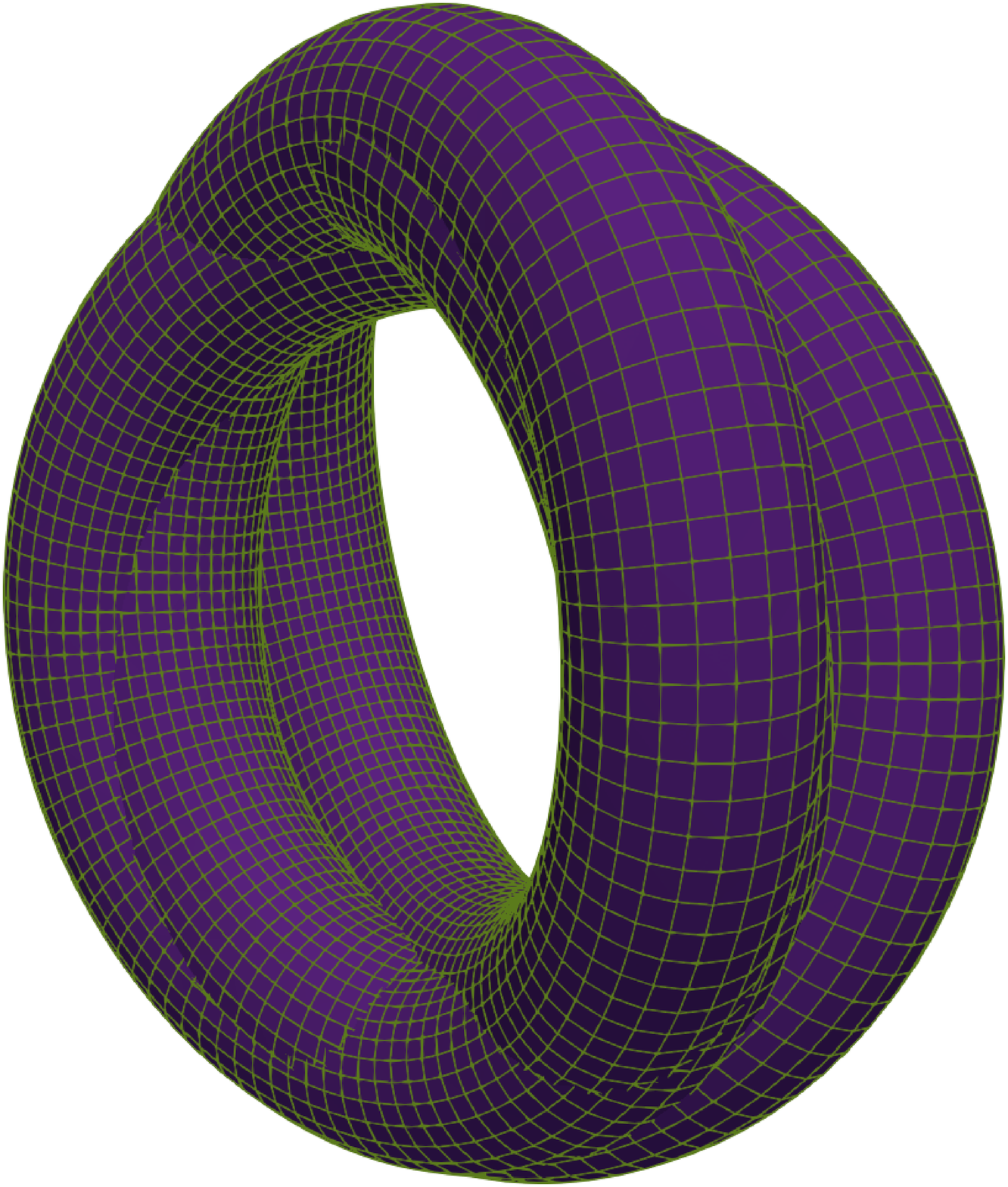}
    \end{minipage}
    \hspace{0.6pc}
    \begin{minipage}[c]{0.22\textwidth}
    \includegraphics[width=\textwidth]{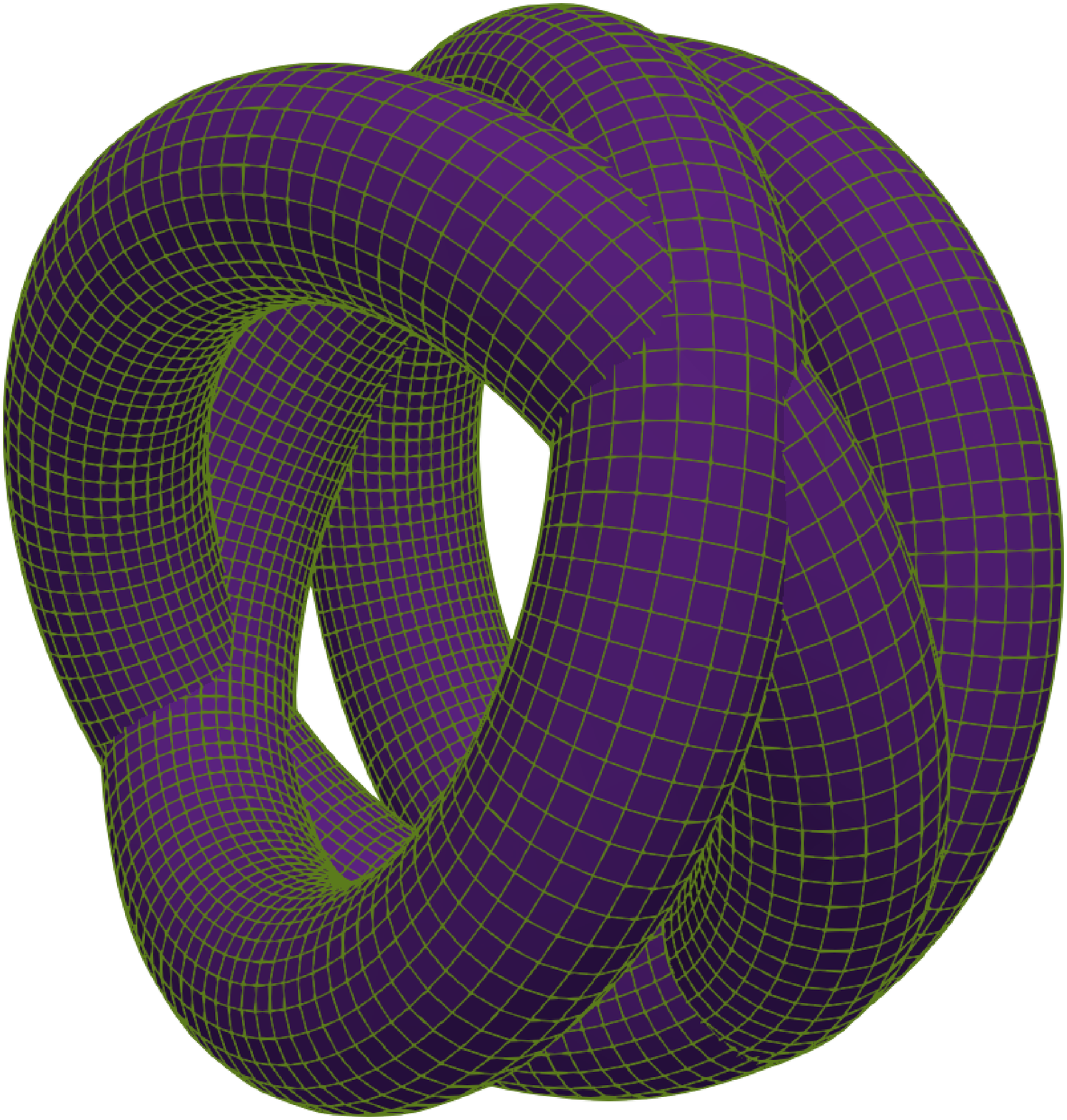}
    \end{minipage}
    \\
    \begin{minipage}[c]{0.22\textwidth}
    \includegraphics[width=\textwidth]{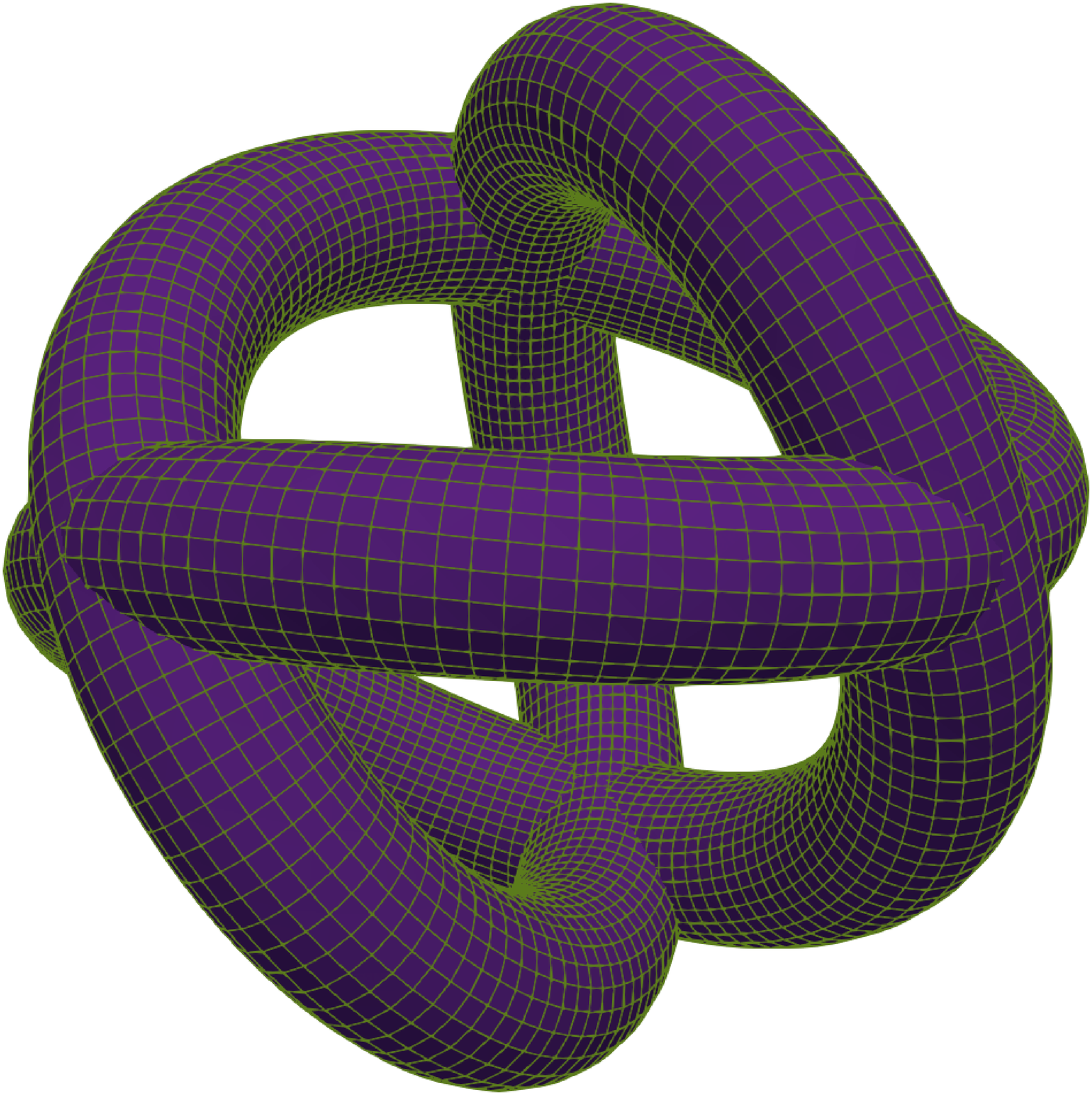}
    \end{minipage}
    \begin{minipage}[c]{0.22\textwidth}
    \includegraphics[width=\textwidth]{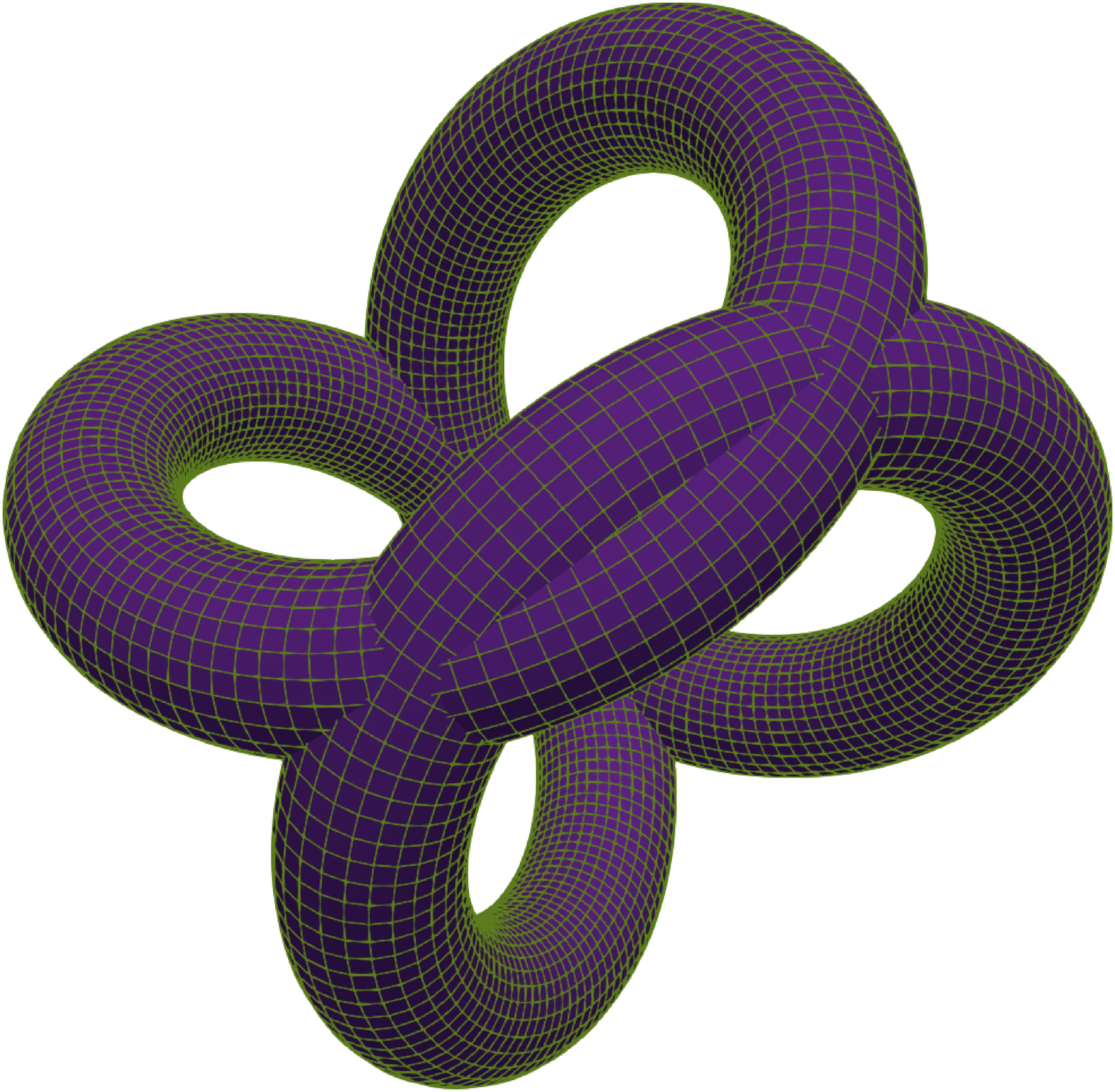}
    \end{minipage}
    \begin{minipage}[c]{0.22\textwidth}
    \includegraphics[width=\textwidth]{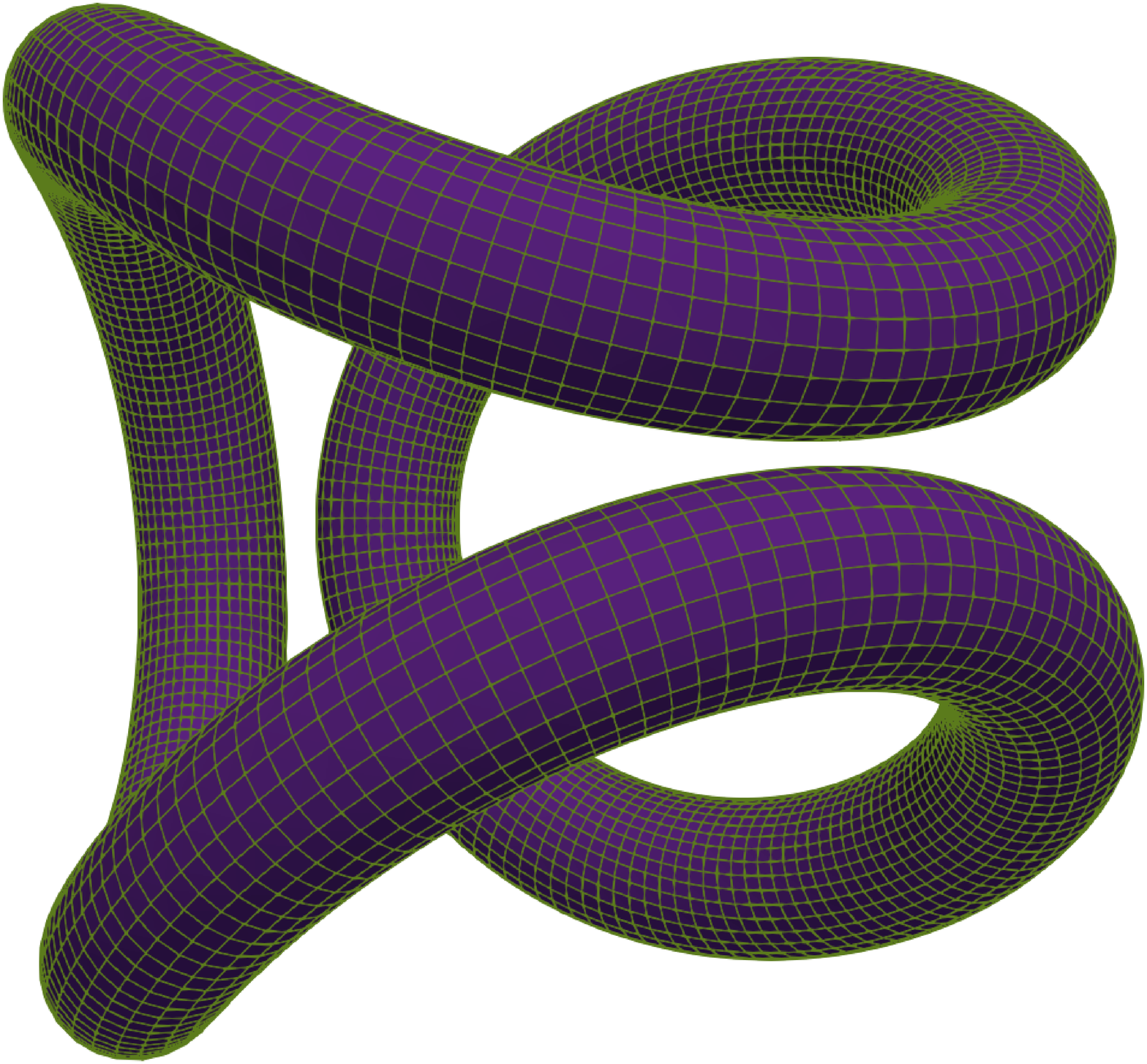}
    \end{minipage}
    \hspace{0.6pc}
    \begin{minipage}[c]{0.22\textwidth}
    \includegraphics[width=\textwidth]{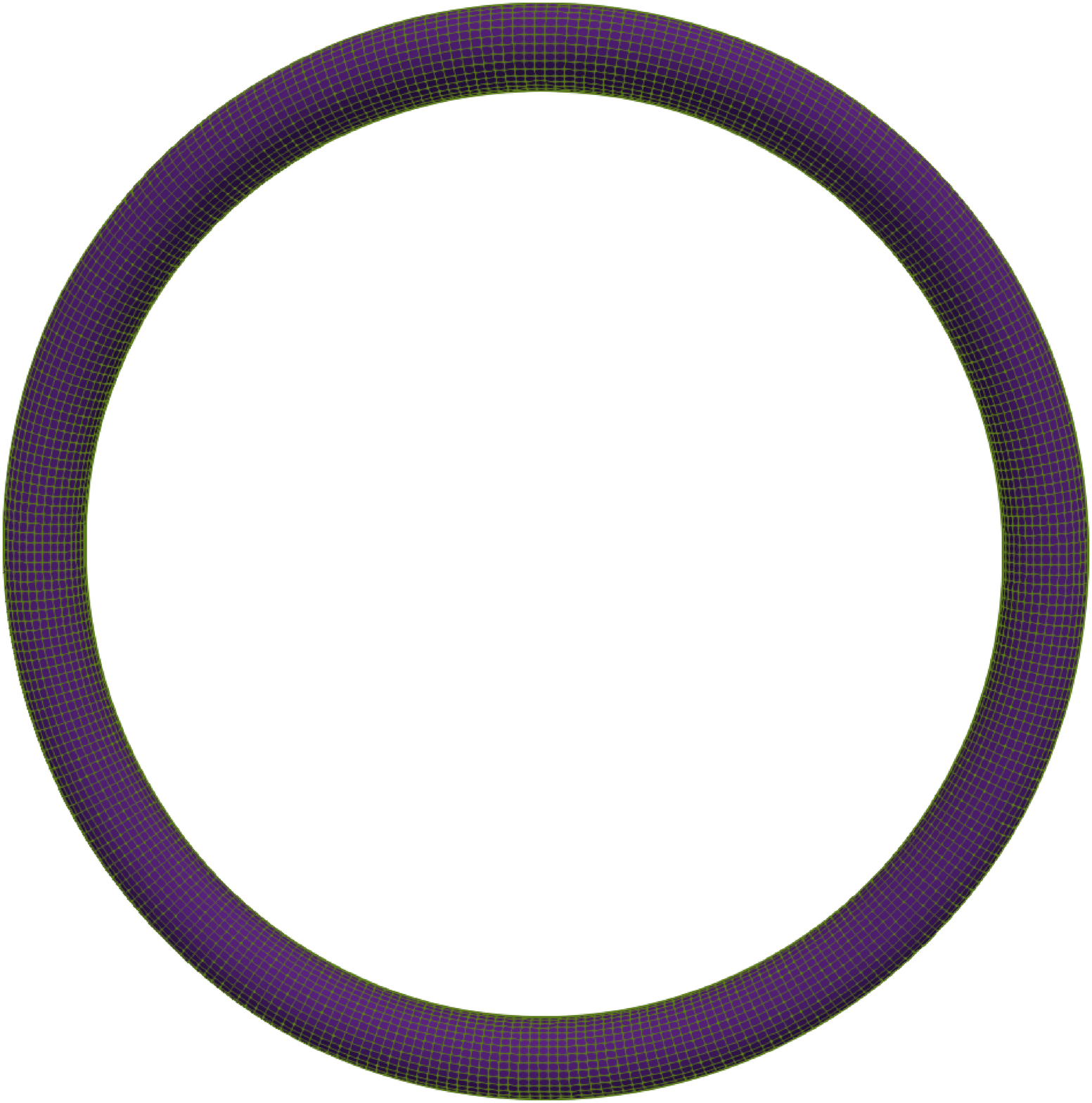}
    \end{minipage}
\end{center}
\caption[Almost-isometric flow of a (3,4)-torus knot]{Surface area and volume constrained 2-Willmore flow with conformal penalty applied to a (3,4)-torus knot.  Though area preservation is only globally enforced, the each surface along the evolution is nearly isometric to the given reference immersion.}
\label{fig:microsoftknot}
\end{figure*}

\begin{figure*}
\begin{center}
\begin{minipage}[c]{0.20\textwidth}
\includegraphics[width=\textwidth]{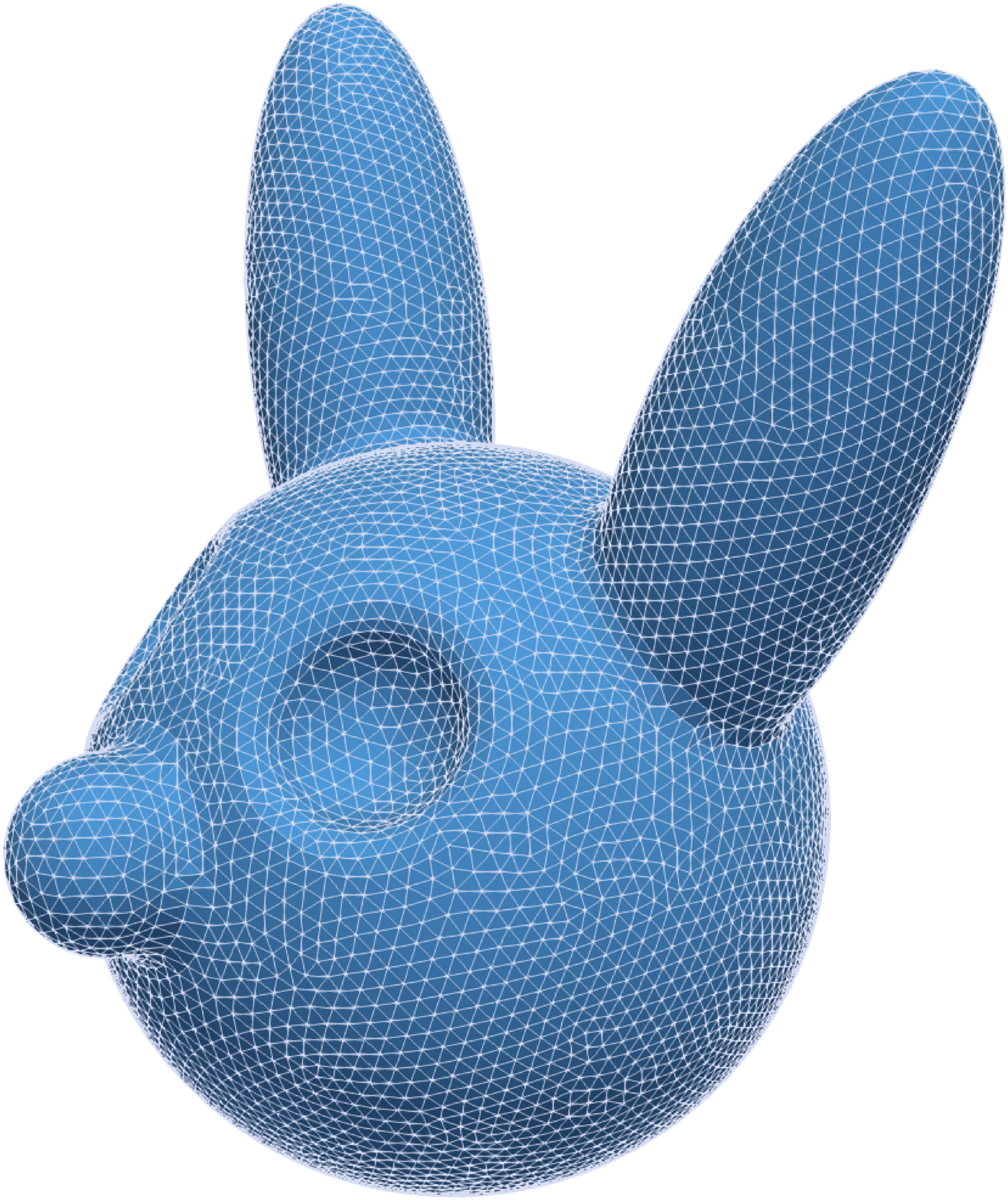}
\end{minipage}
\hspace{-0.6pc}
\begin{minipage}[c]{0.20\textwidth}
\includegraphics[width=\textwidth]{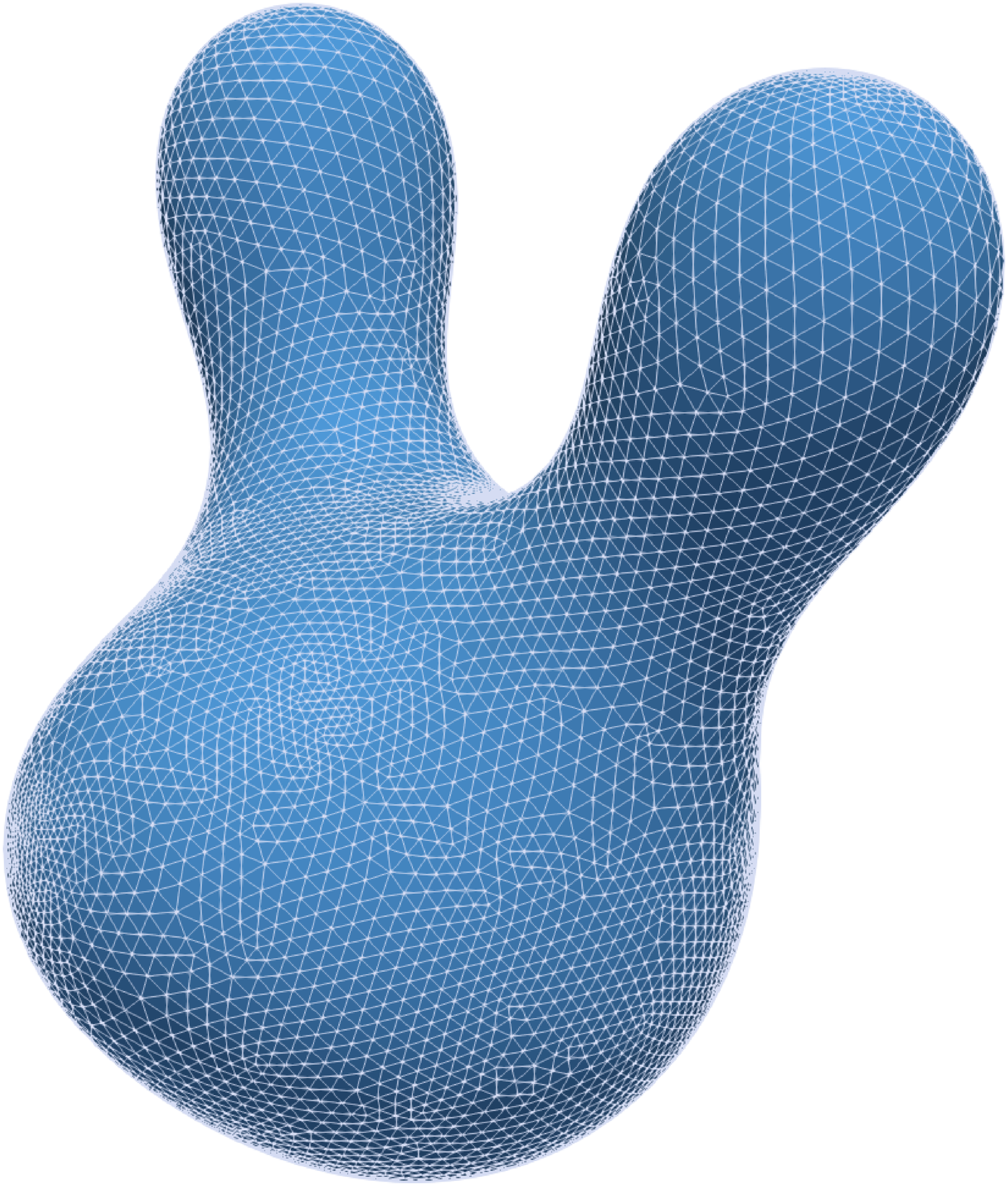}
\end{minipage}
\hspace{-0.5pc}
\begin{minipage}[c]{0.20\textwidth}
\includegraphics[width=\textwidth]{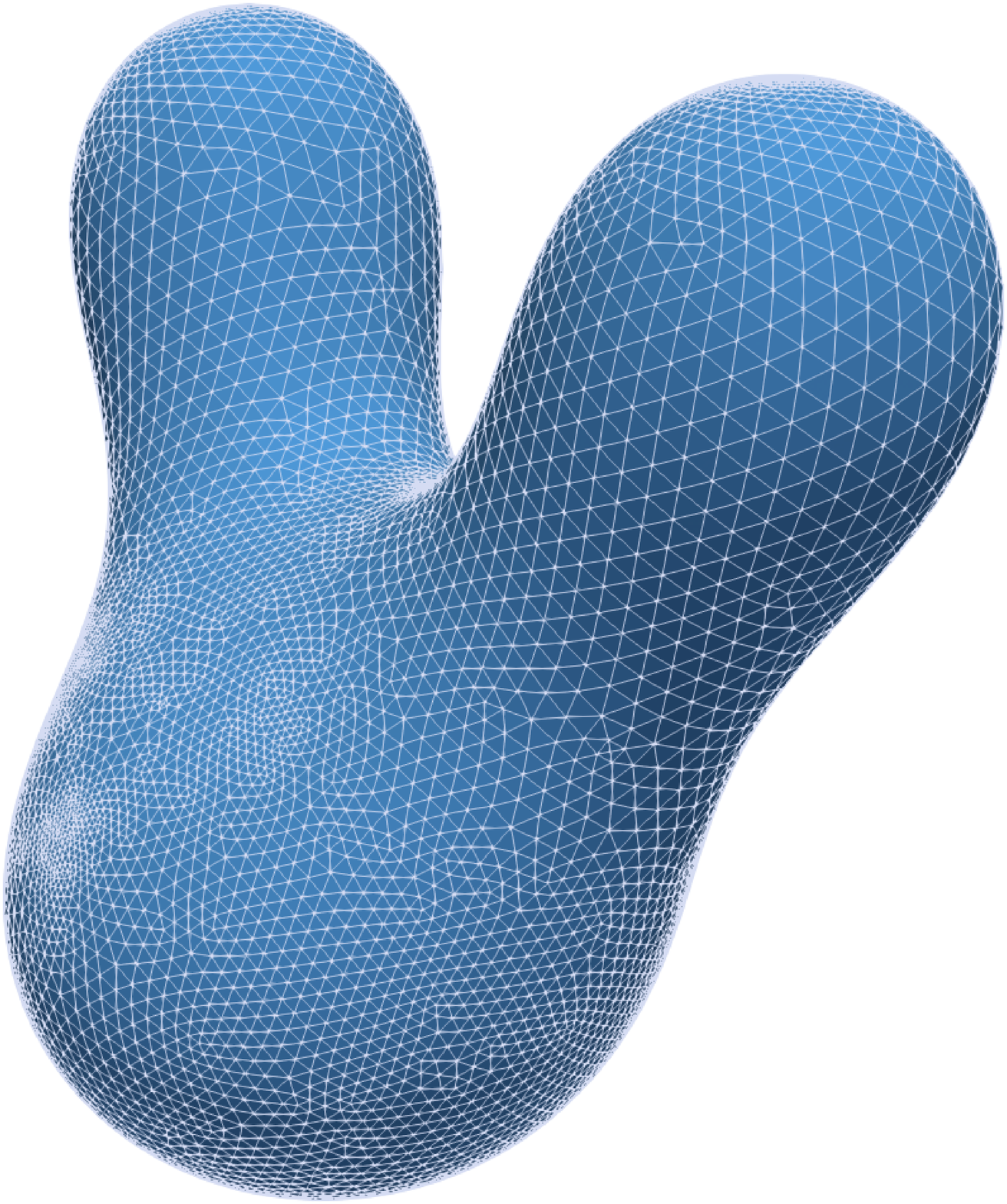}
\end{minipage}
\begin{minipage}[c]{0.20\textwidth}
\includegraphics[width=\textwidth]{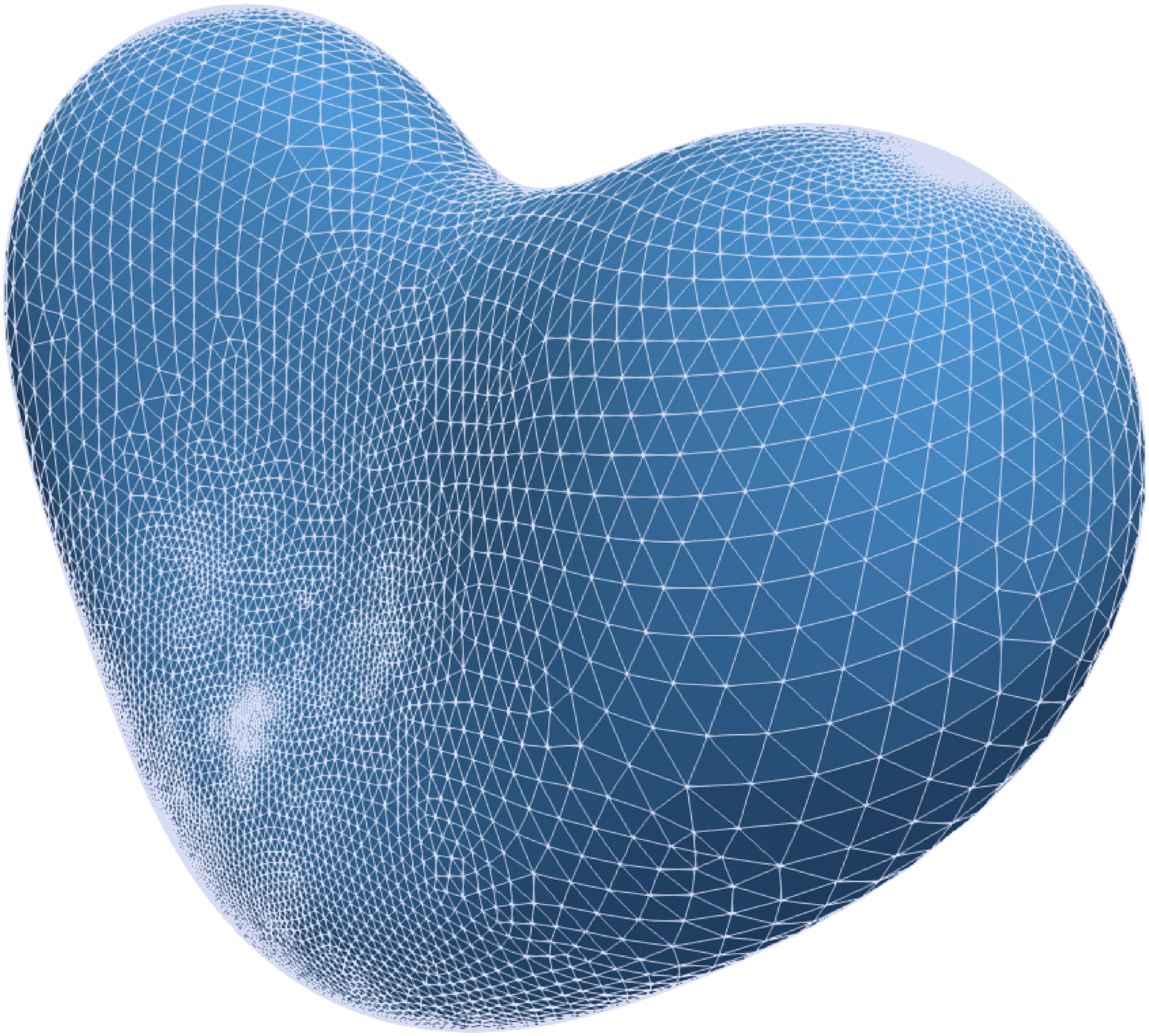}
\end{minipage}
\hspace{-0.5pc}
\begin{minipage}[c]{0.18\textwidth}
\includegraphics[width=\textwidth]{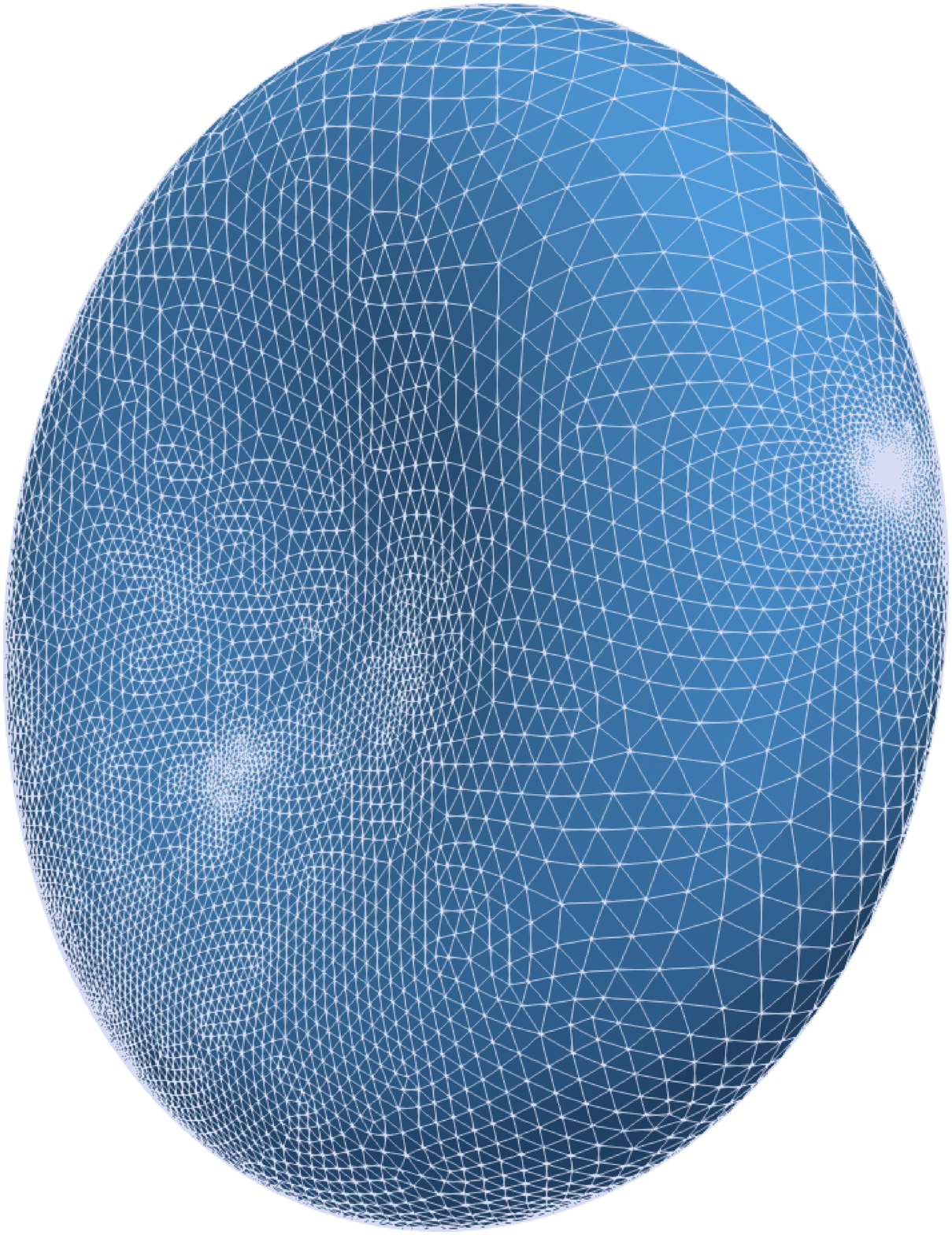}
\end{minipage}
\end{center}
\caption[2-Willmore evolution with constrained area and volume]{2-Willmore evolution of a rabbit with 18k triangles, constrained by both surface area and enclosed volume.  Here appears the biconcave discoid shape characteristic of genus 0 minimizers of the constrained \green{Helfrich-Canham} energy such as red blood cells \cite{ou2014}.  It is further remarkable that the flow behavior here is different than when either constraint is considered on its own, where the rabbit becomes a \greener{globally-minimizing} round sphere instead \greener{due to the scale-invariance of $\mathcal{W}^2$} (c.f. Figure~\ref{fig:doggo}).}
\label{fig:weirddog}
\end{figure*}

It is now pertinent to develop a counterpart to equation (\ref{eq:foru}), so that the operator splitting above can be beneficial.  The resulting equation should reduce to (\ref{eq:gruberflow}) in the normal direction, while also suppressing undesirable non-divergence terms such as $K$.  To accomplish this,
\greener{first note that 
\[ \int_M \IPA{Y}{\psi} d\mu_g = -\int_M \IP{du}{d\psi} = -\delta\mathcal{A}(u)\psi. \]}
Therefore, differentiating with respect to $u$ in the direction $\varphi\in H^1(M;\mathbb{R}^3)$ yields
\greener{
\begin{equation}\label{eq:dziukmagic}
    \int_M \IPA{\delta Y\,\varphi}{\psi} d\mu_g + \int_M \IPA{Y}{\psi}\divg\varphi\,d\mu_g = -\delta^2\mathcal{A}(u)(\varphi,\psi).
\end{equation}}
The goal is to use this expression to develop a weak-form for the p-Willmore equation by choosing an appropriate test function $\psi$.  Moreover, this choice should be made in avoidance of explicit derivatives of the normal vector $N$, since they are not well-suited to discretization using piecewise-linear finite elements. To this end, \greener{similar differentiation of the p-Willmore integrand yields,
\begin{equation*}\label{eq:deltanormYp}
\begin{split}
    \delta |Y|^p\,\varphi &=  \delta \IPA{Y}{Y}^{\frac{p}{2}}\varphi = p|Y|^{p-2}\IPA{\delta Y\,\varphi}{Y} \\
    &= \IPA{\delta Y\,\varphi}{\,p|Y|^{p-2} Y}.
\end{split}
\end{equation*}}
Hence, letting $W:= |Y|^{p-2}Y$ be the weighted mean curvature vector on $M$, choosing $\psi = pW$ in equation (\ref{eq:dziukmagic}), \greener{and using Lemma~\ref{lem:dziuklem}} the variation of the ($2^p$-scaled) Willmore functional $\mathcal{W}^p$ can be calculated as 
\begin{equation*}
    \begin{split}
        \delta\mathcal{W}^p(u)\varphi &= \delta \int_M |Y|^p\, d\mu_g \\
        &= \int_M \IPA{\greener{\delta Y\,\varphi}}{\,p|Y|^{p-2} Y} d\mu_g + \int_M |Y|^p\,\divg\varphi\,\, d\mu_g \\
        &= \int_M \left((1-p)|Y|^p - p\,\divg W\right)\divg\varphi\,\,d\mu_g \\
        &\qquad+ \int_M p\left(\IP{D(\varphi)du}{dW} - \IP{d\varphi}{dW}\right) d\mu_g.
    \end{split}
\end{equation*}
This computation directly implies the following Theorem.
\green{
\begin{theorem}\label{thm:basicpwillflow}
    In the notation above and for $p\geq 1$, the unconstrained p-Willmore flow equation (\ref{eq:pwillflow}) is expressed in weak form by the following system of PDE in the variables u, Y, and W:
    \begin{align}
    0 &= \int_M \IPA{\dot{u}}{\varphi}d\mu_g + \int_M \left((1-p)|Y|^p - p\,\divg W\right)\divg\varphi\,\,d\mu_g \notag \\
    &+ \int_M p \left( \IP{D(\varphi)du}{dW} - \IP{d\varphi}{dW} \right)d\mu_g, \label{eq:forW} \\
    0 &= \int_M \IPA{Y}{\psi}d\mu_g + \int_M \IP{du}{d\psi}d\mu_g, \label{eq:forus} \\
    0 &= \int_M \IPA{W - |Y|^{p-2}Y}{\xi}  d\mu_g = 0, \label{eq:forY}
\end{align}
which must hold for all $t\in (0,T]$ and all $\varphi, \psi, \xi \in H^1\left(M(t);\mathbb{R}^3\right)$.
\end{theorem}
}
\begin{proof}
The proof follows immediately from the definitions of $Y,W,$ and the discussion above.
\end{proof}
\begin{remark}
The reader will notice that this reduces to precisely the system in \cite{dziuk2008} for the case $p=2$, in which case the last equation is not needed.  Additionally, the case $p=0$ (MCF), while not in the domain of the theorem as stated, may be recovered by simply omitting (\ref{eq:forY}) and replacing equation (\ref{eq:forW}) with the equation from Lemma~\ref{lem:dziuklem} for the variation of area:
\begin{equation*}\green{
    \int_M \IPA{\dot u}{\varphi}d\mu_g - \int_M \IP{du}{d\varphi}d\mu_g = 0. }
\end{equation*}
\end{remark}
The system in Theorem~\ref{thm:basicpwillflow} is the primary model for the p-Willmore flow studied here, and provides the basis for the p-Willmore flow algorithm presented later.  Before discussing further modifications, the following theoretical result is presented which guarantees that the p-Willmore energy always decreases along the flow governed by the equations above. \green{Note that this property is well known in the case of MCF (0-Willmore flow), so the proof of this case is omitted.  See e.g. \cite{mantegazza2011} for more details.} Example illustrations of this phenomenon include Figures~\ref{fig:funcow}, \ref{fig:torus}, \ref{fig:microsoftknot}, \ref{fig:weirddog}, and \ref{fig:statueflow}.
\green{
\begin{theorem}\label{thm:paramflow}
The closed surface p-Willmore flow is energy decreasing for $p\geq 1$.  That is, if $W = |Y|^{p-2}Y$ is the weighted mean curvature vector on $M$ and $u: M \times (0,T] \to \mathbb{R}^3$ is family of surface immersions with $M(t) = u(U, t)$ satisfying the weak p-Willmore flow equations (\ref{eq:forW}), (\ref{eq:forus}), and (\ref{eq:forY}),
then the $p$-Willmore flow satisfies
\begin{equation*}
    \int_{M(t)} |\dot{u}|^2\,d\mu_g + \frac{d}{dt}\int_{M(t)}|Y|^p\,d\mu_g = 0.
\end{equation*}
\end{theorem}}
\green{\begin{proof}
Choosing the admissible test functions $\varphi = \dot{u}$ and $\psi = pW$ in (\ref{eq:dziukmagic}), as well as noticing that $\IPA{W}{Y} = |Y|^p$, the following system is observed
\begin{align*}
       0 &= \int_{M} |\dot{u}|^2\,d\mu_g +\int_M \left( (1-p)|Y|^p\divg\dot{u} - p\left(\divg W\right)\left(\divg\dot{u}\right)\right)d\mu_g \notag\\
       &+ \int_M p\left( \IP{D(\dot{u})du}{dW} - \IP{dW}{d\dot{u}}\right)d\mu_g,\\
       0 &= \int_M p\IPA{\delta Y}{W}d\mu_g + \int_M p\left(|Y|^p + \divg W\right)\divg\dot{u}\,\,d\mu_g \notag \\
       &+ \int_M p\left(\IP{dW}{d\dot{u}} - \IP{D(\dot{u})du}{dW}\right)d\mu_g.
\end{align*}
Adding the above equations in view of (\ref{eq:deltanormYp}) then yields
\begin{equation*}
\begin{split}
    0 &= \int_M |\dot{u}|^2\,d\mu_g + \int_M \left( \delta|Y|^p + |Y|^p \divg\dot{u}\right)d\mu_g \\
    &= \int_M |\dot{u}|^2\,d\mu_g + \frac{d}{dt}\int_M |Y|^p\, d\mu_g,
\end{split}
\end{equation*}
completing the argument.
\end{proof}
}
Now, in light of the physical relevance of the functional $\mathcal{W}^p$, it is desirable also to have a model for the p-Willmore flow that is amenable to geometric constraints on surface area and enclosed volume.  This is reasonable not only from a physical point of view (since many curvature-minimizing structures such as biomembranes constrain themselves naturally in these ways) but also in a mathematical sense, as such constraints can serve as a meaningful ``replacement'' for conformal invariance when $p\neq 2$.  More precisely, since the p-Willmore functional is not conformally invariant in general, volume/area preservation ensures that physically-meaningful shapes such as spheres remain locally minimizing for $p\neq 2$, at least among some class of variations.  Practically, this is accomplished through the addition of Lagrange multipliers $\lambda, \mu$ into the model from Theorem~\ref{thm:basicpwillflow}.  More precisely, let $D\subset \mathbb{R}^3$ be a region in space such that $\partial D = M$ \green{and let $d\mu$, $\mathrm{div}$ denote, respectively, the volume element and divergence operator on $\mathbb{R}^3$.  Recall the volume functional, }
\begin{equation*}\green{
    \mathcal{V}(u) = \int_D 1\,d\mu = \frac{1}{3}\int_D \mathrm{div}\,u\,\,d\mu = \frac{1}{3} \int_M \IPA{u}{N}d\mu_g, }
\end{equation*}
where the Divergence Theorem was applied in the last equality.  It is well known (see e.g. \cite{barbosa2012}) that the first variation of volume is given by
\begin{equation*}\green{
    \delta \mathcal{V}(u)\varphi = \frac{1}{3}\int_M \IPA{\varphi}{N}d\mu_g. }
\end{equation*}
On the other hand, recall that Lemma~\ref{lem:dziuklem} implies that the first variation of the area functional can be expressed as
\begin{equation*}\green{
    \delta\mathcal{A}(u)\varphi = \delta \int_M 1\,d\mu_g =  \int_M \IP{du}{d\varphi} d\mu_g. }
\end{equation*}
With these expressions available, it is straightforward to formulate the next problem considered in this work: closed surface p-Willmore flow with constraint.
\begin{problem}[Closed surface p-Willmore flow with constraint]\label{prob:constrainedpwillmore}
Let $p\geq 1$ and \green{ $W := |Y|^{p-2}Y$.} Determine a family $u: M\times (0,T] \to \mathbb{R}^3$ of surface immersions with $M(t)\coloneqq u(M,t)$ such that $M(0)$ has initial volume $V_0$, initial surface area $A_0$, and for all $t \in (0,T]$ the equation
\begin{equation*}\green{
    \dot u = -\delta \, (\mathcal{W}^p + \lambda \mathcal{V} + \greener{\gamma} \mathcal{A}), }
\end{equation*}
is satisfied for some piecewise-constant functions $\lambda, \mu: M(t) \to \mathbb{R}$.  Stated in weak form, the goal is to find functions $u, Y, W, \lambda, \greener{\gamma}$ on $M(t)$ such that the equations
\green{
\begin{align}
    0 &= \int_M \IPA{\dot{u}}{\varphi}d\mu_g + \int_M \greener{\gamma} \IP{du}{d\varphi}d\mu_g + \int_M \lambda\IPA{\varphi}{N}d\mu_g  \notag \\
    &+ \int_M \left((1-p)|Y|^p - p\,\,\divg W\right)\divg\varphi\,\,d\mu_g \notag \\
    &+\int_M p \left( \IP{ D(\varphi)du}{dW} - \IP{d\varphi}{dW}\right)d\mu_g, \label{eq:forWconstrained} \\
    0 &= \int_M \IPA{Y}{\psi}d\mu_g + \int_M \IP{du}{d\psi}d\mu_g, \notag \\
    0 &= \int_M \IPA{W - |Y|^{p-2}Y}{\xi}d\mu_g,\notag \\
    A_0 &= \int_M 1\,d\mu_g, \label{eq:forarea} \\
    3V_0 &= \int_M \IPA{u}{N}d\mu_g, \label{eq:forvol}
\end{align}
are satisfied for all $t\in (0,T]$ and all $\varphi, \psi, \xi \in H^1\left(M(t); \mathbb{R}^3\right)$.
}
\end{problem}
\begin{remark}
The case where $p=0$ may again be considered by replacing the equation (\ref{eq:forWconstrained}) with the simpler relationship
\begin{equation*}\green{
        0= \int_M \IPA{\dot u}{\varphi}d\mu_g + \int_M \lambda\IPA{\varphi}{N}d\mu_g - \int_M \IP{du}{d\varphi}d\mu_g. }
\end{equation*}
Of course, area preservation makes no sense in this context (since the objective of MCF is to decrease area), so equation (\ref{eq:forarea}) should also be disregarded in this case.  In addition,
note that the system of Problem~\ref{prob:constrainedpwillmore} can also be used to study the p-Willmore flow with fixed volume or fixed surface area separately.  In particular, fixed volume is obtained by setting $\mu\equiv 0$ and ignoring (\ref{eq:forarea}), and fixed area is accomplished similarly with $\lambda\equiv 0$ and omission of (\ref{eq:forvol}).  In practice, Boolean variables were implemented to enable switching between the different constrained/unconstrained cases.  
\end{remark}
Problem~\ref{prob:constrainedpwillmore} provides a way to examine the p-Willmore flow subject to geometric constraints on surface area or enclosed volume.  This is a highly interesting situation, since minimizing surfaces can vary widely with the type of constraint that is considered.  For example, when beginning with the embedded surface of genus 0 seen in Figure~\ref{fig:weirddog}, enforcing either volume or area preservation separately during the 2-Willmore flow produces a spherical minimizer.   On the other hand, Figure~\ref{fig:weirddog} displays the behavior when this flow is constrained by both surface area and enclosed volume together.  This scenario arises frequently in mathematical biology when considering membrane behavior in an external solution, and minimizing surfaces often realize familiar shapes\textemdash such as the biconcave discoid seen here, which is typical of red blood cells.  See \cite{ou2014} for further details.

\green{
\begin{remark}
It is not difficult to show that the constrained p-Willmore flow in Problem~\ref{prob:constrainedpwillmore} enjoys the same stability property demonstrated in Theorem~\ref{thm:paramflow}.  To see this, repeat the argument from that proof using (\ref{eq:forWconstrained}) instead of (\ref{eq:forW}), and recall the derivatives of the area and volume functionals given previously.
\end{remark}}

\section{Building the mesh regularization equations}
One of the main questions that arises in the computer implementation of curvature flows is how to preserve the quality of the surface mesh as it evolves.  If the initial mesh becomes sufficiently degenerate along the flow, it will crash the simulation\textemdash sometimes well before any troublesome behavior occurs in the actual surface geometry (c.f. Figure~\ref{fig:brokeknot}).  Since curvature flows often alter the initial surface quite dramatically, this can present a serious issue for accurately modeling flow behavior.  Several different techniques have been developed to combat this issue e.g. \cite{crane2011,levy2002,desbrun2002,gu2003,floater2005}, each with their own strengths and weaknesses.  A common challenge present in all methods of mesh regularization is striking a good balance between area-preservation and conformality, or angle-preservation.  Of course, area-preserving maps can be arbitrarily ugly (think following the flow of a vector field tangent to the surface) and conformal maps frequently distort area in undesirable ways. Therefore, the technique employed in this work is inspired by the least-squares conformal mapping procedure of \cite{levy2002} as well as the following result from \cite{kamberov1996}, which can be thought of as a generalization of the Cauchy-Riemann equations from classical complex analysis \green{(see Appendix~\ref{app:CR} for more details).}
\begin{theorem}(Kamberov, Pedit, Pinkall)\label{thm:conformal}
Let $X: M \to \text{Im}\,\mathbb{H}$ be an immersion of \green{the orientable surface} $M$ into the imaginary part of the quaternions, and let $J$ be a complex structure (rotation operator $J^2 = -\text{Id}_{TM}$) on $TM$.  Then, if $\ast\alpha = \alpha \circ J$ is minus the usual Hodge star on differential forms, it follows that $X$ is conformal if and only if there is a Gauss map (unit normal field) $N: M \to \text{Im}\,\mathbb{H}$ such that $\ast dX = N\, dX$.
\end{theorem}
\begin{figure*}
\begin{center}
    \begin{minipage}[c]{0.4\textwidth}
    \includegraphics[width=\textwidth]{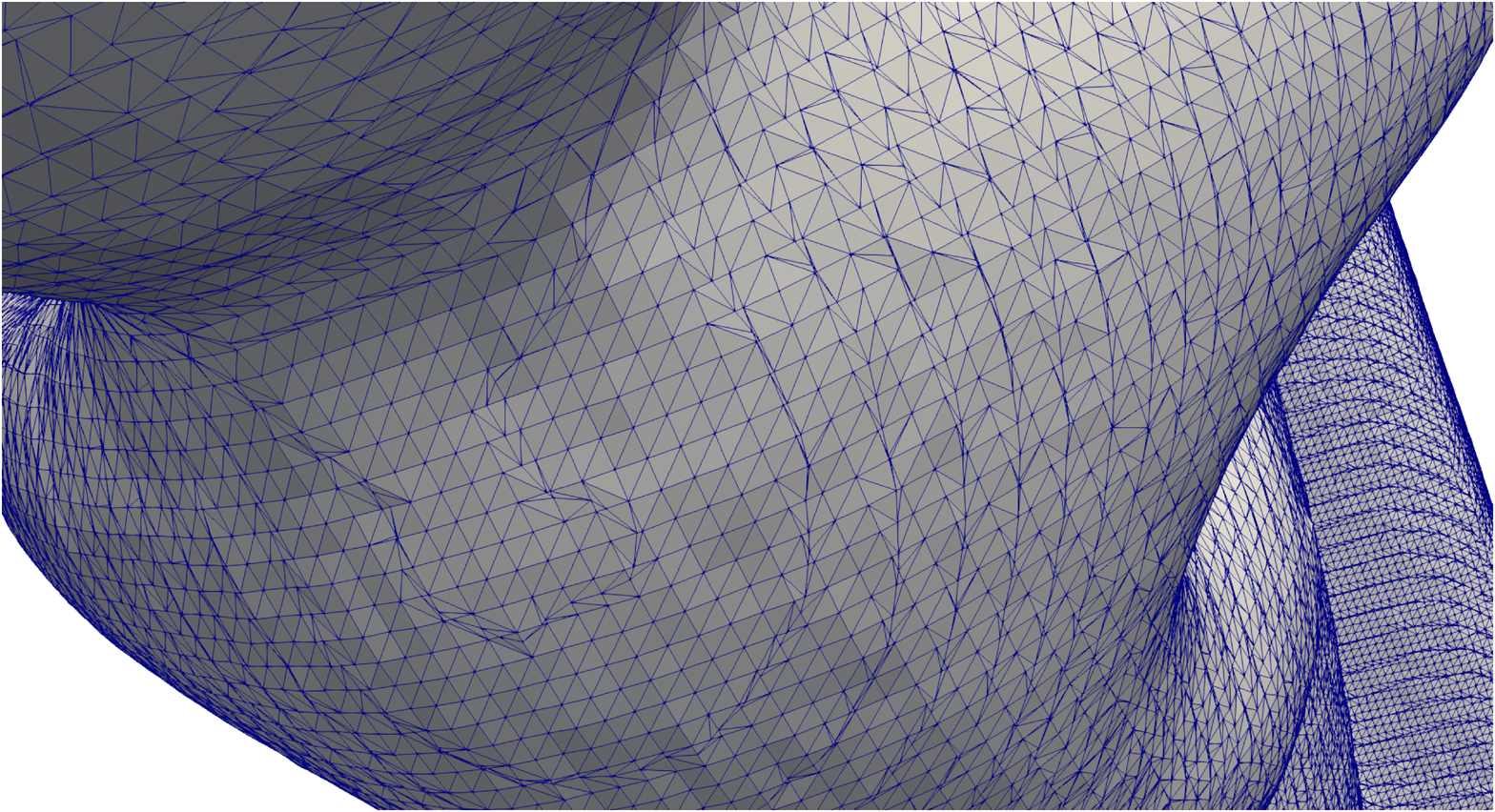}
    \end{minipage}
    \hspace{2pc}
    \begin{minipage}[c]{0.4\textwidth}
    \includegraphics[width=\textwidth]{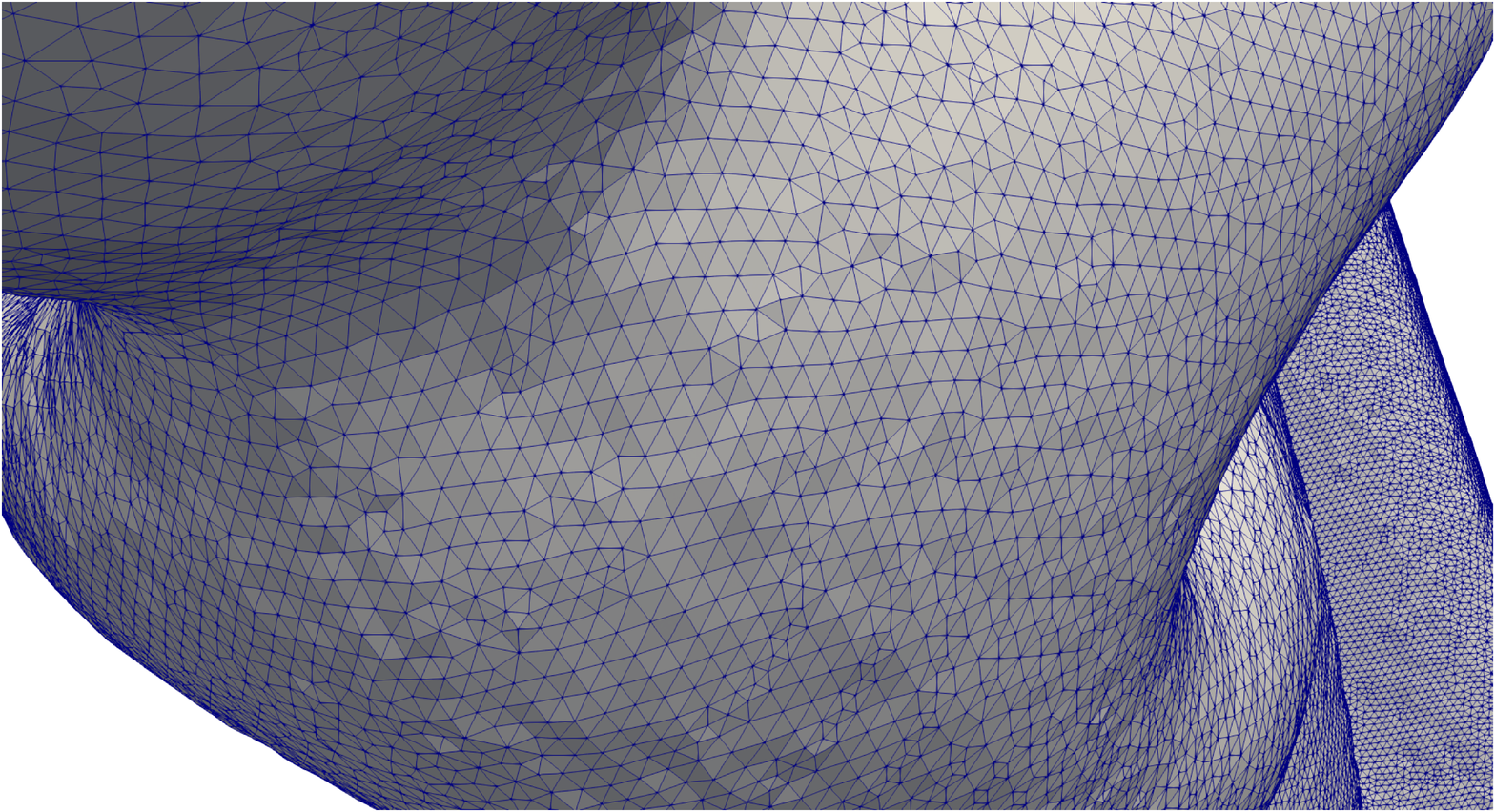}
    \end{minipage} 
\end{center}
\caption[Mesh regularization applied to a statue]{The (linear) procedure described in Problem~\ref{prob:meshreg} applied to a statue mesh of 483k triangles, with close-up on the back of one figure.  Before (left) and after (right).}
\label{fig:meshcomp}
\end{figure*}
Since $\text{Im}\,\mathbb{H}$ is canonically isomorphic to $\mathbb{R}^3$ as a vector space, this gives a criterion for conformality that can be weakly enforced during the p-Willmore flow.  More precisely, recall that $N\perp dX(v)$ for all tangent \green{vector fields} $v \in TM$, and that multiplication of $v,w\in\text{Im}\,\mathbb{H}$ obeys the rule $vw = -\langle v, w\rangle + v\times w$ where $\times$ is the usual vector cross product.  It follows that $N\, dX(v) = N \times dX(v)$ for all $v\in TM$, and the conformality condition in Theorem~\ref{thm:conformal} can be expressed \green{pointwise} as
\begin{equation}\label{eq:confcond}
    \ast dX(v) = N \times dX(v) \qquad \text{for all}\,\,v\in TM.
\end{equation}
For the present purpose of mesh regularization, it suffices to enforce this condition weakly through a minimization procedure.  Define the conformal distortion functional
\begin{equation*}
    \mathcal{CD}(u) = \frac{1}{2}\int_M \left|du\,J - N \times du \right|^2 d\mu_g.
\end{equation*}
In view of Theorem~\ref{thm:conformal}, $\mathcal{CD}(u)$ is identically 0 if and only if $u$ is a conformal immersion of $M$ into $\mathbb{R}^3$.
\greener{Assuming the surface metric is held fixed,} minimization of the conformal distortion leads to the necessary condition
\begin{equation}\label{eq:conformal}
    \delta\mathcal{CD}(u)\varphi = \int_M \big\langle \left(du\,J - N \times du\right), \left(d\varphi\,J - N\times d\varphi\right) \big\rangle_g\, d\mu_g = 0,
\end{equation}
which must hold for all $\varphi \in H^1(M,\mathbb{R}^3)$.

\greener{
\begin{remark}
The reader will notice that the evolution-dependent nature of the mesh regularization problem has been ignored. As the goal is compute a map very close to $u(M)$ itself, the metric (and hence the volume element) associated to the new immersion are approximated using the respective quantities coming from the present immersion.
\end{remark}}

To make use of this equation in the computational framework considered here, it must first be expressed in terms of the local coordinates on $U \subset \mathbb{R}^2$.  To that end, a particular section $v\in TM$ of the ``complex line bundle'' $TM$ is first chosen; it is advantageous to consider the parametrization domain $U$ and choose $X^*v = e_1$ where $e_1$ represents the first standard basis vector for $TU$. Then, (abusing $J$ to denote the pullback complex structure on $TU$ as well as the original structure on $TM$), $e_2 = J(e_1)$ and the integrand of (\ref{eq:conformal}) applied to the basis $dX(e_i)$ for $TM$ can be pulled back through $X$ to yield the coordinate expression 
\begin{equation}\label{eq:coordconf}
    g^{ij}\IPA{\left(dX\,J(e_i) - N \times dX(e_i)\right)}{\left(dX\,J(e_j) - N \times dX(e_j)\right)}, 
\end{equation}
where $N$ is understood to mean the pullback (outer) normal field
\begin{equation*}
    N = \frac{X_1 \times X_2}{\lvert X_1 \times X_2\rvert},
\end{equation*}
which is valid on the parametrization domain $U$.  To write this in a more compact form, first notice that
\begin{align*}
    dX\,J(e_1) - N \times dX(e_1) &= \,dX(e_2) - N \times dX(e_1), \\
    dX\,J(e_2) - N \times dX(e_2) &= -\left(dX(e_1) + N \times dX(e_2)\right).
\end{align*}
So, letting $N = (n_1, n_2, n_3)^T$ and denoting $X^i_j = \left\langle e_i, dX(e_j)\right\rangle$ (where we have abused notation again by referring to $e_i$ as a standard basis field on both $U$ and $M$), one can form the vectors 
\begin{align*}
&V = dX(e_2) - N \times dX(e_1) = 
  \begin{pmatrix}
    X^1_2 - n_2\, X^3_1 + n_3\, X^2_1 \\ 
    X^2_2 - n_3\, X^1_1 + n_1\, X^3_1 \\
    X^3_2 - n_1\, X^2_1 + n_2\, X^1_1
  \end{pmatrix}, \\
&W = dX(e_1) + N \times dX(e_2) = 
  \begin{pmatrix}
    X^1_1 + n_2\, X^3_2 - n_3\, X^2_2 \\ 
    X^2_1 + n_3\, X^1_2 - n_1\, X^3_2 \\
    X^3_1 + n_1\, X^2_2 - n_2\, X^1_2
  \end{pmatrix}.
\end{align*}
With this, careful reorganization of (\ref{eq:coordconf}) yields the $\mathbb{R}^{3x2}$ dyadic product $\tilde{Q}:\mathrm{Jac}\,\varphi$, where $\mathrm{Jac}$ denotes the usual Jacobian and $\tilde{Q}$ is given in components as (indices $1\leq i \leq 3$, mod 3)
\begin{align*}
    \hat{Q}^i_1 &= g^{22}\,W_i + g^{11}\left( n_{i+1}\,V_{i+2} - n_{i+2}\,V_{i+1} \right) \\
    &+ g^{12}\left( n_{i+2}\,W_{i+1} - n_{i+1}\,W_{i+2} - V_i \right), \\
    \hat{Q}^i_2 &= g^{11}\,V_i + g^{22}\left( n_{i+2}\,W_{i+1} - n_{i+1}\,W_{i+2} \right) \\
    &+ g^{12}\left( n_{i+1}\,V_{i+2} - n_{i+2}\,V_{i+1} - W_i \right).
\end{align*}
For cleanliness of presentation, note that there is a tensor $\greener{Q = Q(u)} \in T^*U \otimes TM$ such that $\hat{Q}^{Ij} = g^{kj}Q_{k}^I$, so we may write this product (at least formally) as $\IP{Q}{d\varphi}$.  Therefore, equation (\ref{eq:conformal}) can now be expressed concisely as
\begin{equation}\label{eq:conciseconf}
\int_M \IP{\greener{Q(u)}}{d\varphi}d\mu_g = 0.
\end{equation}
Equation \ref{eq:conciseconf} is used to ensure that the surface mesh finds a configuration that is ``as conformal as possible'' to a specified discretization of the reference \greener{surface $M \subset \mathbb{R}^3$.}  In practice, best results (particularly for triangle meshes) are achieved if this reference is \greener{defined implicitly based on an adjustment of the initial mesh data.  In particular, it is useful to choose the reference configuration to be the starting surface with interior angles adjusted relative to the number of elements sharing a vertex.  More precisely, there is the following procedure.  For each vertex $i$ in the triangulation, first compute the number $m_i$ of elements with $i$ as a vertex.  Then, if $T$ is a triangle adjacent at $i$ with vertices $i_1,i_2,i_3$ and interior angles $\alpha_{i_1},\alpha_{i_2},\alpha_{i_3}$, rescale $\alpha_{i_k} \mapsto \alpha_{i_k}/m_i$ to generate ``ideal'' interior angles.}  It is clear that, in general, $\alpha_{i_1}+\alpha_{i_2}+\alpha_{i_3} \neq \pi$, which is necessary for closure.  The strategy is therefore to use the largest interior angle of each triangle to adjust the others.  In particular, suppose $\alpha_{i_1} > \alpha_{i_2}$ and $\alpha_{i_1} > \alpha_{i_3}$.  Then, set 
\[ \alpha_{i_k} \leftarrow \frac{\pi - \alpha_{i_1}}{\alpha_{i_2} + \alpha_{i_3}}\,\alpha_{i_k}, \qquad k = 2,3. \]
Moreover, in the case that a given triangle has two or three leading angles, each angle is set to $\pi/3$.  Pseudocode for this procedure is given in Algorithm~\ref{alg:envn}.

\begin{algorithm}
\caption{Generation of reference angles}
\label{alg:envn}
\green{
\begin{algorithmic}
\REQUIRE Reference triangulation $\mathcal{T}$ of the closed surface $M$.
\FOR{$T\in\mathcal{T}$}
    \FOR{vertex $1 \leq i \leq 3$}
        \STATE Compute $m_i =$ \# of adjacent elements
        \STATE $\alpha_i \leftarrow \alpha_i / m_i$
    \ENDFOR
    \STATE Determine maximum vertex angle $\alpha_i$.
    \IF{$\alpha_i > \alpha_j$ for all $j \neq i$}
        \FOR{vertices $j \neq i$}
            \STATE $\alpha_j \leftarrow \alpha_j \left(\pi - \alpha_{i}\right) / \left(\sum_{k\neq i} \alpha_{k}\right)$
        \ENDFOR
    \ELSE
        \FOR{vertices $1 \leq j \leq 3$}
            \STATE $\alpha_j \leftarrow \alpha_j\,\pi\, / \left(\sum_{k=1}^3 \alpha_{k}\right)$.
        \ENDFOR
    \ENDIF
\ENDFOR
\end{algorithmic}}
\end{algorithm}

Of course, if equation (\ref{eq:conciseconf}) is to be useful as an effective tangential  reparametrization of a surface evolving by  p-Willmore flow, it should not be solved without constraint. For example, it is clear that any constant function $u$ will satisfy this equation as stated, \green{so some care must be taken to prevent trivialities.  Moreover, it is also necessary to constrain the regularization in (\ref{eq:conciseconf}) so that it does not destroy the current surface geometry by moving the surface too far in the normal direction.  One potential solution to this issue is motivated by the following observation}:  if the aim is to recover a new immersion $\hat{u}$ which is ``close'' to the current immersion $u$, then the difference $\hat{u}(x) - u(x)$ at any $x\in M$ should be tangential to first order, hence orthogonal to the surface normal $N(x)$.  \green{ Said differently, a first-order approximation to tangential motion along the surface can be obtained by requiring that the pointwise equation
\begin{equation*}\label{eq:confmult}
    \IPA{\left(\hat{u}(x) - u(x)\right)}{N(x)}=0,
\end{equation*}
hold for all $x\in M$ during the above minimization. In fact, since exact conformality is not required for the present purpose of mesh regularization, it is advantageous to weaken this requirement further using a penalty term.  Since saddle-point problems with a mixture of linear and piecewise-constant finite elements can exhibit unstable discretizations, the inclusion of such a term helps to prevent numerical artifacts from appearing during the implementation. Precisely, the conformal penalty regularization procedure is presented as the following problem.}


\begin{problem}[Conformal penalty regularization]\label{prob:meshreg}
\green{
Given a fixed $\varepsilon>0$ and an oriented surface immersion $u: M \to \mathbb{R}^3$ with outward unit normal field $N$, solving the conformally-penalized mesh regularization problem amounts to finding a function $v:M \to \mathbb{R}^3$ and a Lagrange multiplier $\rho:M \to \mathbb{R}$, so that the new immersion $\hat{u} = u + v$ is the solution to
\begin{align*}
    &\min_v \left(\mathcal{CD}(u + v) + \greener{ \frac{\varepsilon}{2} \int_M \rho^2\,d\mu_g + \int_M \rho \IPA{v}{N}\,d\mu_g} \right)
\end{align*}}
Formulated weakly, the goal becomes to find a new immersion $\hat{u}$ and a multiplier $\rho$ which satisfy the system 
\green{
\begin{align*}
    0 &= \int_M \rho\IPA{\varphi}{N}d\mu_g + \int_M \IP{\greener{Q(u + v)}}{d\varphi} d\mu_g, \\
    0 &= \int_{M} \psi\IPA{v}{N}d\mu_g + \varepsilon\int_M \psi\rho\,d\mu_g,
\end{align*}
for all $\varphi,\psi\in H^1(M; \mathbb{R}^3)$.}
\end{problem}


Solving Problem~\ref{prob:meshreg} at each step of the p-Willmore flow inhibits the computer simulation from breaking arbitrarily, at the expense \green{of potentially altering the flow solution at each time step}.  To be sure, without such a procedure in place it is not unusual for global minimizers to remain computationally out-of-reach, as is shown in Figure~\ref{fig:brokeknot}.  Moreover,  Figure~\ref{fig:meshcomp} demonstrates how this regularization is useful even for stationary surfaces, as it greatly improves the quality of (perhaps very irregular) surface meshes.  The practical discretization of both this system and the p-Willmore system in Problem~\ref{prob:constrainedpwillmore} will be discussed in the next sections.

\section{Discretization of model systems}
Discretization of the models in Problem~\ref{prob:constrainedpwillmore} and Problem~\ref{prob:meshreg} will now be discussed.  In particular, specifics of the spatial and temporal discretization are presented, leading to appropriate discrete versions of the continuous problems above.  Moreover, some insight is given into the treatment of nonlinearities, and the main algorithm of this work is given.

\begin{figure}
\begin{center}
\begin{minipage}[c]{0.11\textwidth}
\includegraphics[width=\textwidth]{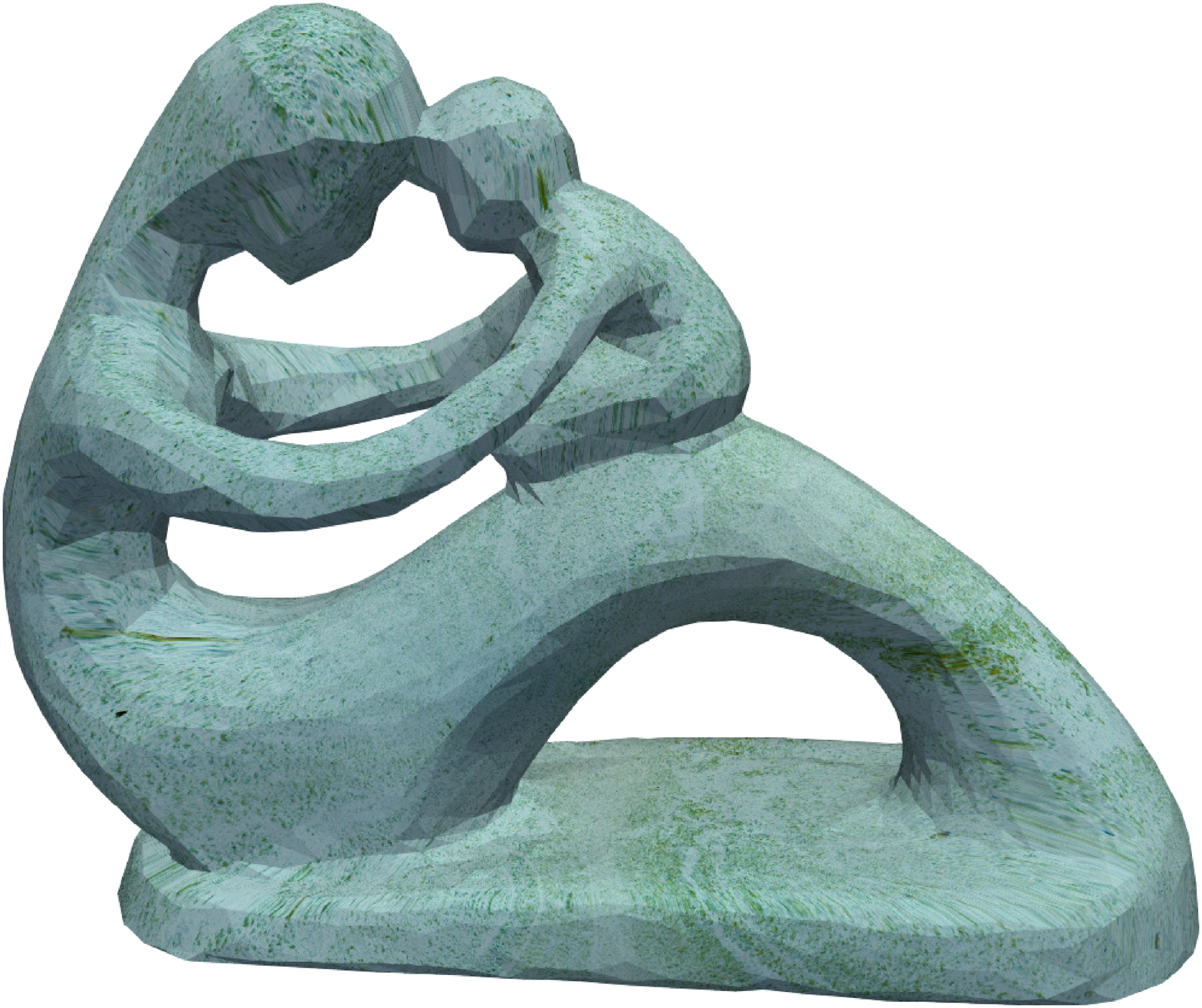}
\end{minipage}
\begin{minipage}[c]{0.11\textwidth}
\includegraphics[width=\textwidth]{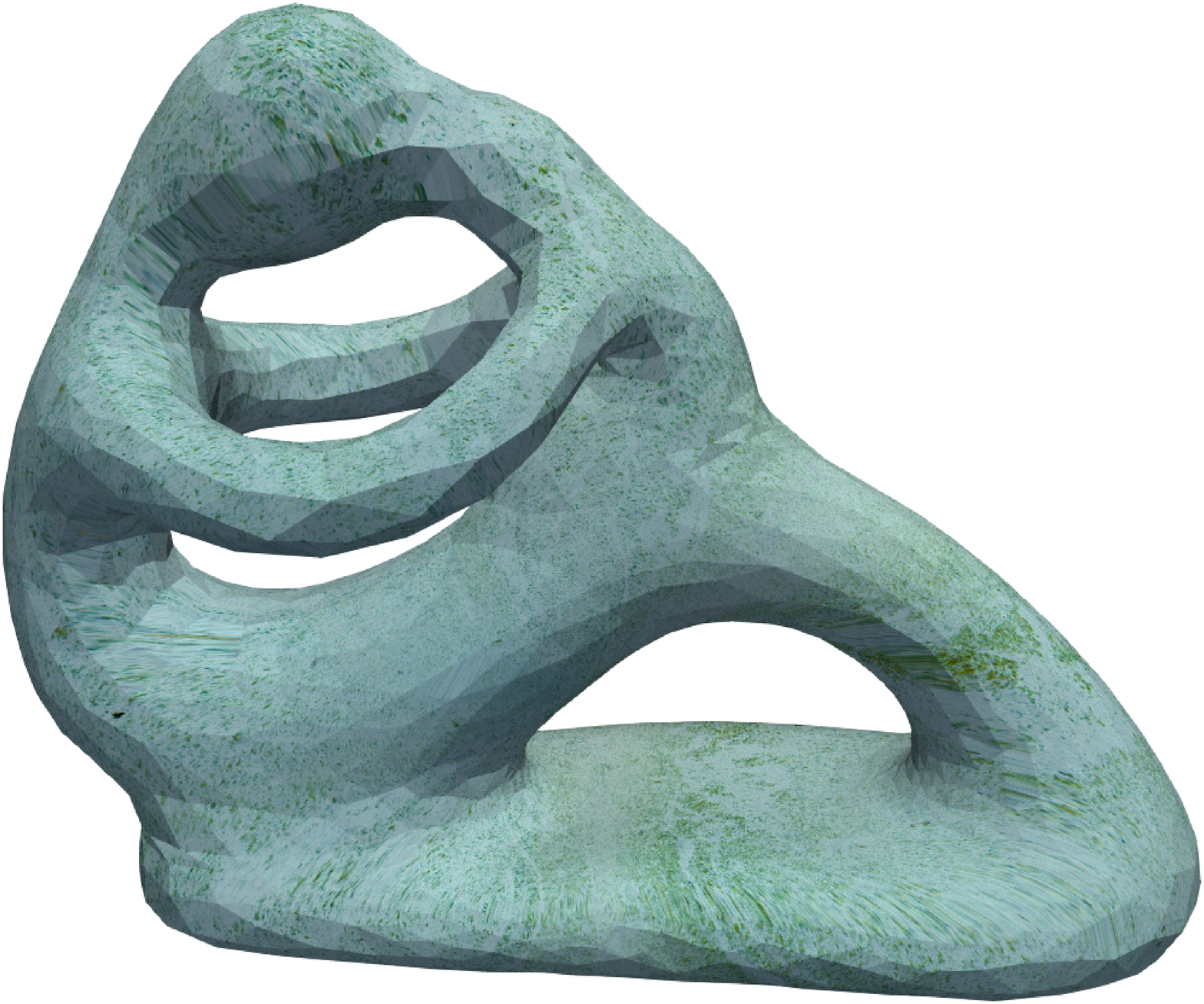}
\end{minipage}
\begin{minipage}[c]{0.11\textwidth}
\includegraphics[width=\textwidth]{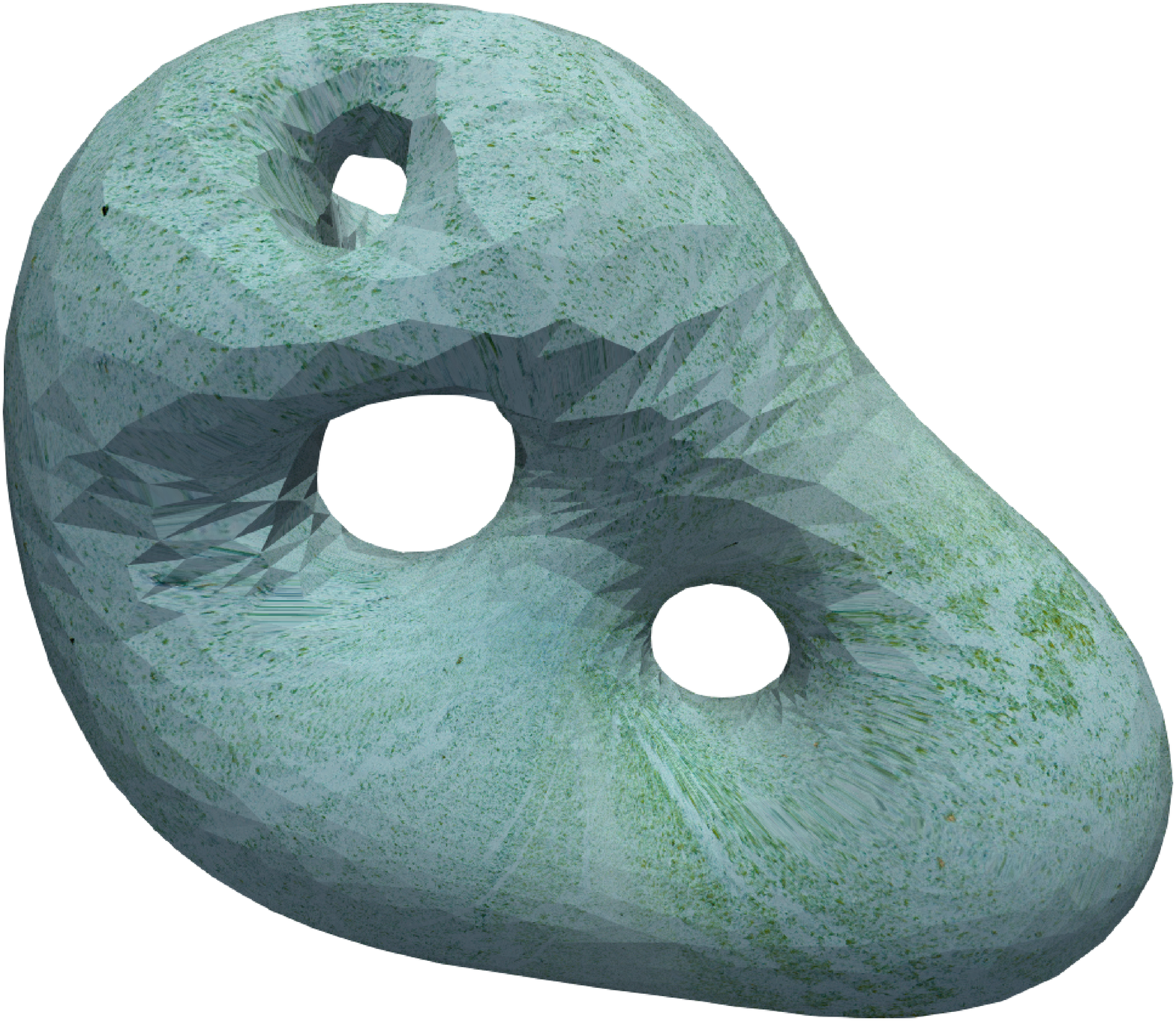}
\end{minipage}
\begin{minipage}[c]{0.11\textwidth}
\includegraphics[width=\textwidth]{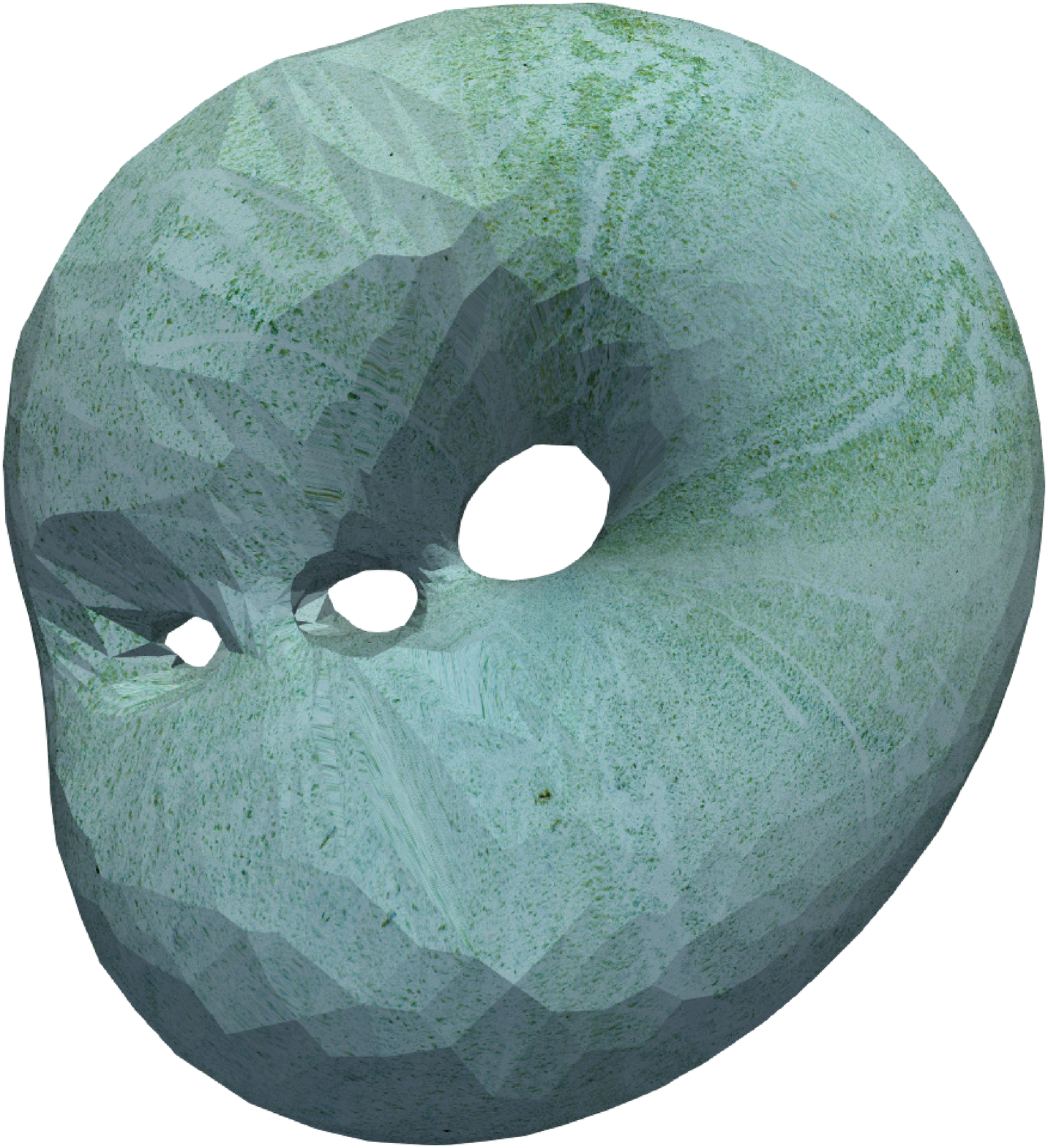}
\end{minipage}
\end{center}
\caption[2-Willmore flow of a statue]{Unconstrained 2-Willmore flow of a genus 4 statue mesh (last image rotated 90 degrees counterclockwise).}
\label{fig:statueflow}
\end{figure}

As in \cite{dziuk2013}, the smooth surface $M \subset \mathbb{R}^3$ is assumed to be approximated by a polygonal surface $M_h$ consisting of 2-simplices (triangles) $T_h$ that are not degenerate, so that
\begin{equation*}
    M_h = \bigcup_{T_h \in \mathcal{T}_h} T_h,
\end{equation*}
forms an admissible triangulation of $M$.  Denoting the nodes of this triangulation by $\{a_j\}_{j=1}^N$, the standard nodal basis $\{\phi_i\}$ on $M$ satisfies $\phi_i(a_j,t) = \delta_{ij}$.  The space of piecewise-linear finite elements on $M_h(t)$ is then denoted
\begin{equation*}
    S_h(t) = \text{Span}\{\phi_i\} = \{\phi \in C^0(M_h(t))\, : \, \phi|_{T_h} \in \mathbb{P}_1(T_h), T_h \in \mathcal{T}_h\},
\end{equation*}
where $\mathbb{P}_1(T_h)$ denotes the space of linear polynomials on $T_h$.

\begin{figure*}
\begin{center}
    \begin{minipage}[c]{0.22\textwidth}
    \includegraphics[width=\textwidth]{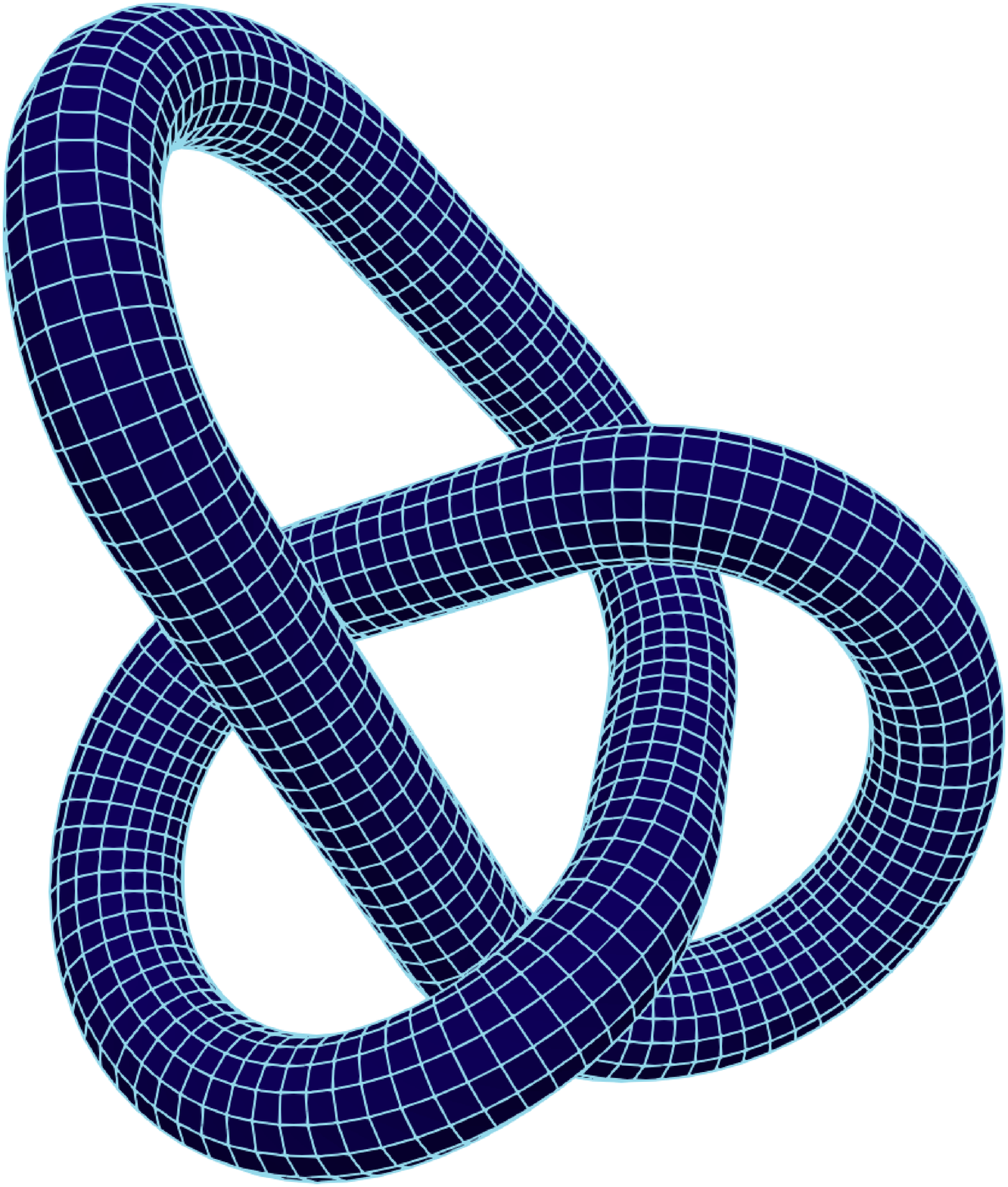}
    \end{minipage}
    \begin{minipage}[c]{0.20\textwidth}
    \includegraphics[width=\textwidth]{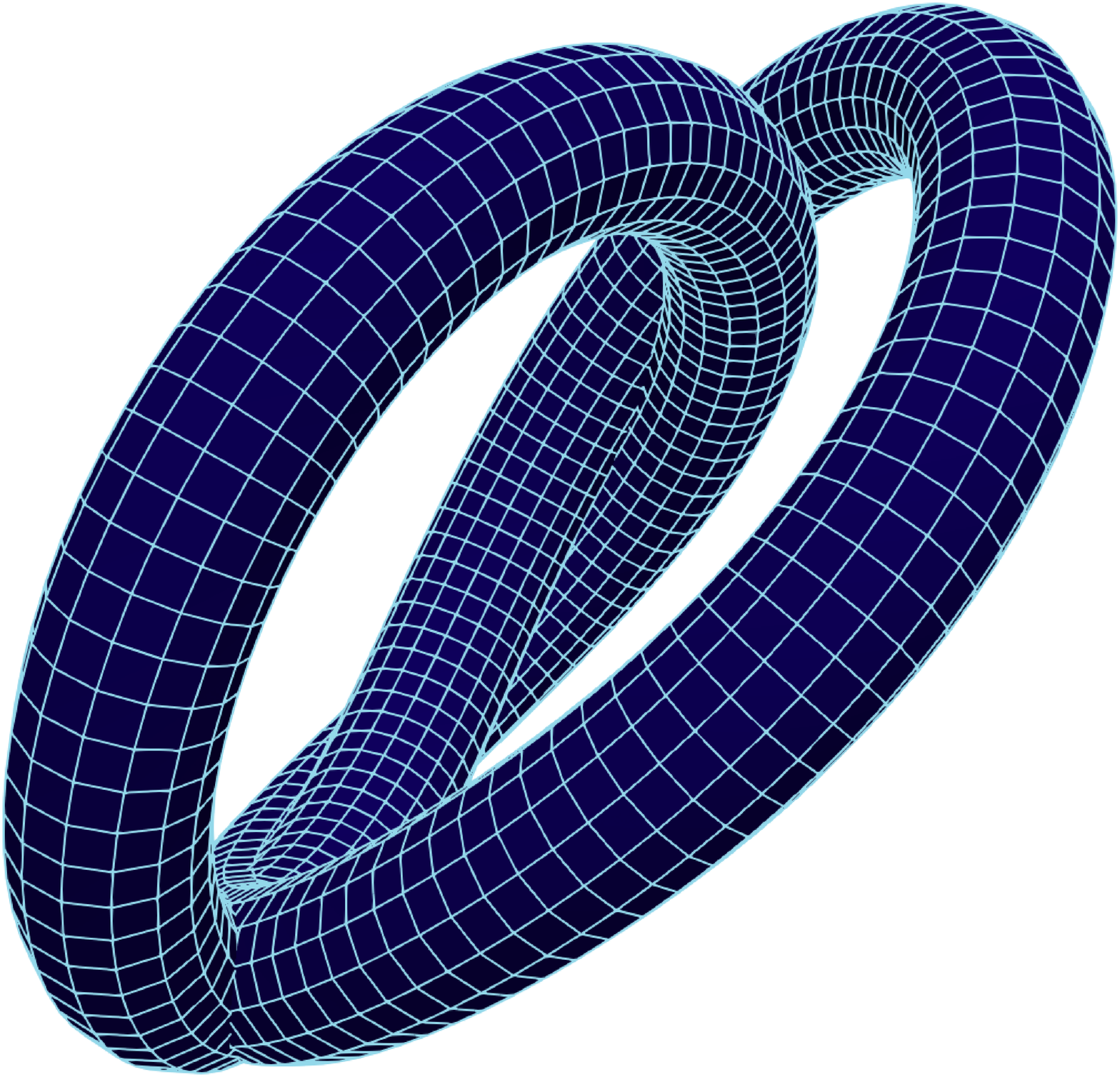}
    \end{minipage}
    \begin{minipage}[c]{0.22\textwidth}
    \includegraphics[width=\textwidth]{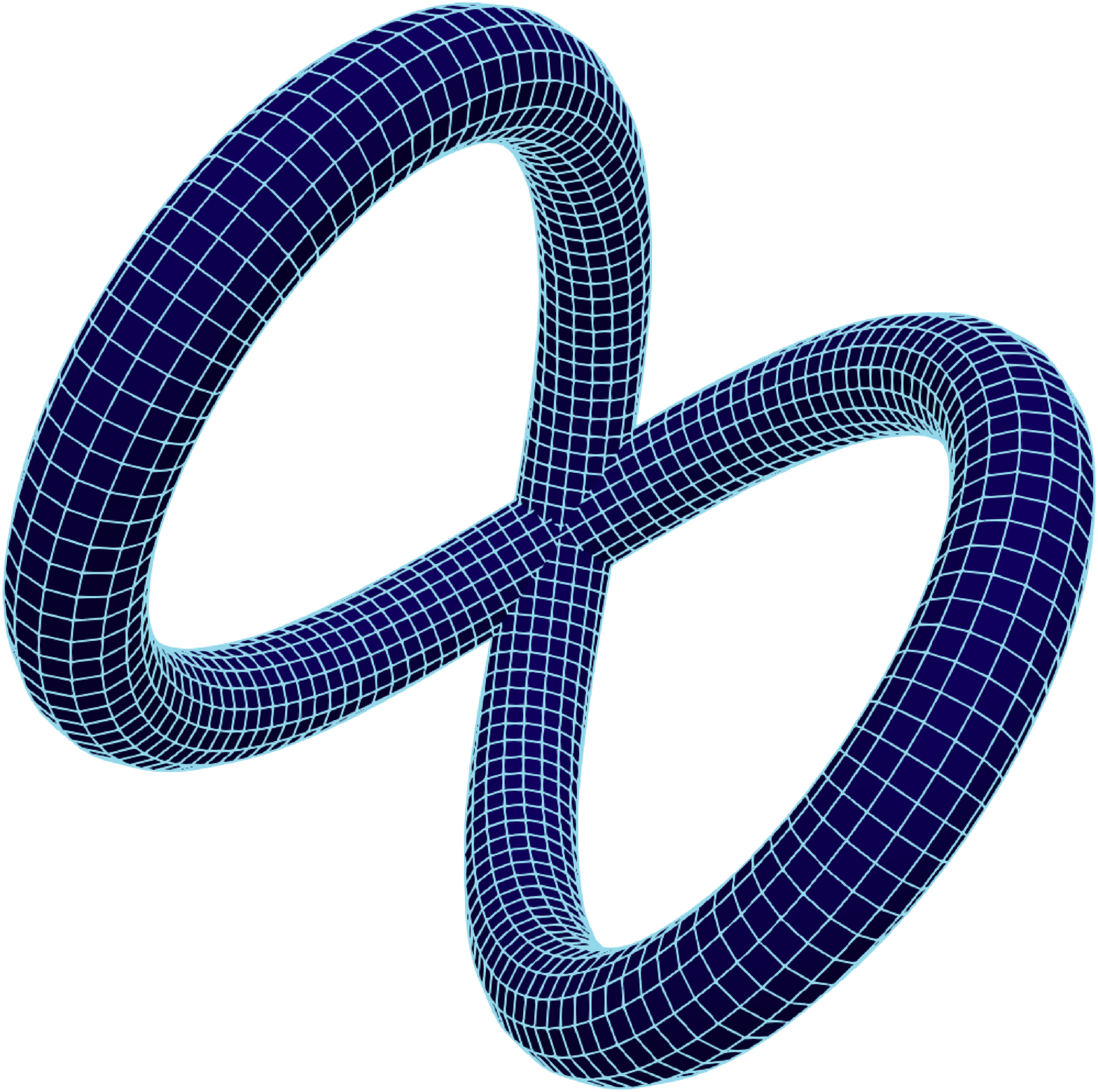}
    \end{minipage}
    \hspace{0.6pc}
    \begin{minipage}[c]{0.22\textwidth}
    \includegraphics[width=\textwidth]{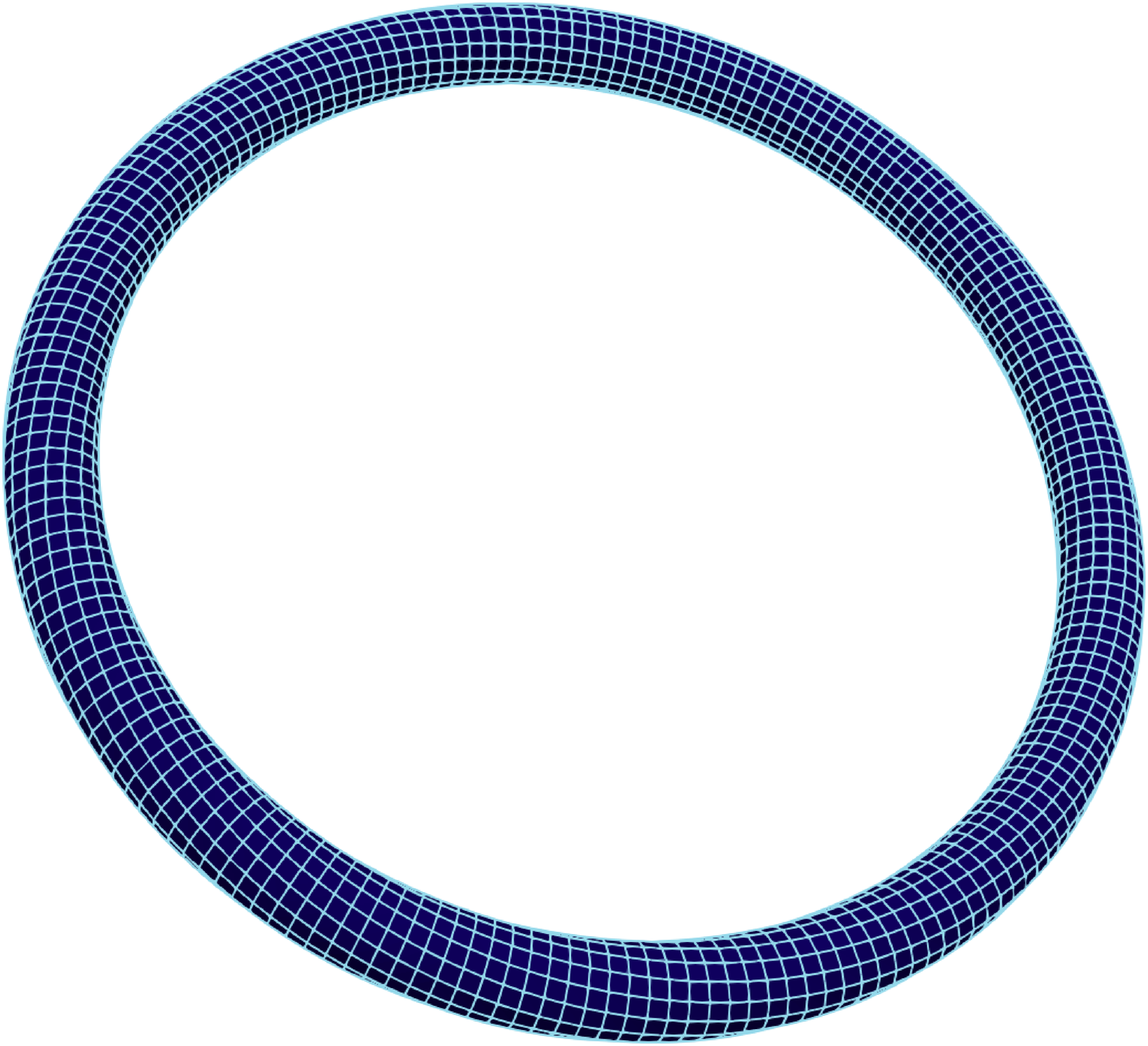}
    \end{minipage}
    \\
    \begin{minipage}[c]{0.22\textwidth}
    \includegraphics[width=\textwidth]{figs/knot.eps}
    \end{minipage}
    \begin{minipage}[c]{0.20\textwidth}
    \includegraphics[width=\textwidth]{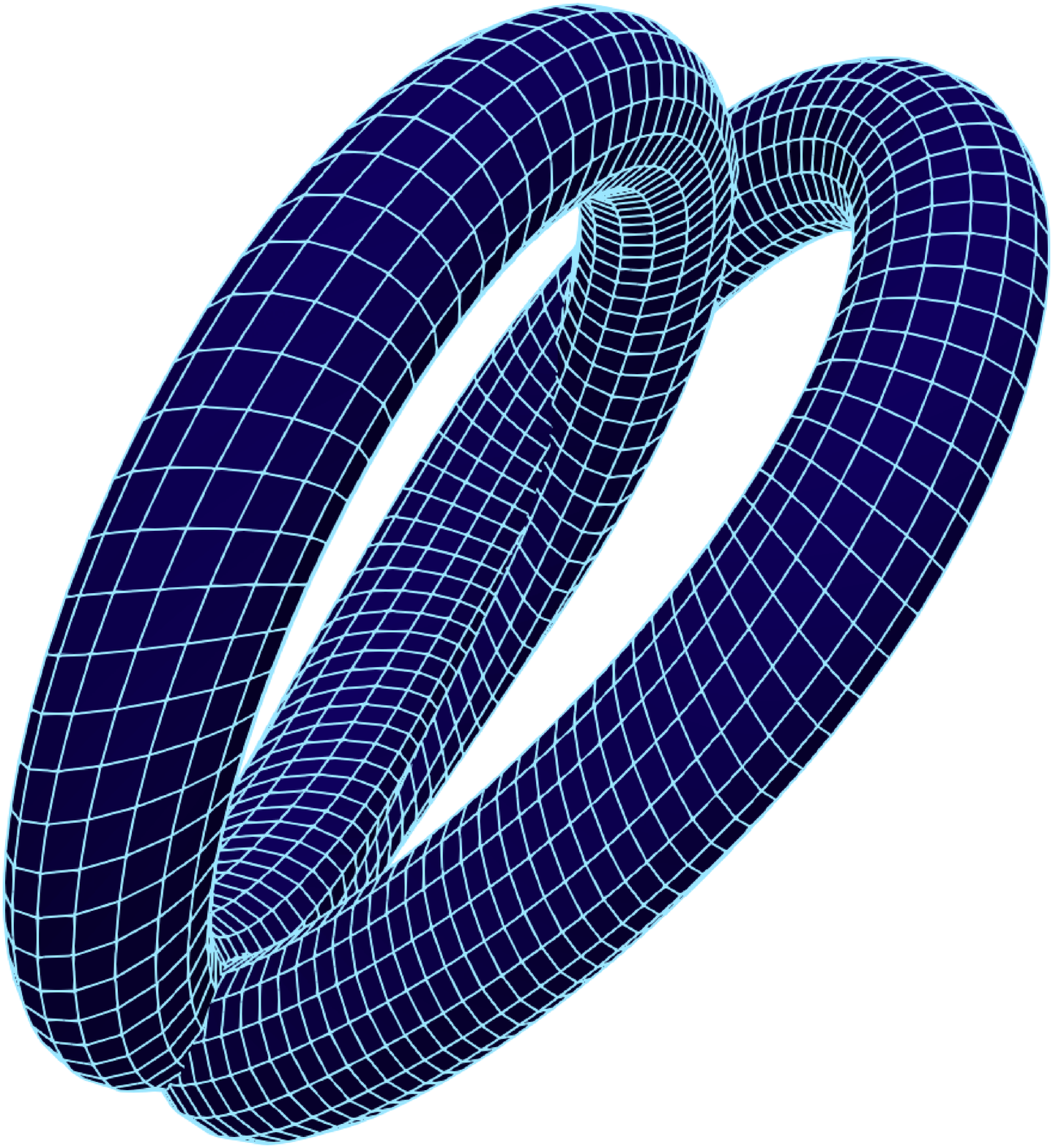}
    \end{minipage}
    \begin{minipage}[c]{0.22\textwidth}
    \includegraphics[width=\textwidth]{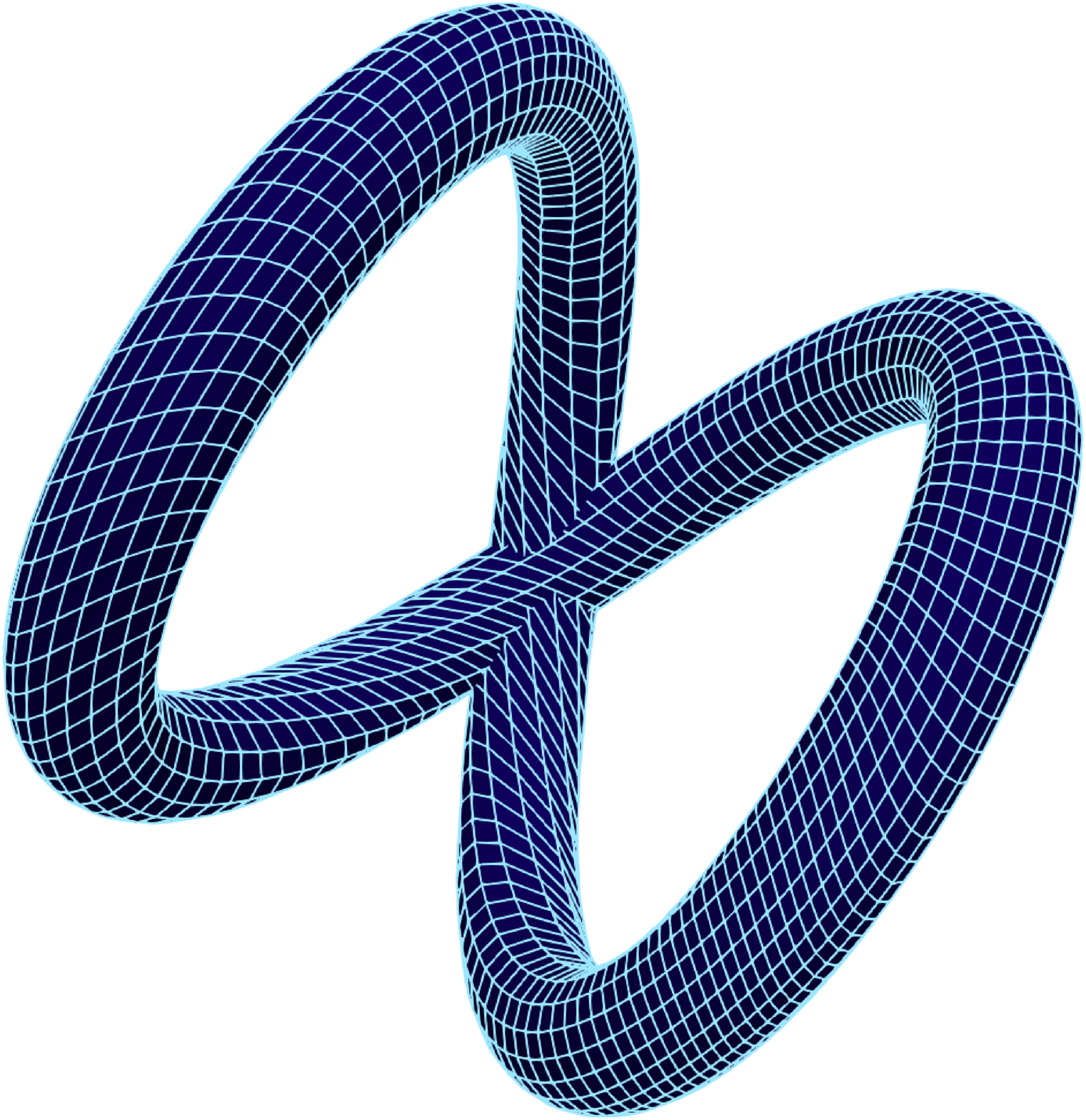}
    \end{minipage}
    \hspace{0.6pc}
    \begin{minipage}[c]{0.22\textwidth}
    \includegraphics[width=\textwidth]{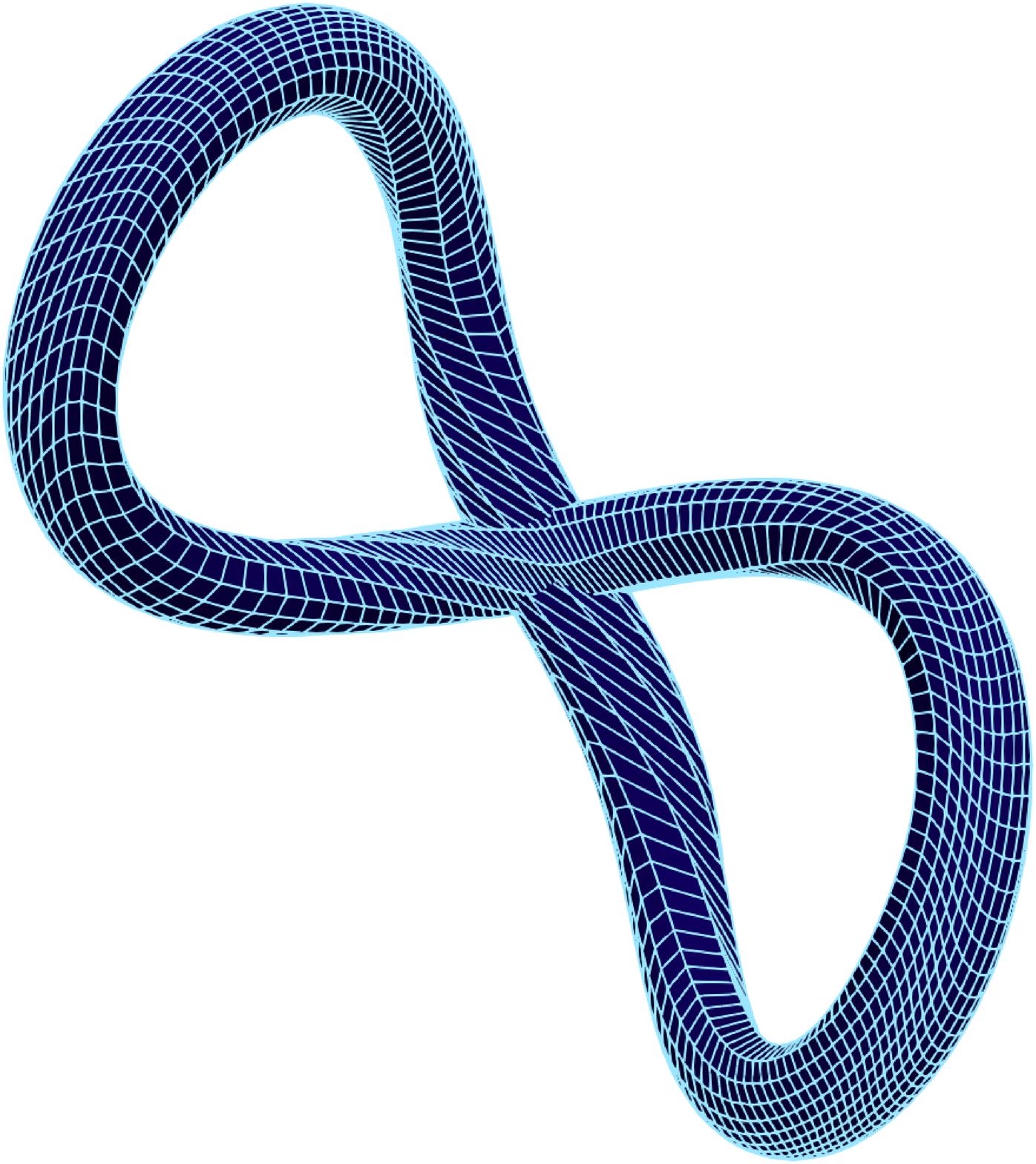}
    \end{minipage}
\end{center}
\caption[Constrained 2-Willmore evolution of a trefoil knot]{Surface area and volume constrained 2-Willmore evolution of a trefoil knot, with conformal mesh regularization (top) and without (bottom).  \greener{Due to a coarse discretization,} the mesh degenerates when not regularized, preventing movement to the minimizing surface.  Conversely, element quality remains nearly perfect along the flow when regularized at each step, and evolves completely to the desired minimum \greener{despite the mesh resolution.}}
\label{fig:brokeknot}
\end{figure*}

Now, suppose the simplices $T_h \in \mathcal{T}_h$ have maximum diameter $h$ with inner radius bounded from below by $ch$ for some $c>0$. Then, for any choice of unit normal $N$ there is $\delta >0$ so that points $x\in\mathbb{R}^3$ in the tube of radius delta around $M$ can be expressed as 
\begin{equation*}
    x = a(x) + d(x)N(x),
\end{equation*}
where $a(x)$ lies on $M$ and $|d|<\delta$.  It is assumed that $M_h$ is contained in this tube, so that any function $f_h$ defined on the discrete surface $M_h$ can be lifted to a function $f_h^l$ on the smooth surface $M$ by requiring
\begin{equation*}
    f_h^l(a(x)) = f_h(x), \qquad x\in M_h.
\end{equation*}
Denote the inverse process by $f^{-l}$.  
In this notation, the lift of the finite element space $S_h$ is denoted 
\begin{equation*}
    S_h^l = \{\phi^l \, |\, \phi\in S_h\},
\end{equation*}
and it is possible to compare geometric quantities on $M$ and $M_h$.  \green{In particular, let $g_h$ denote the induced metric on $M_h$, $u_h$ denote identity on $M_h$, and $Y_h = \Delta_{g_h} u_h$ denote the discrete  mean curvature vector defined through the relationship
\begin{equation*}
    0 = \int_{M_h} \IPA{Y_h}{\psi_h} d\mu_{g_h} + \int_{M_h} \IPH{du_h}{d\psi_h} d\mu_{g_h}, \quad\forall \psi_h \in S_h,
\end{equation*}
where $d = d_h$ is understood to mean the derivative operator with respect to the surface $M_h$ and $d\mu_{g_h}$ is the area element with respect to the  metric $g_h$.
Then, there are the discrete area and volume functionals
\begin{align*}
    \mathcal{A}(u_h) &= \int_{M_h} 1\,d\mu_{g_h}, \\
    \mathcal{V}(u_h) &= \int_{M_h} \IPA{u_h}{N_h} d\mu_{g_h},
\end{align*}
and the discrete $p$-Willmore functional is defined to be
\begin{equation*}
    \mathcal{W}^p(u_h) = \frac{1}{2^p}\int_{M_h} \left\lvert Y_h\right\rvert^p.
\end{equation*}
Moreover, if $J$ is understood to be a linear operator on the tangent space $TM_h$, the discrete conformal distortion functional can be defined as
\begin{equation*}
    \mathcal{CD}(u_h) = \frac{1}{2}\int_{M_h} \left|du_h\,J - N_h \times du_h \right|^2 d\mu_{g_h}.
\end{equation*}}

\green{These definitions are seen to yield a reasonable spatial discretization of the continuous Problems~\ref{prob:constrainedpwillmore} and \ref{prob:meshreg}, but it remains to discuss their time-dependent aspect as well.  A reasonable strategy for the temporal discretization of complicated nonlinear PDEs on evolving surfaces is to linearize the problem at each time step, effectively pushing the nonlinearities into the temporal domain.  This is the strategy of \cite{dziuk2008,dziuk2013}.  Though such formulations enjoy many of the benefits of their linear counterparts, they typically require a very small time step and are less robust to round-off errors as well as other sources of numerical instability.  Therefore, the present strategy for discretizing the p-Willmore flow Problem~\ref{prob:constrainedpwillmore} is to ``center'' the discretization in time, except for some isolated terms for which this is problematic.  Though this approach will still produce a numerical scheme which is first-order in time, it is seen to significantly improve the stability of the fully-discrete p-Willmore flow.  To describe this idea more precisely,  let $\tau > 0$ be a fixed temporal stepsize, and denote $u_h^k \coloneqq u_h(\cdot, k\tau)$.  Consider that at time step $k$, the images of $u_h^k$ and $u_h^{k+1}$ are the old resp. new surfaces $M_h^k$ and $M_h^{k+1}$.  Therefore, $M_h^\half$ will denote the central surface defined by the immersion $u_h^\half = (1/2)\left(u_h^k + u_h^{k+1}\right)$, and a field quantity $T_h$ defined on each surface will have the central counterpart $T_h^\half = (1/2)\left(T_h^k + T_h^{k+1}\right)$.  Additionally, $d$ and $D$ will now be used to denote discretized operators with respect to the metric on the central surface.  With this, running the discrete p-Willmore flow with constraint is presented as the following problem.
\begin{problem}[Discrete p-Willmore flow with constraint]\label{prob:discpwillmore}
Let $u,Y,W,\lambda,\mu$ be as in Problem~\ref{prob:constrainedpwillmore}.  Given the discrete data $u_h^k,Y_h^k,W_h^k$ at time $t = k\tau$, the p-Willmore flow problem is to find functions $u_h^{k+1},Y_h^{k+1},W_h^{k+1},\lambda_h,\greener{\gamma_h}$ which satisfy the system of equations
{\allowdisplaybreaks
\begin{align}
    0 &= \int_{M_h^\half} \IPA{Y_h^\half}{\psi_h}d\mu_{g_h} + \int_{M_h^\half} \IPH{du_h^{k+1}}{d\psi_h} d\mu_{g_h}, \notag \\
    0 &= \int_{M_h^\half} \IPA{\left(W_h^\half - \left|Y_h^\half\right|^{p-2}Y_h^\half\right)}{\xi_h} d\mu_{g_h}, \label{eq:forWdisc} \\
    0 &= \int_{M_h^\half} \IPH{du_h^\half}{\left(du_h^{k+1}-du_h^k\right)} d\mu_{g_h}, \label{eq:areacon}\\
  0 &= \int_{M_h^\half} \IPA{\left(u_h^{k+1} - u_h^k\right)}{N_h^\half} d\mu_{g_h}, \label{eq:volcon} \\
    0 &= \int_{M_h^\half} \frac{\IPA{\left(u_h^{k+1}-u_h^k\right)}{\varphi_h}}{\tau} \,d\mu_{g_h} \notag \\
    &+ \int_{M_h^\half} \,\lambda_h\IPA{\varphi_h}{N_h^\half} d\mu_{g_h} \notag \\
    &+ \int_{M_h^\half} \,\greener{\gamma_h}\IPH{du_h^\half}{d\varphi_h} d\mu_{g_h} \notag \\
    &+ (1-p)\int_{M_h^\half} \left|Y_h^\half\right|^p \IPH{du_h^\half}{d\varphi_h} d\mu_{g_h} \notag \\
    &- p\int_{M_h^\half} \left(\divgh W_h^\half\right)\IPH{du_h^\half}{d\varphi_h} d\mu_{g_h} \notag \\
    &- p\int_{M_h^\half} \IPH{dW_h^{k+1}}{d\varphi_h} d\mu_{g_h} \notag \\
    &+ p\int_{M_h^\half} \IPH{D(\varphi_h)du_h^k}{dW_h^k} d\mu_{g_h}, \label{eq:forYdisc}
\end{align}}
for all $\varphi_h,\psi_h, \xi_h \in S_h$.
\end{problem}
\begin{remark}
Similar to the continuous situation, the case $p=0$ can be handled by omitting equations (\ref{eq:forWdisc}), (\ref{eq:areacon}) and replacing equation~(\ref{eq:forYdisc}) with 
\begin{align*}
        0 &= \int_{M_h^\half} \frac{\IPA{\left(u_h^{k+1}-u_h^k\right)}{\varphi_h}}{\tau} \,d\mu_{g_h} \\
        &+ \int_{M_h^\half} \,\lambda_h\IPA{\varphi_h}{N_h^\half} d\mu_{g_h} \\
        &- \int_{M_h^\half} \IPH{du_h^{k+1}}{d\varphi_h}d\mu_{g_h}.
\end{align*}
\end{remark}}

\green{
To further explain the heuristic behind Problem~\ref{prob:discpwillmore}, first recall that the identity $\divg\varphi = \IP{du}{d\varphi}$ holds in the continuous setting.  Therefore, as in \cite{dziuk2008} the divergence of $\varphi$ has been discretized as
\begin{equation*}
    \divgh\varphi_h = \IPH{du_h^\half}{d\varphi_h},
\end{equation*}
which improves numerical stability.  Moreover, the terms $\IP{dW}{d\varphi}$ and $\IP{du}{d\psi}$ have been discretized fully-implicitly in the interest of moving closer to a second-order time discretization.  On the other hand, the term $\IP{D(\varphi)du}{dW}$ is ill-behaved when not taken explicitly, so it has been discretized with respect to the old data $u_h^k, W_h^k$.  Note that, in any case, differentiation and integration are done with respect to the central surface $M_h^\half$, which greatly improves the results.}

\green{
Moreover, it is reasonable to consider (\ref{eq:areacon}) and (\ref{eq:volcon}) as a central discretization of the constraint equations (\ref{eq:forarea}) and (\ref{eq:forvol}).  To see this, first observe that 
\begin{align*}
    3\mathcal{A}(u) &= \int_M 3\,d\mu_g = \int_M \mathrm{div}\,u \,d\mu_g \\
    &=\int_M \left(\divg u + \IPA{N}{\nabla_N u}\right) d\mu_g = \int_M \left(1+ |du|^2\right) d\mu_g,
\end{align*}
where the second line uses the tangential/normal decomposition of the ambient divergence operator and the fact that $\nabla u = I$ is the identity matrix. Therefore, since area preservation can be enforced by requiring zero change between the areas of $M_h^{k+1}$ and $M_h^k$, a reasonable condition for area-preservation is
\begin{align*}
    0 &= 2\mathcal{A}\left(u_h^{k+1}\right) - 2\mathcal{A}\left(u_h^k\right) = \int_{M_h^{k+1}}\left\lvert du^{k+1}_h\right\rvert^2 d\mu_{g_h} - \int_{M_h^k} \left\lvert du^k_h\right\rvert^2 d\mu_{g_h} \\
    &\approx \int_{M_h^\half} \IPH{du_h^\half}{\left(du_h^{k+1}-du_h^k\right)} d\mu_{g_h},
\end{align*}
which is (\ref{eq:areacon}).  Similarly, the condition for volume preservation becomes 
\begin{align*}
    0 &= 3\mathcal{V}\left(u^{k+1}_h\right) - 3\mathcal{V}\left(u^k_h\right) \\
    &= \int_{M_h^{k+1}} \IPA{u^{k+1}_h}{N_h^{k+1}} d\mu_{g_h} - \int_{M_h^k} \IPA{u^k_h}{N_h^{k}} d\mu_{g_h} \\
    &\approx \int_{M_h^\half} \IPA{\left(u_h^{k+1} - u_h^k\right)}{N_h^\half} d\mu_{g_h},
\end{align*}
which is expression (\ref{eq:volcon}).}

\green{
Of course, the conformal penalty regularization Problem~\ref{prob:meshreg} can also be discretized in a similar fashion.  Though this problem is not necessarily time-dependent, it is advantageous to treat it somewhat implicitly so that the tangent bundle of the regularized mesh ``fits together'' in a smoother fashion.  Indeed, there is no reason to expect that sliding the mesh points of the initial surface along individual tangent planes will produce a new distribution which is itself integrable, so a semi-implicit discretization tends to produce better results in this case.  To that end, let $u_h^{k+1}, \hat{u}_h^{k+1}$ be the old resp. new immersions of $M_h^{k+1}$ as in Problem~\ref{prob:meshreg}, and let $N_h^{k+1}$ resp. $\hat{N}_h^{k+1}$ be their respective normal vector fields. Moreover, let $\tilde{N}_h^{k+1} := (1/2)\left(N_h^{k+1} + \hat{N}_h^{k+1}\right)$ denote the ``central'' normal field.  The discretized mesh regularization procedure then proceeds as follows.
\begin{problem}[Discrete conformal penalty regularization]\label{prob:discreteconf}
Let $\varepsilon>0$ be fixed, let $\hat{u},u,N,\rho$ be as in Problem~\ref{prob:meshreg}, and let $\tilde{N}$ be as above. Given $u_h^{k+1},N_h^{k+1}$, solving the discrete conformal penalty regularization problem means finding functions $\hat{u}_h^{k+1}, \rho_h$ which satisfy the system 
\begin{align*}
    0 &= \int_{M_h^{k+1}} \rho_h \IPA{\varphi_h}{\tilde{N}_h^{k+1}} d\mu_{g_h} + \int_{M_h^{k+1}} \IPH{\hat{Q}_h^{k+1}}{d\varphi_h} d\mu_{g_h}, \\
    0 &= \int_{M_h^{k+1}} \psi_h \IPA{\left(\hat{u}_h^{k+1} - u_h^{k+1} \right)}{\tilde{N}_h^{k+1}} d\mu_{g_h} + \varepsilon\int_{M_h^{k+1}}\psi_h\,\rho_h \,d\mu_{g_h},
\end{align*}
for all $\varphi_h,\psi_h \in S_h$ and  where $\IPH{\hat{Q}_h^{k+1}}{d\varphi_h}$ refers to the discretization on the known surface $M_h^{k+1}$ of the analogous quantity in \greener{Problem \ref{prob:meshreg}),} which involves components of the known normal $N_h^{k+1}$ and derivatives of the unknown immersion $\hat{u}_h^{k+1}$, computed with respect to $M_h^{k+1}$.
\end{problem}
Note that Problem~\ref{prob:discreteconf} is nonlinear, but only because the normal vector field arising from the surface preservation constraint has been taken centrally.  Therefore, as mentioned in Remark~\ref{rem:linvsnonlin}, the conformal penalty regularization procedure can be easily modified by replacing $\tilde{N}_h^{k+1}$ with $N_h^{k+1}$, yielding a linear system of equations.  This provides a tradeoff between mesh quality and computational time, as illustrated in Figure~\ref{fig:moocompare}.}

\green{
Now that the relevant continuous problems have been discretized, it is appropriate to give the full algorithm for running the p-Willmore flow with conformal penalty.  First, recall that when given an immersion $u_h^k$ at time step $k$ it is necessary to compute the curvature data $Y_h^k,W_h^k$.  This is accomplished through the solution of two consecutive linear systems:
\begin{align}
   0 &= \int_{M_h^k} \IPA{Y_h^k}{\psi_h}d\mu_{g_h} + \int_{M_h^k} \IPH{du_h^k}{d\psi_h} d\mu_{g_h}, \label{eq:forYinit} \\
   0 &=\int_{M_h^k} \IPA{\left(W_h^k - \left|Y_h^k\right|^{p-2}Y_h^k\right)}{\xi_h} d\mu_{g_h} \label{eq:forWinit}.
\end{align}
Note that, although the weighted mean curvature vector is defined pointwise as $W = |Y|^{p-2}Y$, equation (\ref{eq:forWinit}) computes $W$ in a weak sense.  The reason for this is because the solution to (\ref{eq:forWinit}) is often more smoothly distributed across the surface than the interpolated quantity stored pointwise at vertices, which leads to better numerical behavior during simulation.}

\begin{algorithm}
\caption{p-Willmore flow with conformal penalty}
\label{alg:main}
\green{
\begin{algorithmic}
\REQUIRE Closed, oriented surface immersion $u_h^0: M_h^0 \to \mathbb{R}^3$; real numbers $\varepsilon, \tau > 0$, integer $k_{\mathrm{max}} \geq 1$.
\WHILE{$0\leq k \leq k_{\mathrm{max}}$}
    \STATE Solve (\ref{eq:forYinit}) for $Y_h^k$
    \STATE Solve (\ref{eq:forWinit}) for $W_h^k$
    \STATE Solve Problem~\ref{prob:discpwillmore} for  $u_h^{k+1},Y_h^{k+1},W_h^{k+1},\lambda_h,\greener{\gamma_h}$
    \STATE Solve Problem~\ref{prob:discreteconf} for $\hat{u}_h^{k+1},\rho_h$
    \STATE $u_h^{k+1} = \hat{u}_h^{k+1}$
    \STATE $k = k + 1$
\ENDWHILE
\end{algorithmic}}
\end{algorithm}

\begin{figure}
\begin{center}
    \begin{minipage}[c]{0.45\textwidth}
        \includegraphics[width=\textwidth]{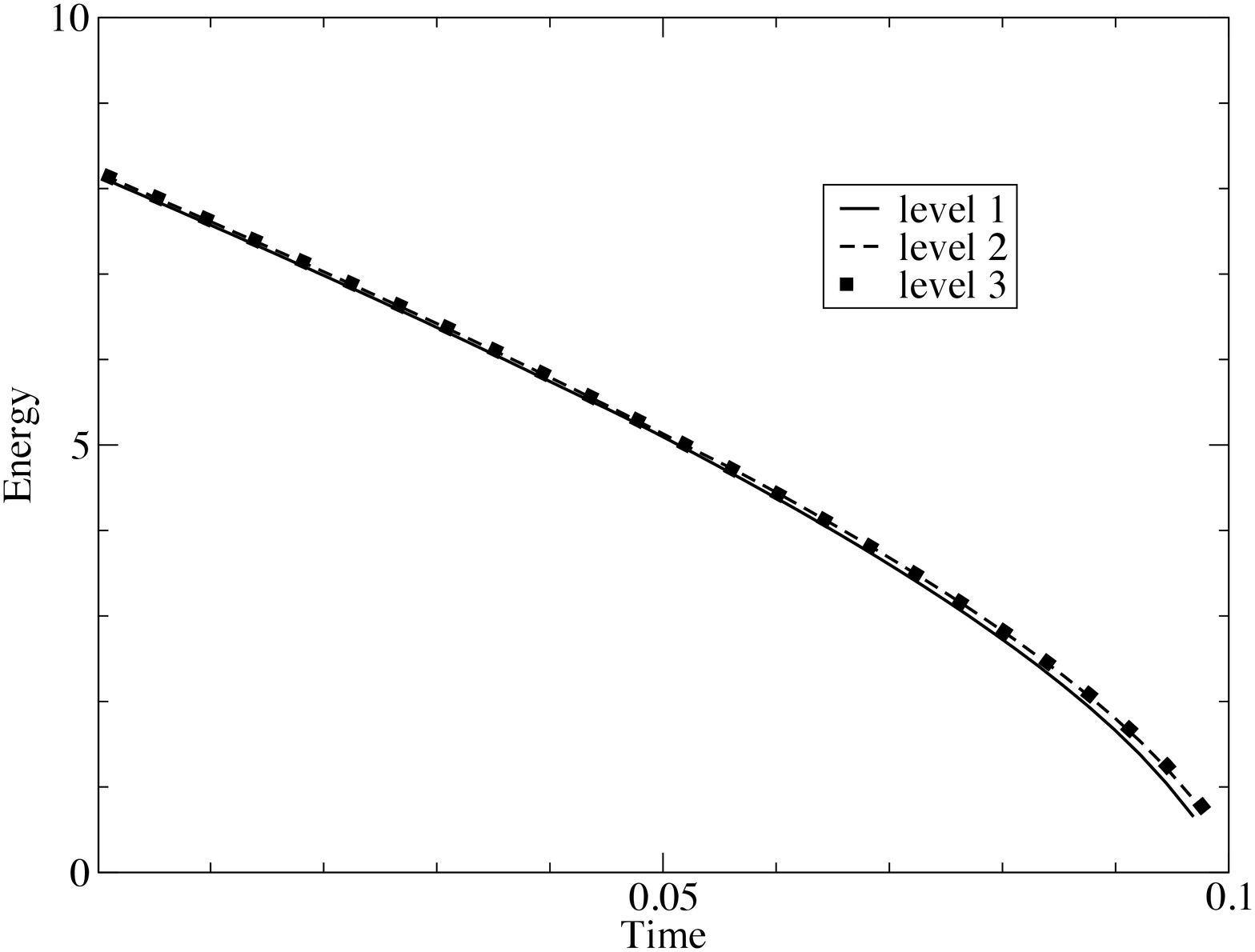}
    \end{minipage}
    \\
    \begin{minipage}[c]{0.45\textwidth}
        \includegraphics[width=\textwidth]{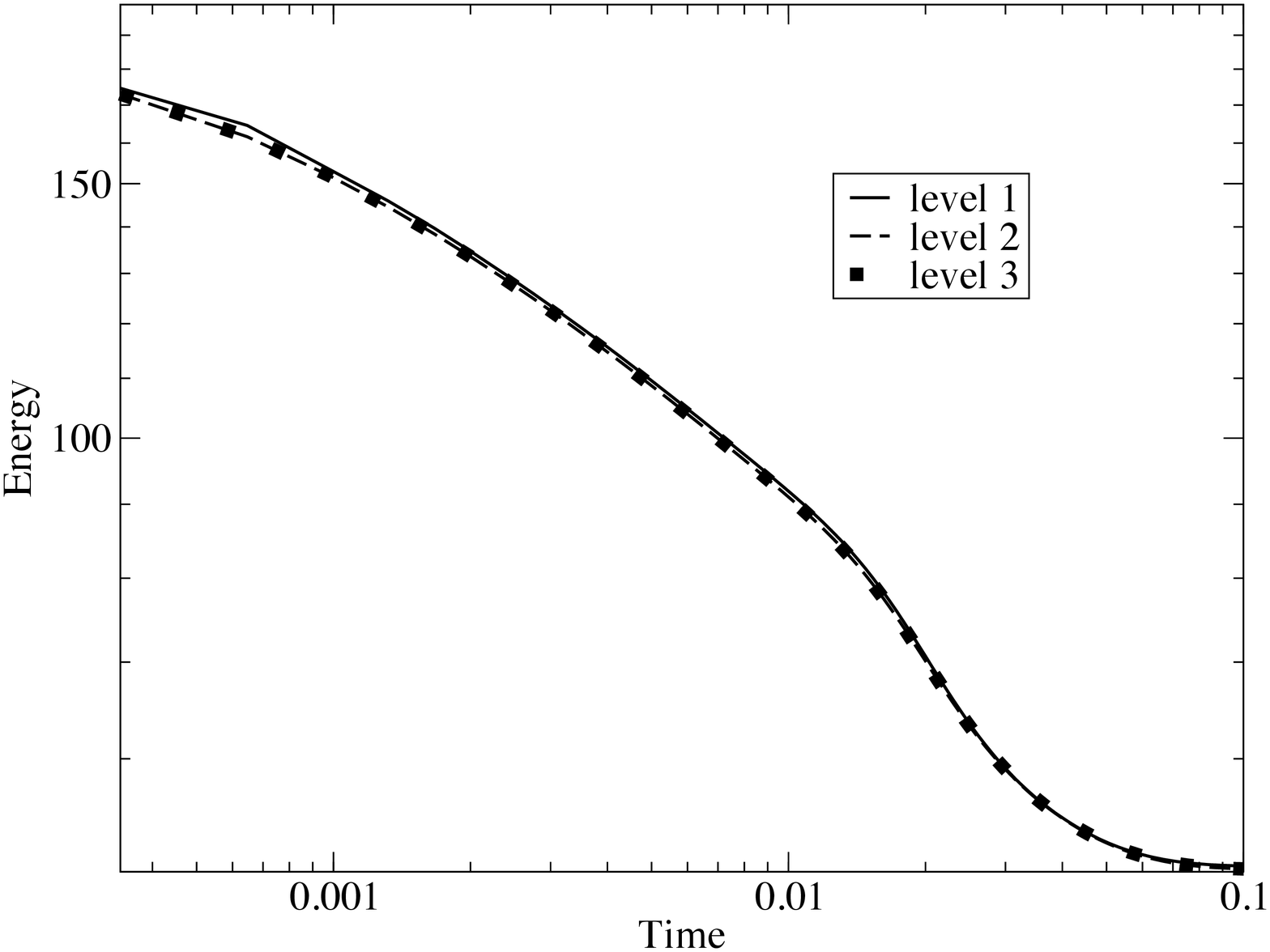}
    \end{minipage}
    \\
    \begin{minipage}[c]{0.45\textwidth}
        \includegraphics[width=\textwidth]{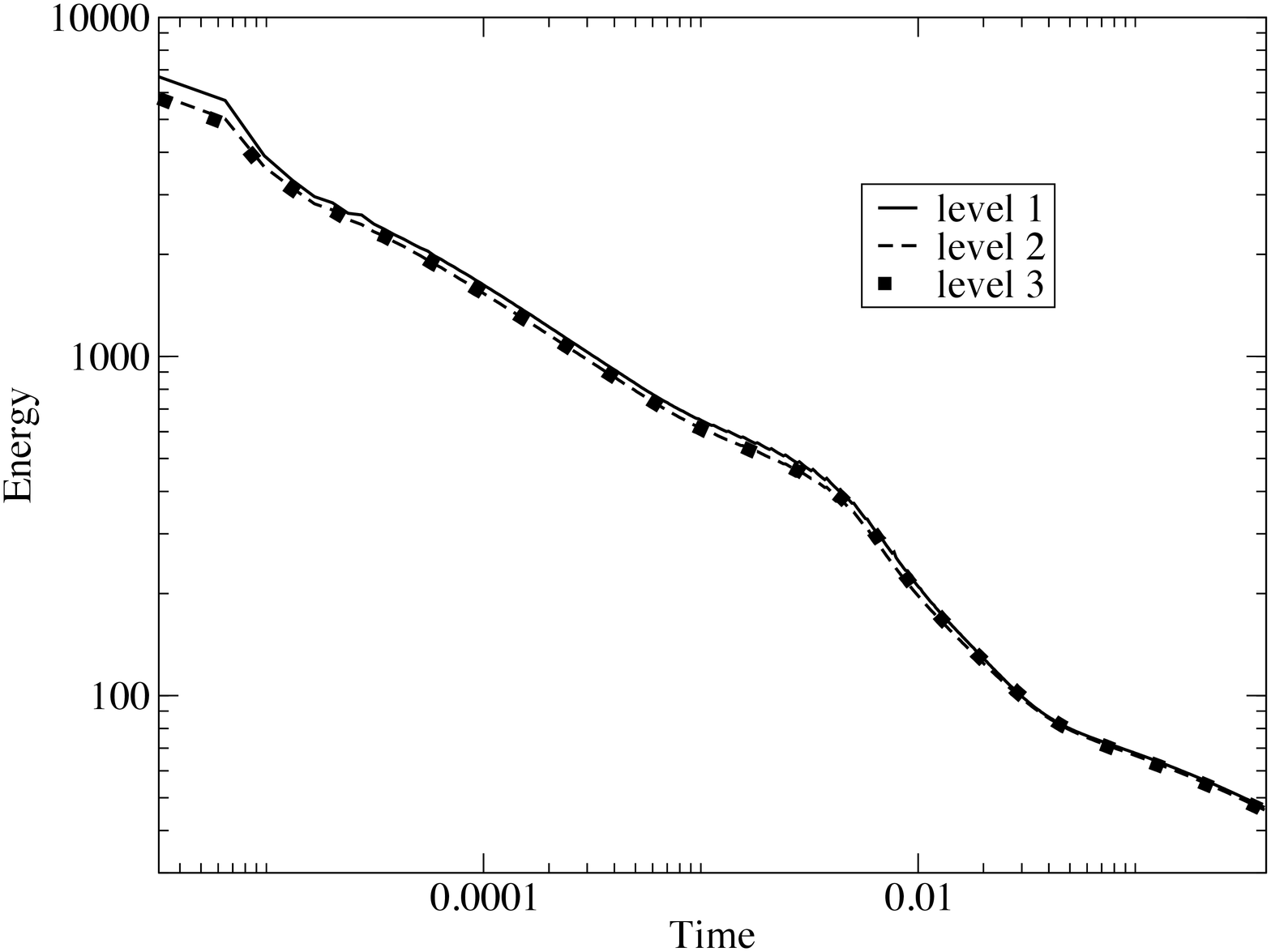}
    \end{minipage}
\end{center}
\caption{Algorithm~\ref{alg:main} applied to various refinements of the letter ``C'' from Figure~\ref{fig:C} when $p=0 \,\mathrm{(top)}, 2 \,\mathrm{(mid)}, 4 \,\mathrm{(bot)}$.  Time domain displayed logarithmically for $p=2,4$.  Note that the energy decreases in every case, and appears to converge with mesh refinement.}  
\label{fig:stability}
\end{figure}

\green{
Algorithm~\ref{alg:main} is the full procedure developed here for studying the computational p-Willmore flow of closed surfaces.  Though Problems~\ref{prob:discpwillmore} and \ref{prob:discreteconf} can certainly be used independently of each other, their combination as above provides a nice tool which leads to the variety of flow simulations seen presently.  Before discussing potential applications of this algorithm, it is important to discuss its stability.  Though precise analysis of the fully-discrete systems in Algorithm~\ref{alg:main} has not yet been done, it is easy to verify empirically that the energy-decreasing property of the continuous system (c.f. Theorem~\ref{thm:paramflow}) is preserved by the chosen discretization. An example of this is demonstrated in Figure~\ref{fig:stability}, which shows experimental results for the p-Willmore flow applied to the mesh in Figure~\ref{fig:C}.  As expected, the energy decreases monotonically in every case, suggesting that the fully-discrete flow is indeed numerically stable.  Moreover, notice that the conformally-invariant 2-Willmore flow levels off at $16\pi$ ($2^2$ times the theoretical minimum, \greener{c.f. Remark~\ref{rem:scaled}}), while the MCF and 4-Willmore flows decrease indefinitely.}

\green{
On the other hand, it is clear from experimentation that stability for Algorithm~\ref{alg:main} \greener{can hold} only in a conditional sense.  Even though the discrete scheme used is essentially implicit, a restriction on the temporal step size is necessary for reasonable results. This is expected, in general, since integration is performed on a surface (mesh) whose evolution is solution-dependent; a large time step will easily generate entanglements which cause the mesh to crumple and invert. In fact, this is precisely what motivated the regularization Problem~\ref{prob:discreteconf}.  Since an evolving surface may change its shape dramatically over time, intermediate regularization is necessary to prevent failure caused by numerical degeneration along the flow.}

\green{
Empirically, it is seen that $p=0,2$ are quite robust to changes in the temporal step-size, but values of $p>2$ produce simulations that are noticeably more sensitive.  As illustrated in Table~\ref{tab:imp}, the simulations for 4-Willmore flow tend to require a much smaller step-size and a much larger amount of iterations to converge.  However, it is also observed that the p-Willmore flow with conformal penalty is relatively independent of mesh resolution; for a fixed $p$ and a temporal step-size $\tau$ that is stable on a coarse mesh, $\tau$ appears to remain stable on any refinement of that mesh.  This desirable property is suspected to come from the regularization in Problem~\ref{prob:discreteconf}, which ensures that mesh elements do not become too heavily distorted during the p-Willmore flow. \greener{Indeed, this is reflected in Figure~\ref{fig:brokeknot}, where the regularization prevents mesh elements from degenerating, even as the area and volume are both constrained.} The next section will discuss some specifics regarding the implementation of Algorithm~\ref{alg:main}, as well as how it can be applied to navigate a common problem in computer graphics.}



\section{Implementation and application}
\green{
\subsection{Implementation}
In brief, the nonlinear systems in Problems~\ref{prob:discpwillmore} and \ref{prob:discreteconf} are solved through a two-iteration Newton scheme, using a $7^{th}$-order tensor product quadrature rule to evaluate the relevant integrals.  To elaborate, consider the process of solving Problem~\ref{prob:discpwillmore}, and denote $\mathbf{v}_h = (u_h,Y_h, W_h,\lambda_h,\gamma_h)$. Then, a solution to the discrete p-Willmore flow system can be formally represented as a solution to the equation 
\begin{align*}
\mathcal{R}\left(\mathbf{v}_h\right) = 0,
\end{align*}
where $\mathcal{R}$ is an operator representing the nonlinear residual of the p-Willmore system.  Let $\mathbf{v}_h^0$ be a trial solution, and let $\mathcal{J}\left(\mathbf{v}_h\right)$ represent the tangent operator (or Jacobian) of $\mathcal{R}$ in $\mathbf{v}_h$, evaluated through
\[\mathcal{J}\left(\mathbf{v}_h^i\right)=\frac{\partial \mathcal{R}}{{\partial \mathbf{v}_h}}\left(\mathbf{v}_h^i\right).\]
Then, Newton iteration involves the procedure 
\begin{equation*}
\mathbf{v}_h^i = \mathbf{v}_h^{i-1} - \mathcal{J}^{-1}\left(\mathbf{v}_h^{i-1}\right) \mathcal{R}\left(\mathbf{v}_h^{i-1}\right) \quad \mbox{ for } \quad i \ge 1,
\end{equation*}
which is typically repeated until the residual quantity $\| \mathcal{R}(\mathbf{v}_h^i) \|$ drops below a predefined tolerance value.  Newton iteration is known to exhibit quadratic convergence provided that the initial guess $v_h^0$ is sufficiently close to the true solution. In the case of Algorithm~\ref{alg:main}, only two iterations of each nonlinear system in Problems~\ref{prob:discpwillmore} and \ref{prob:discreteconf} are performed at each time step, which is sufficient to produce a small residual and negligible change between successive solutions.}


\green{
Moreover, note that symbolic differentiation of $\mathcal{J}$ is cumbersome for the particular problems considered here, due to the presence of integrals evaluated on the evolving surface $M_h$. Though approximate evaluation of $\mathcal{J}$ is of course possible (by e.g. neglecting the motion-dependent nature of some terms or using approximate differentiation methods), it is advantageous to compute the exact Jacobian so that less error is introduced at each step.  This is accomplished presently with fast reverse automatic differentiation as described in \cite{hogan2014}.  Automatic differentiation techniques use the chain rule along with backpropagation to numerically evaluate the derivatives of a specified function. In particular, since the derivative of a composite function involves a product of terms which are sequentially computable through elementary arithmetic operations, repeated application of the chain rule can be used to accurately evaluate derivatives of arbitrary order. The implementation here uses the Adept library, which enables algorithms written in C and C++ to be automatically differentiated with an operator overloading strategy. In addition, the solution of all linear systems necessary for Algorithm~\ref{alg:main} is performed using the direct solver found in the MUMPS library \cite{MUMPS:1,MUMPS:2}.}

\green{
It is worth mentioning that a viable alternative to this approach would be to pull every integral expression in Problems~\ref{prob:discpwillmore} and \ref{prob:discreteconf} back to an ``original'' parametrization domain $U \times \{0\}$, which avoids differentiation on a moving surface.  While this greatly complicates the formulation, it has the advantage of allowing for the Jacobian $\mathcal{J}$ to be evaluated using purely symbolic differentiation. However, for Algorithm~\ref{alg:main} this approach is not at all necessary, and automatic differentiation is found to be optimal for producing good results.  Therefore, the present approach has been chosen for its higher clarity of formulation, as well as its relative ease of numerical implementation. }

\begin{figure}
\begin{center}
\begin{minipage}[c]{0.22\textwidth}
\includegraphics[width=\textwidth]{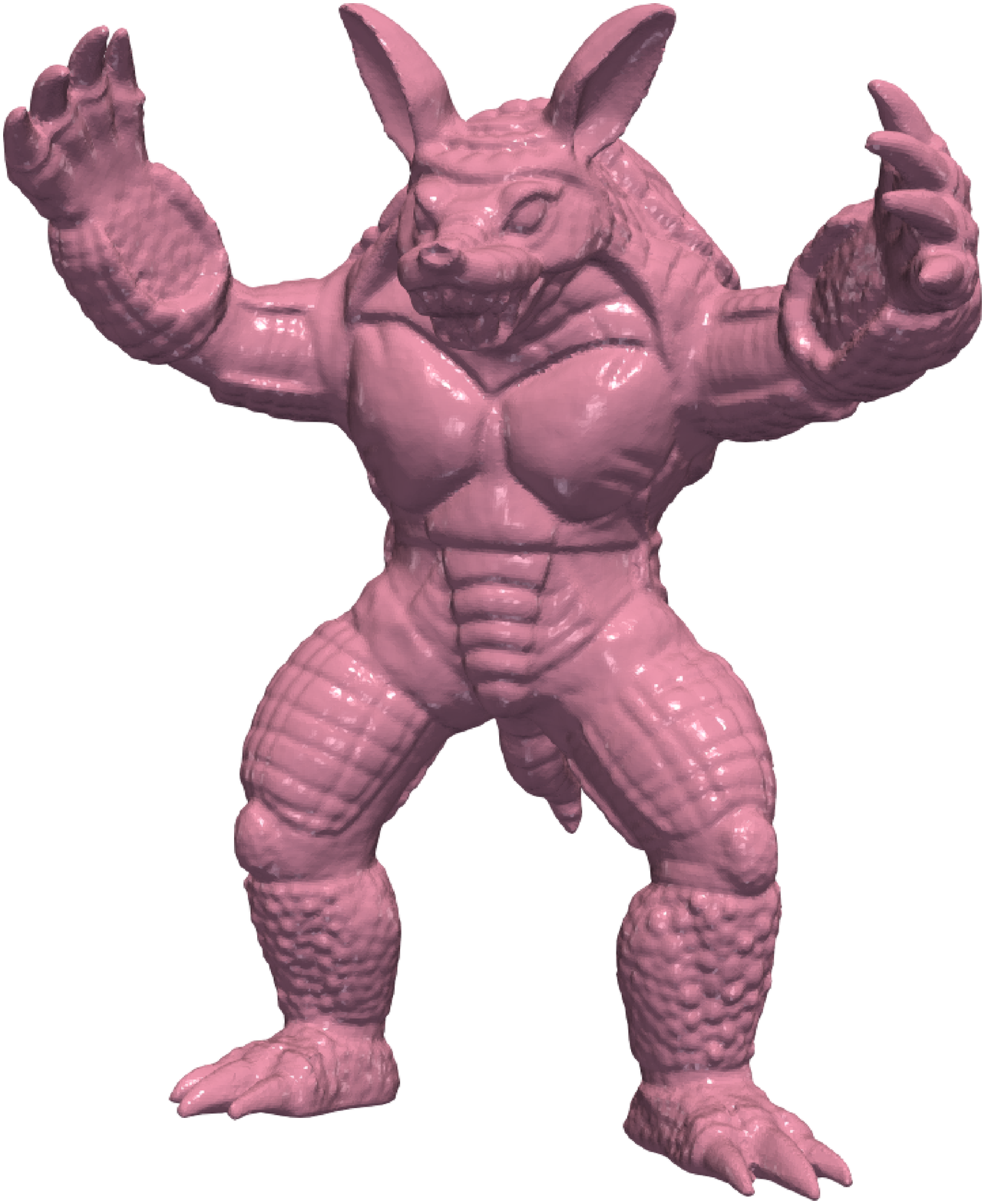}
\end{minipage}
\hspace{0.4pc}
\begin{minipage}[c]{0.22\textwidth}
\includegraphics[width=\textwidth]{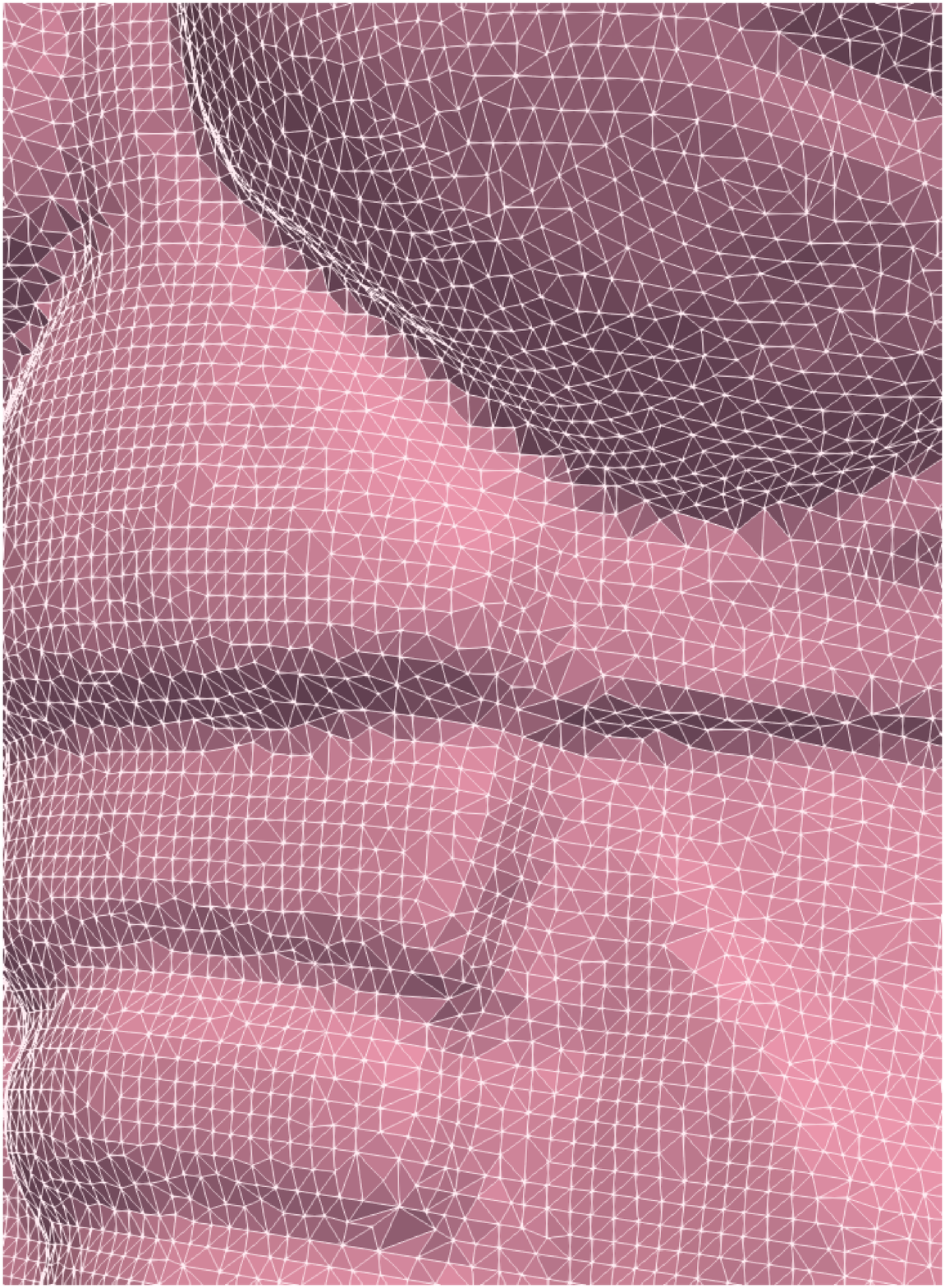}
\end{minipage}
\\
\begin{minipage}[c]{0.22\textwidth}
\includegraphics[width=\textwidth]{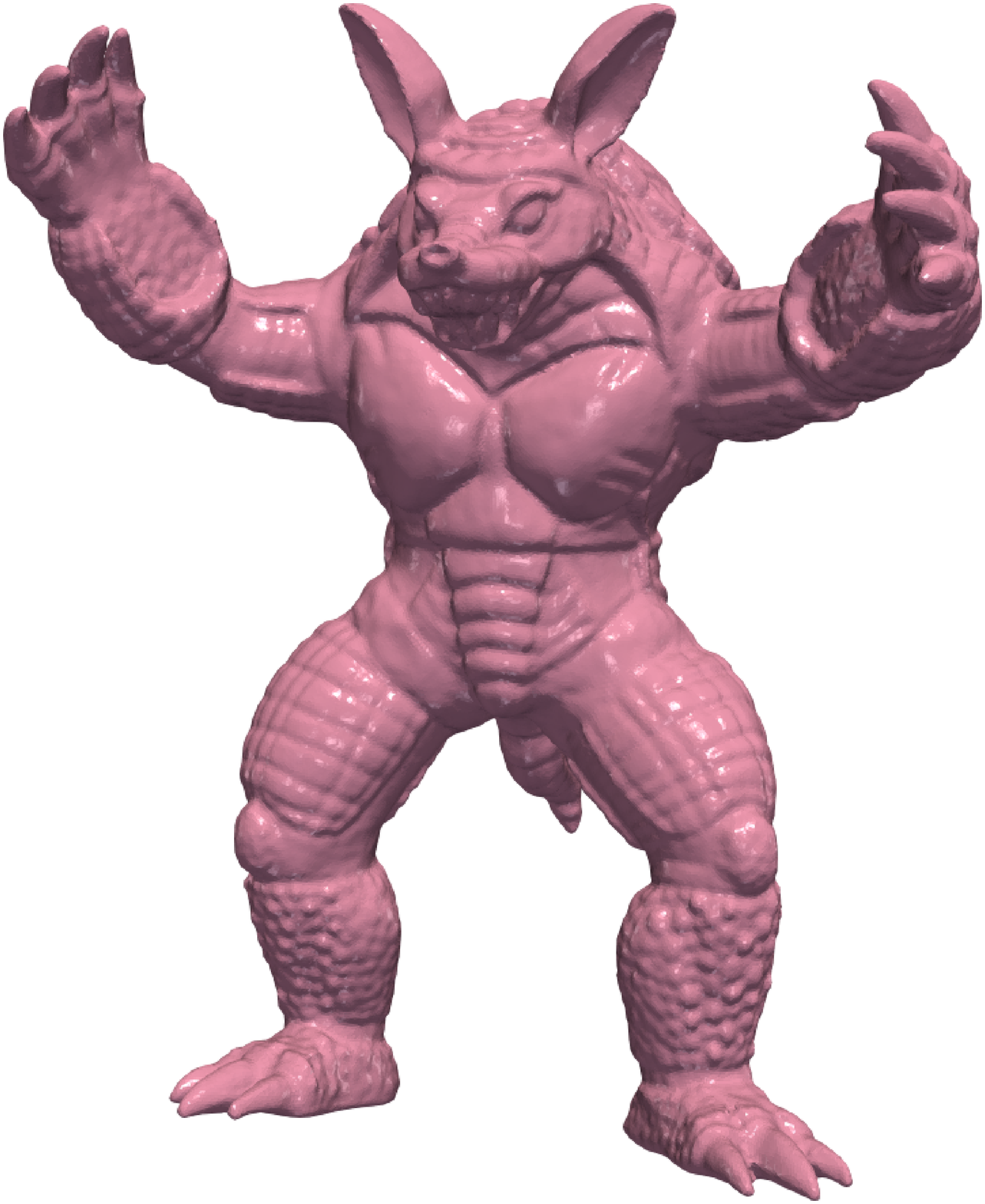}
\end{minipage}
\hspace{0.4pc}
\begin{minipage}[c]{0.22\textwidth}
\includegraphics[width=\textwidth]{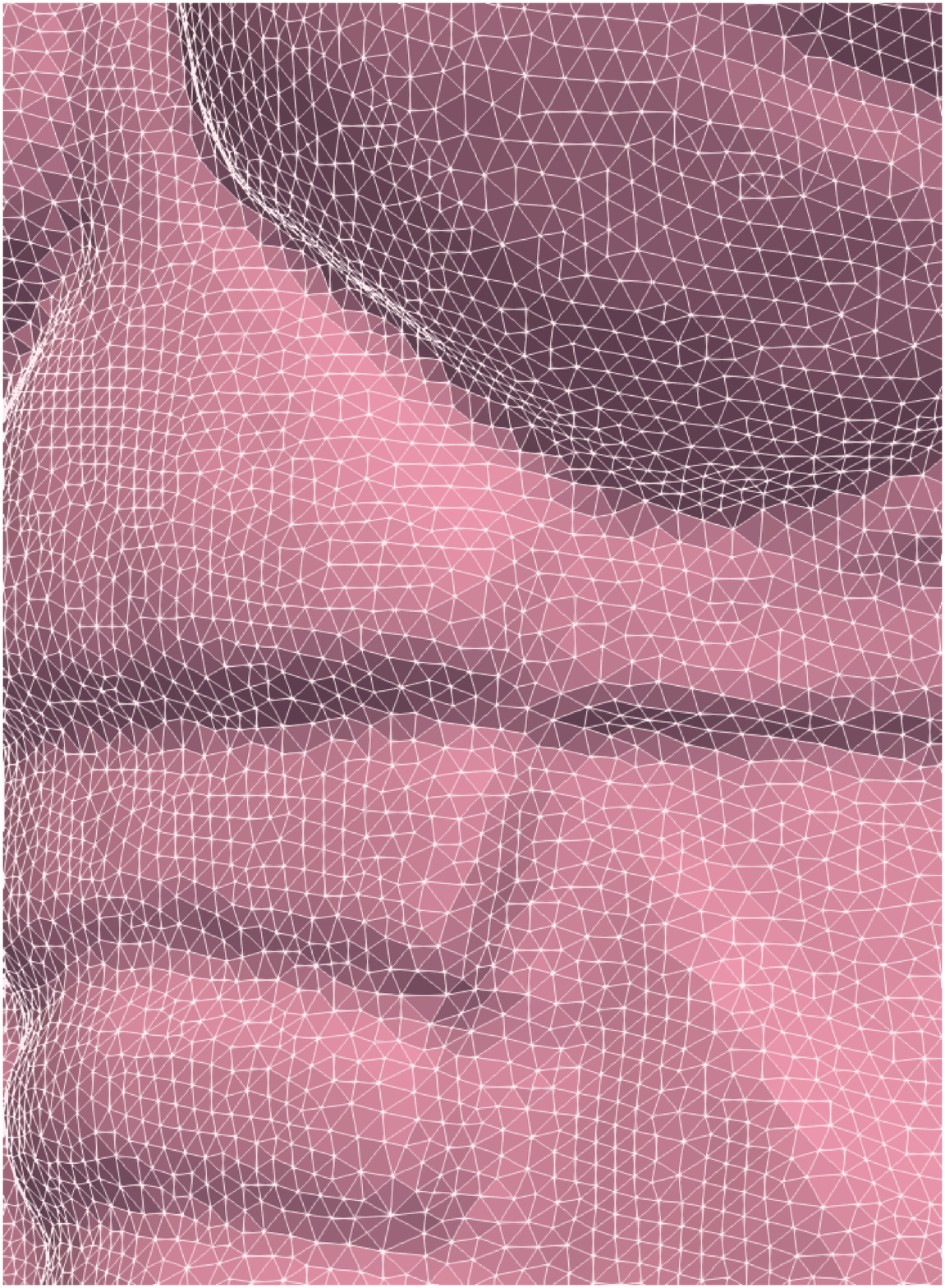}
\end{minipage}
\end{center}
\caption[Mesh Regularization of an Armadillo]{A cartoon armadillo with 346k triangles edited by our method.  Requires roughly 12 minutes of solver time on a 2.7GHz Intel Core i5 with 8GB of RAM.}
\label{fig:dillo}
\end{figure}

\green{
\subsection{Application: mesh editing}
Many algorithms in computer graphics are sensitive to the quality of their initial surface data, and (as seen with the p-Willmore flow) a poor mesh can frequently cause numerical failure independent of the actual geometry involved.  To add to the library of techniques which address this problem, consider the application of Algorithm~\ref{alg:main} with the goal of improving mesh quality.  It is seen that running a short p-Willmore flow followed by the conformal penalty regularization procedure will often produce a surface that is very close to the original, but with a better quality triangulation.  For example, Figure~\ref{fig:funcow} (picture 2) shows the result of one iteration of Algorithm~\ref{alg:main} with $p=2$ and a very small stepsize.  Notice that the change in surface geometry is quite small, while the mesh has been significantly improved.  Similarly, Figure~\ref{fig:dillo} shows the result after one \emph{linear} 2-Willmore iteration followed by two-step nonlinear conformal penalty regularization.  Here the original and remeshed surfaces can hardly be distinguished by eye, though the new triangulation is again much more regular.  Further, Figures~\ref{fig:moocompare} and \ref{fig:meshcomp} show the effects of conformal penalty regularization without any p-Willmore flow, which requires much less compute time.  In every case, the initial mesh is significantly improved with only slight changes to the surface geometry. }

\begin{figure*}
\begin{center}
    \begin{minipage}[c]{0.15\textwidth}
    MCF \\ (0-Willmore)
    \end{minipage}
    \hspace{-0.4pc}
    \begin{minipage}[c]{0.2\textwidth}
    \includegraphics[width=\textwidth]{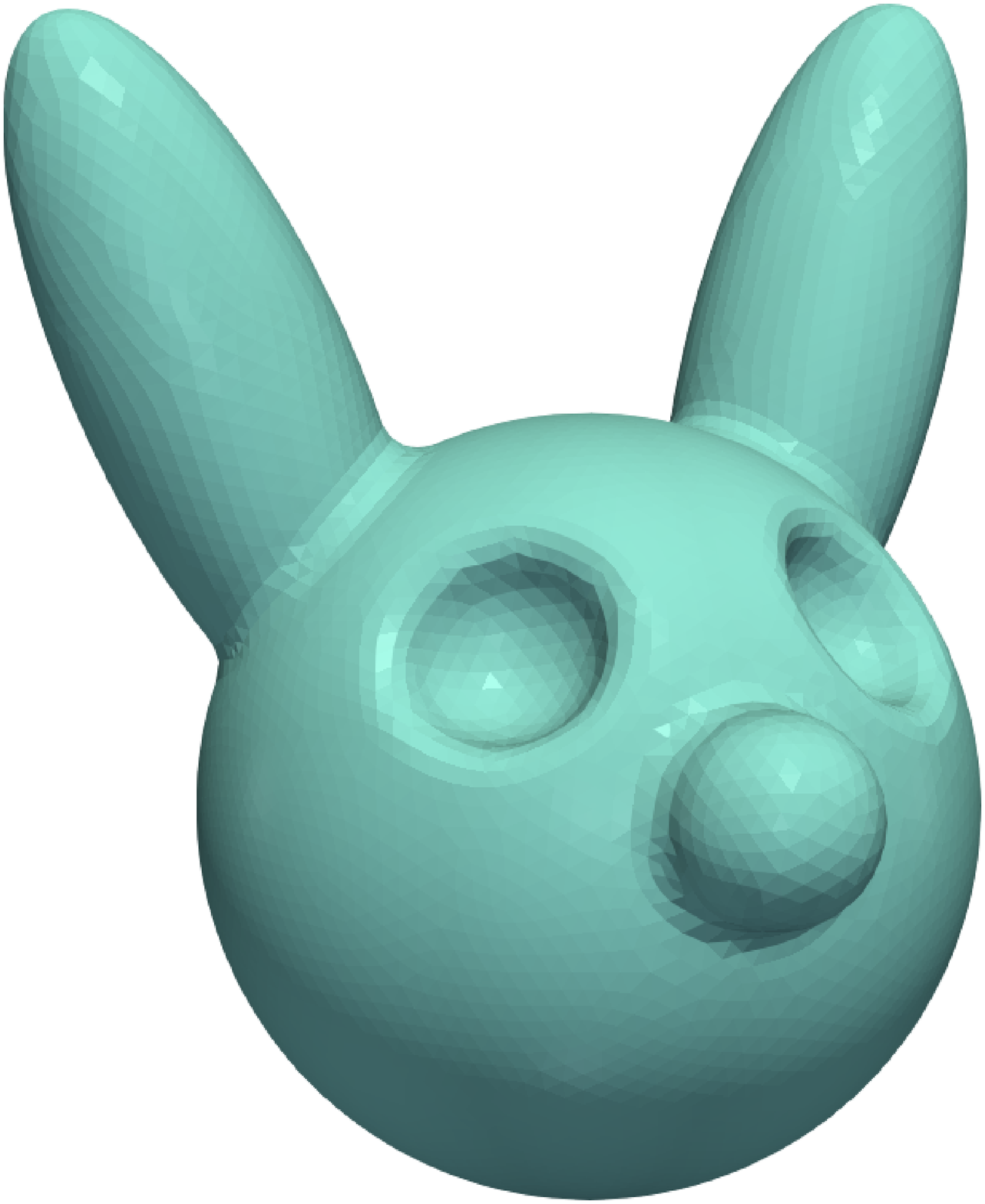}
    \end{minipage}
    \begin{minipage}[c]{0.20\textwidth}
    \includegraphics[width=\textwidth]{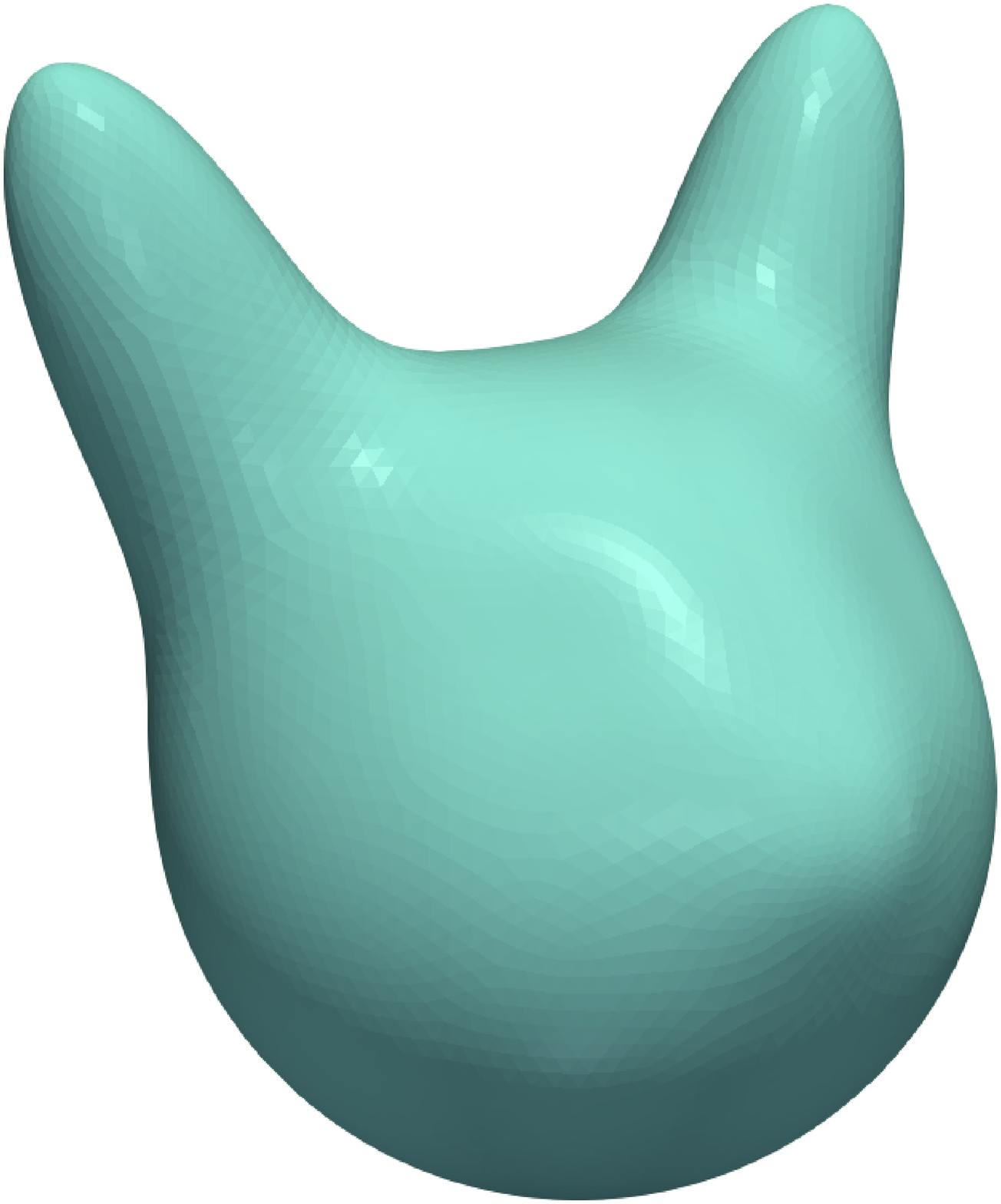}
    \end{minipage}
    \begin{minipage}[c]{0.20\textwidth}
    \includegraphics[width=\textwidth]{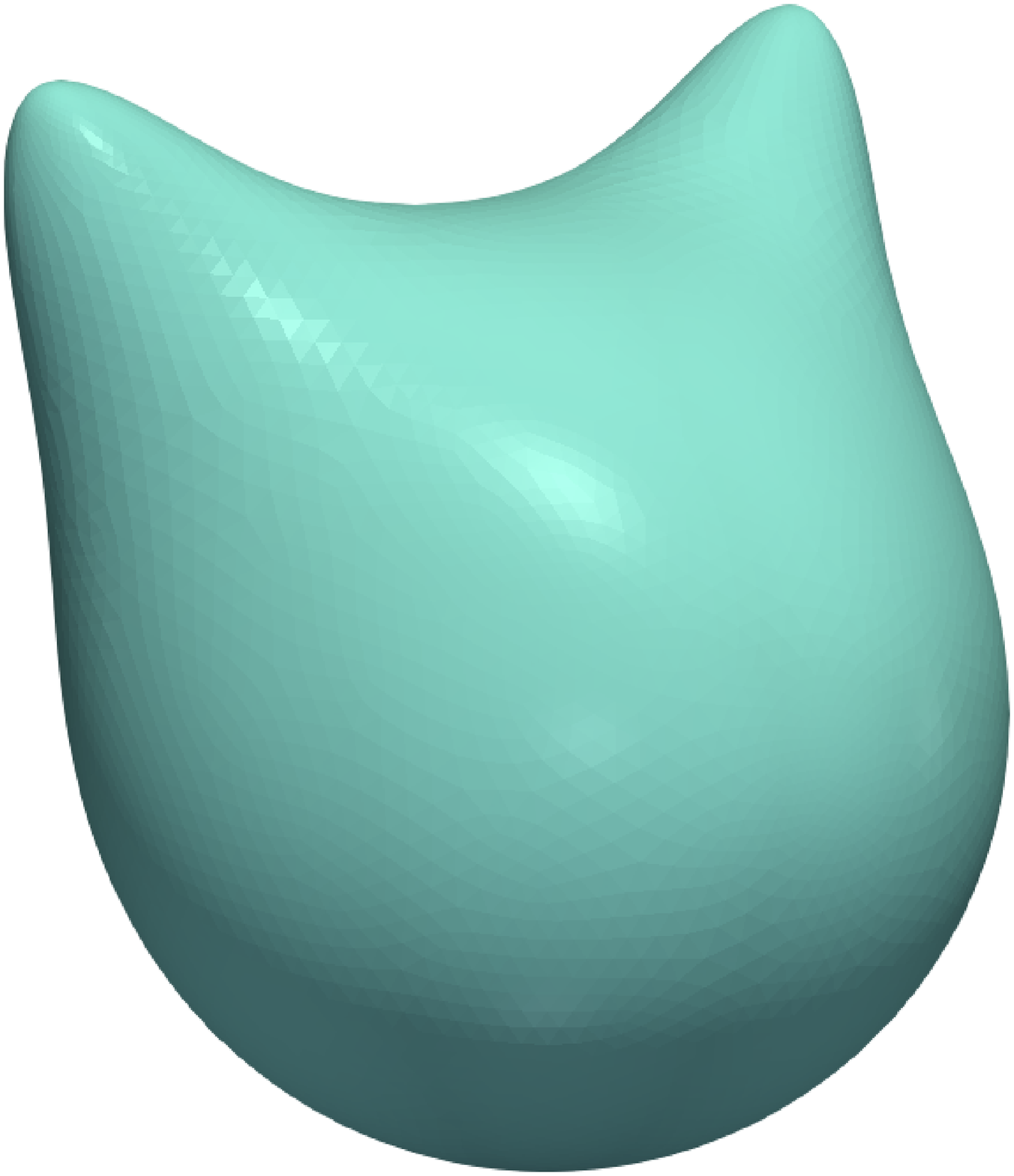}
    \end{minipage}
    \begin{minipage}[c]{0.20\textwidth}
    \includegraphics[width=\textwidth]{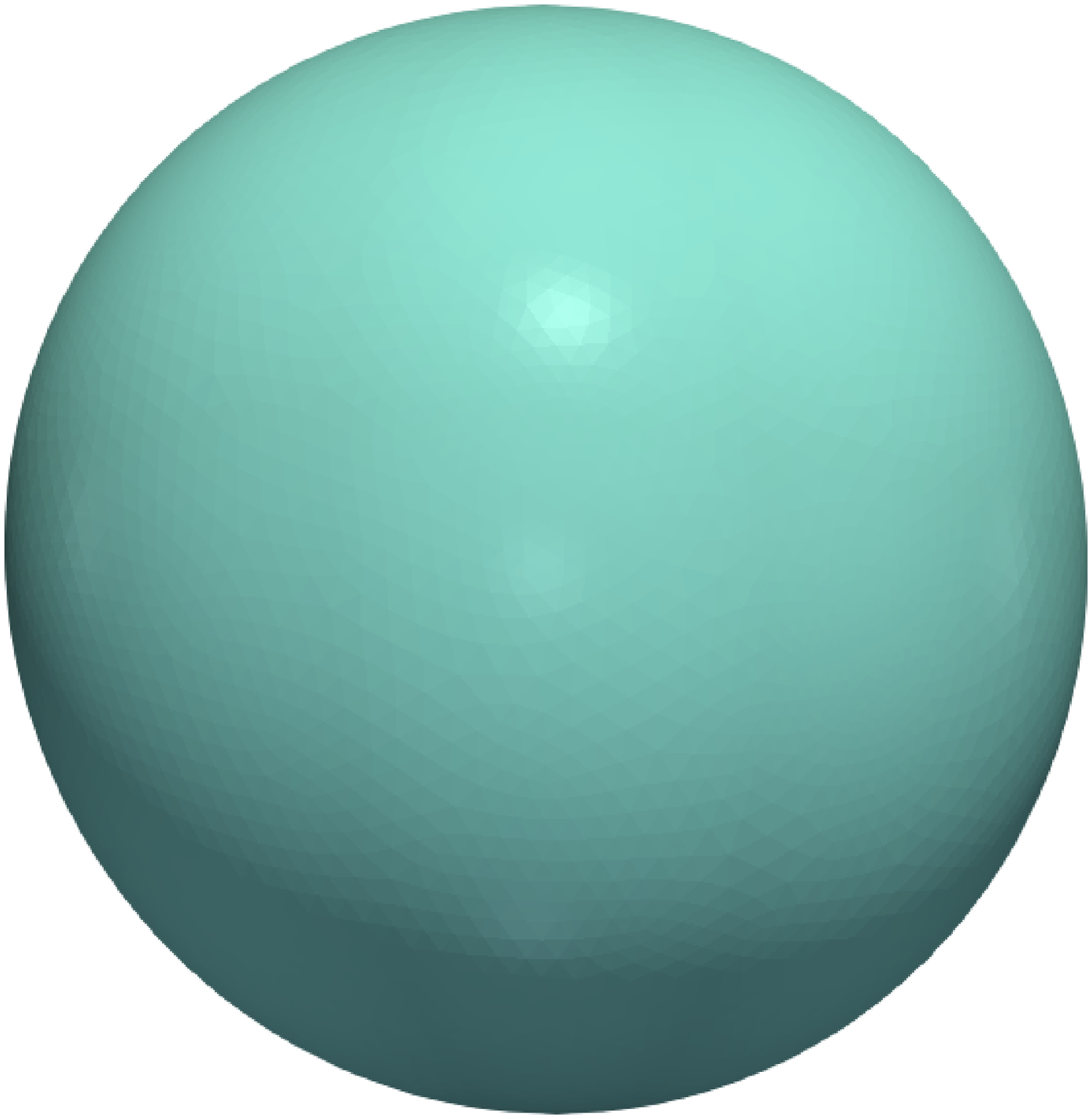}
    \end{minipage}
    \\
    \begin{minipage}[c]{0.15\textwidth}
    Willmore flow \\ (2-Willmore)
    \end{minipage}
    \hspace{-0.4pc}
    \begin{minipage}[c]{0.2\textwidth}
    \includegraphics[width=\textwidth]{figs/dog_start.eps}
    \end{minipage}
    \begin{minipage}[c]{0.2\textwidth}
    \includegraphics[width=\textwidth]{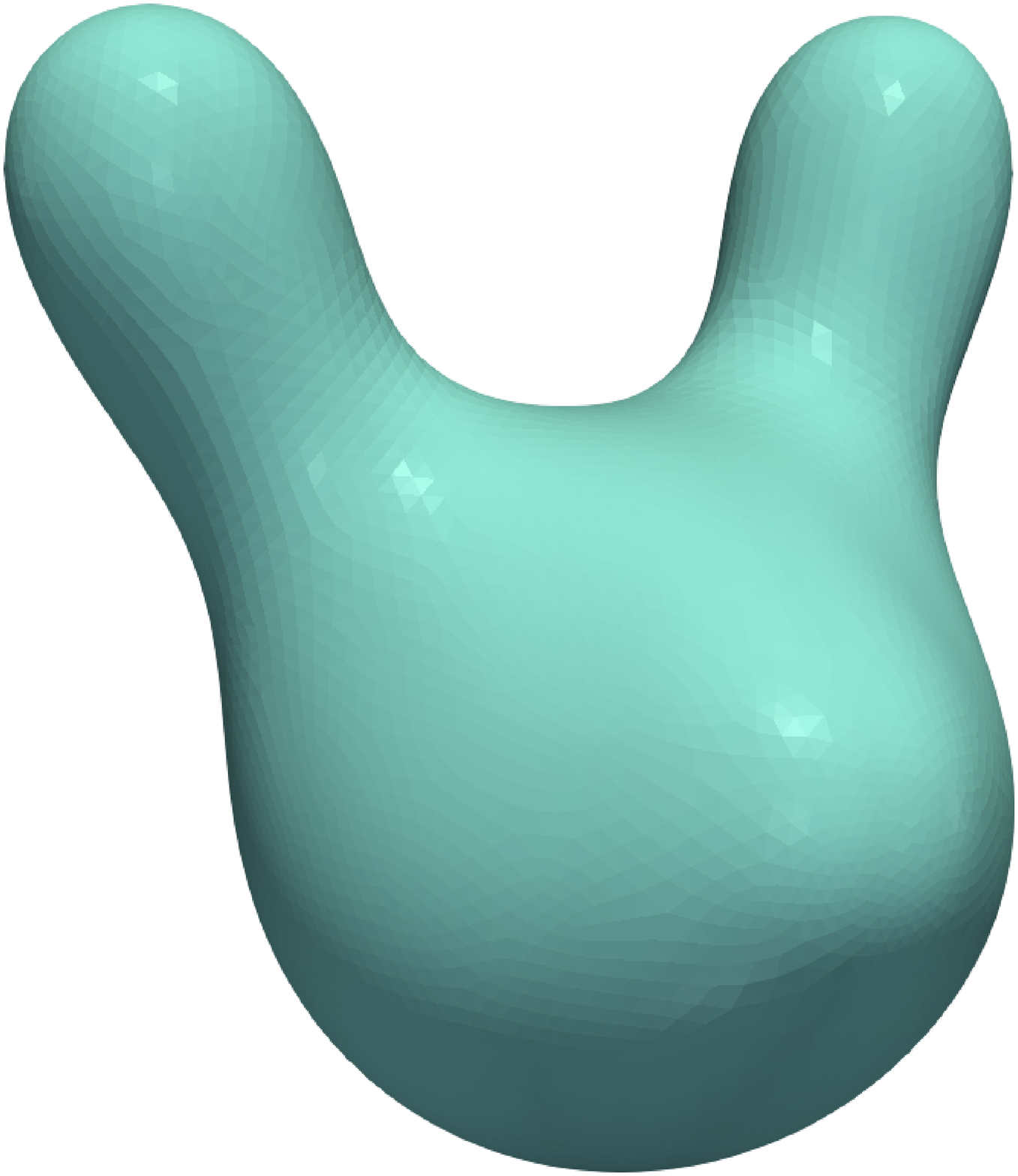}
    \end{minipage}
    \hspace{0.1pc}
    \begin{minipage}[c]{0.20\textwidth}
    \includegraphics[width=\textwidth]{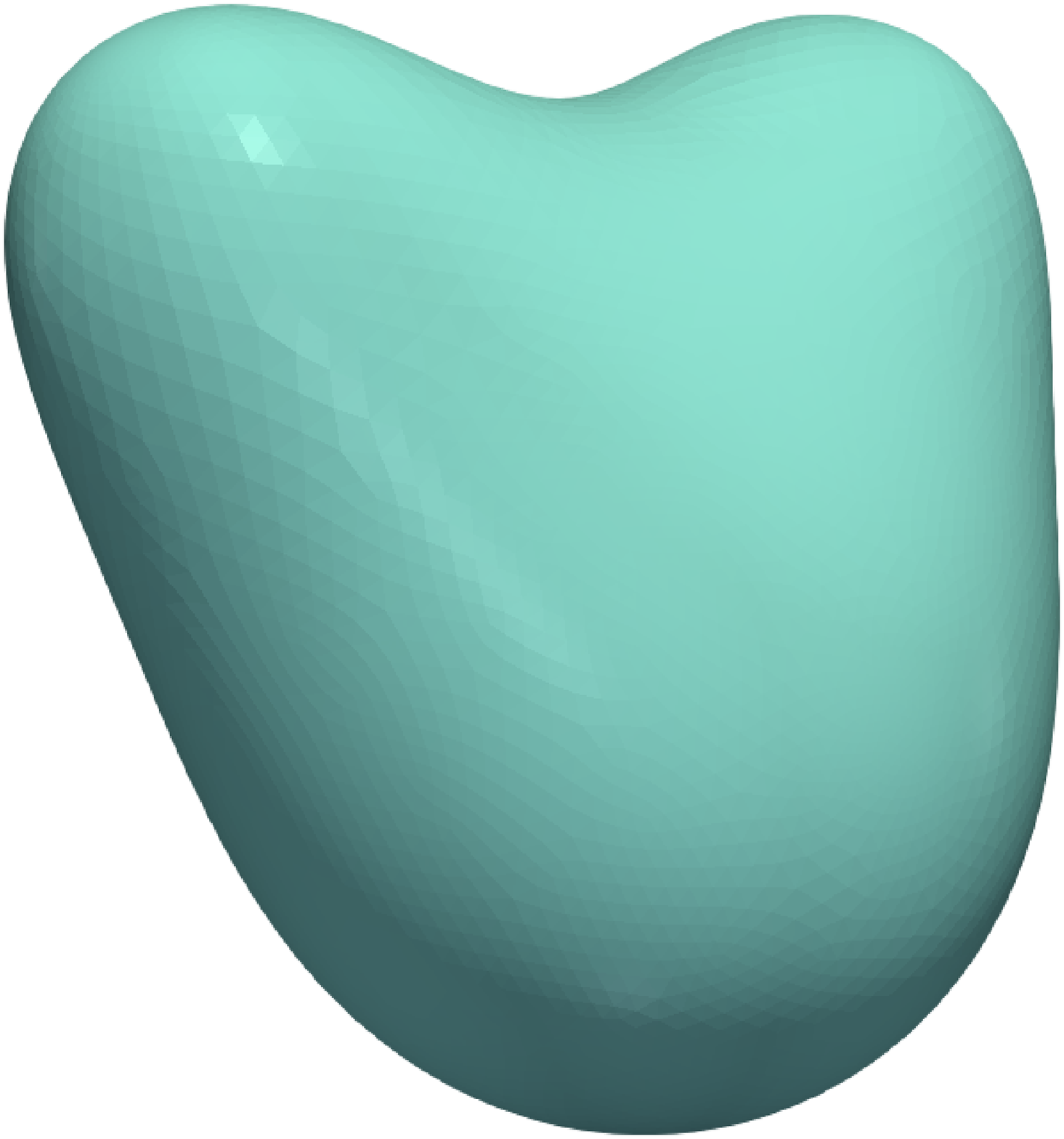}
    \end{minipage}
    \begin{minipage}[c]{0.20\textwidth}
    \includegraphics[width=\textwidth]{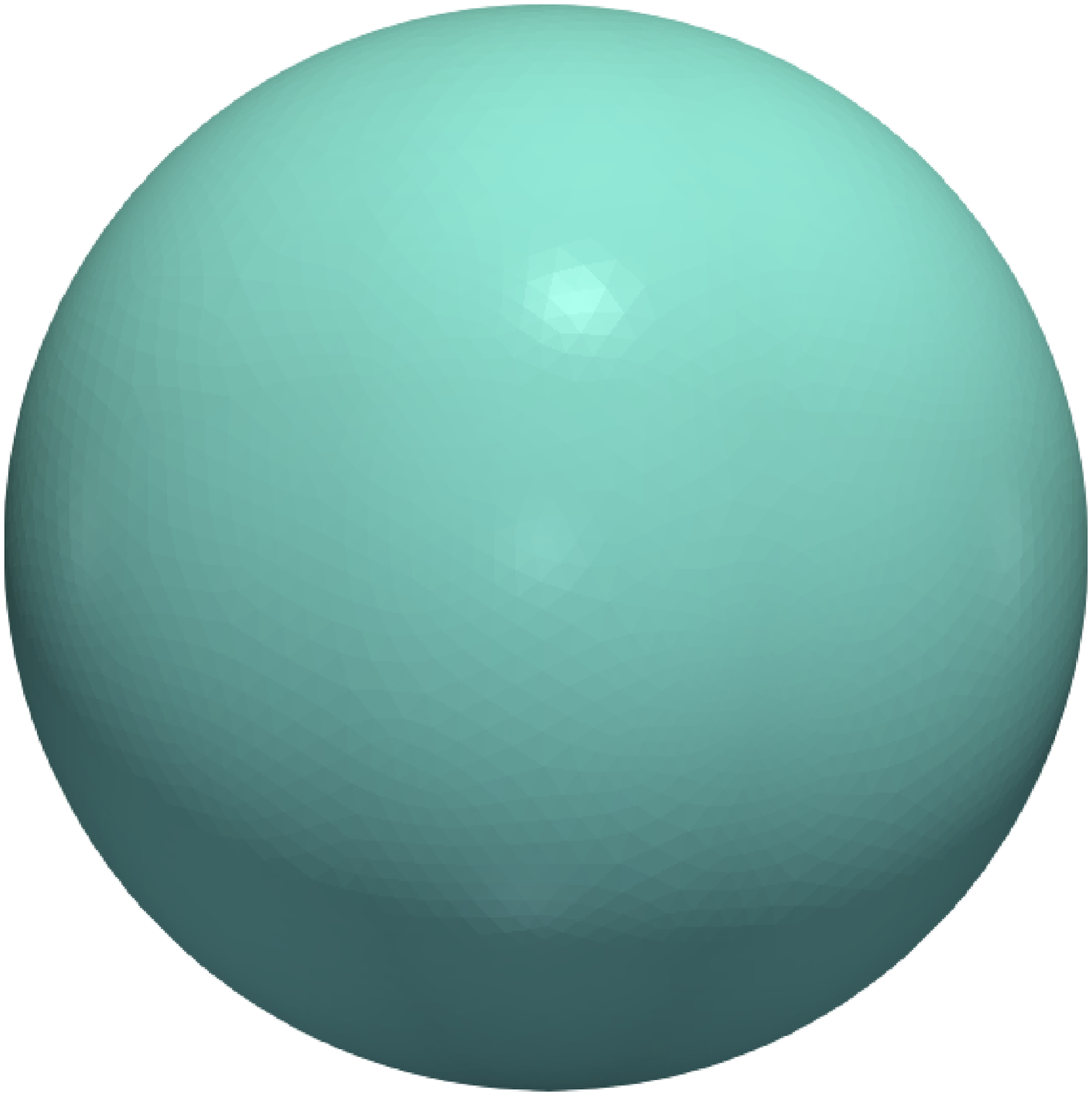}
    \end{minipage}
    \\
    \begin{minipage}[c]{0.15\textwidth}
    4-Willmore flow \\
    \end{minipage}
    \hspace{-0.4pc}
    \begin{minipage}[c]{0.2\textwidth}
    \includegraphics[width=\textwidth]{figs/dog_start.eps}
    \end{minipage}
    \begin{minipage}[c]{0.2\textwidth}
    \includegraphics[width=\textwidth]{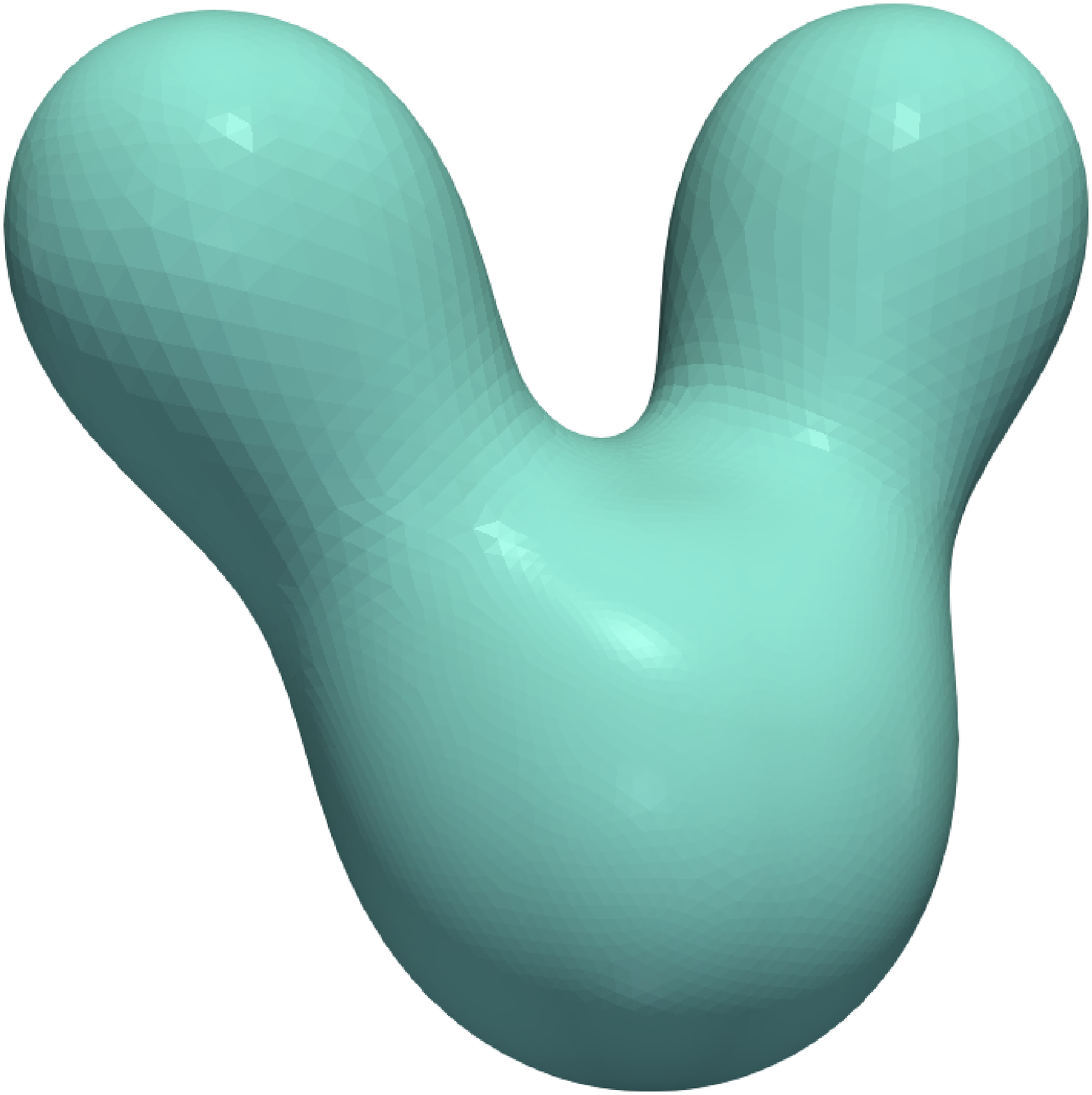}
    \end{minipage}
    \hspace{0.4pc}
    \begin{minipage}[c]{0.20\textwidth}
    \includegraphics[width=\textwidth]{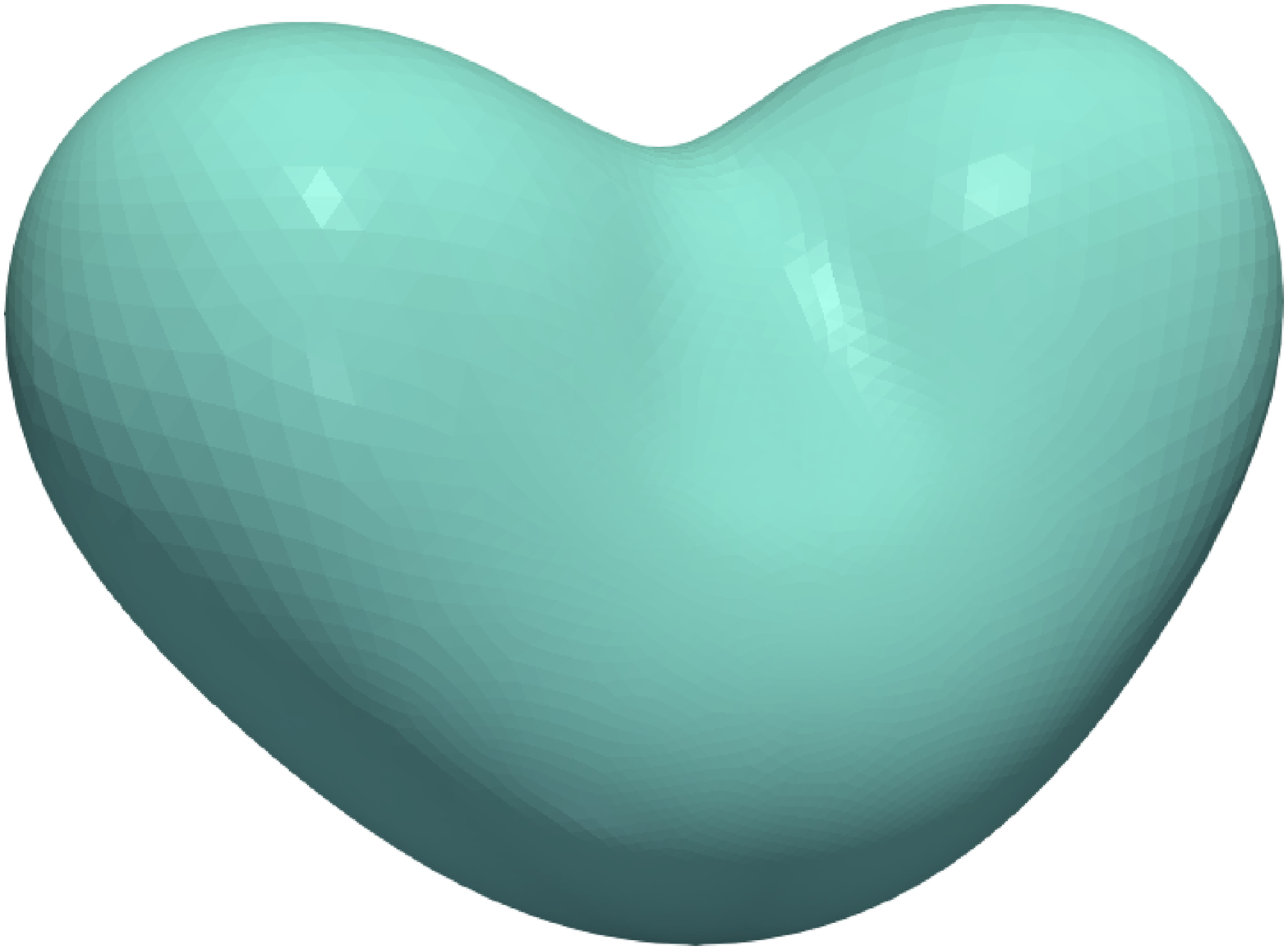}
    \end{minipage}
    \begin{minipage}[c]{0.20\textwidth}
    \includegraphics[width=\textwidth]{figs/p2after87.eps}
    \end{minipage}
\end{center}
\caption[Comparison of transient surfaces when $p=0,2,4$]{Volume-preserving p-Willmore flow with conformal penalty applied to the mesh in Figure~\ref{fig:weirddog} when $p = 0$ (top), 2 (mid), 4 (bottom).} 
\label{fig:doggo}
\end{figure*}

\green{
On the other hand, the p-Willmore flow may also find utility in computer animation, as it can be used to dramatically alter the geometry of an object in a prescribed way. In particular, detailed objects with sharp features will evolve under the p-Willmore flow to minima that are as round as possible, which could be desirable when modeling fluids.  Moreover, Figures~\ref{fig:C} and \ref{fig:doggo} show that the value of $p$ has a significant effect on the flow behavior, though this is not surprising.  Since the functional $\mathcal{W}^p$ measures the $p^{th}$ power of $|H|$, regions of high curvature are weighted increasingly with the value of $p$.  This is why regions of high curvature tend to ``round out'' faster when $p$ is large (c.f. Figure~\ref{fig:doggo}), which may be desirable if the goal is to evolve more prominent features before affecting others that are less pronounced.
}

\green{
\paragraph{\textbf{Conclusion:}}
The p-Willmore flow with conformal penalty provides a unified computational treatment of gradient flows which arise from functionals which depend exponentially on the unsigned mean curvature.  The algorithm presented here provides a new device for visualizing the p-Willmore flow of closed surfaces subject to natural constraints on surface area and enclosed volume.  Besides extending known methods for computing the Willmore flow, it is seen that the conformal penalty regularization procedure inherent in this algorithm allows for certain computational surfaces to evolve to minima that would otherwise be unreachable.  Moreover, this regularization can be applied to stationary immersions as well, significantly improving mesh quality with only minor changes to the surface itself.  Avenues for future work include a more rigorous study of the consistency, stability, and convergence of the p-Willmore flow under mesh refinement, as well as a computer implementation that is more robust to rough data.
}


\begin{table*}
\caption{Implementation details for the figures.  Note that $l^2 = (1/3)^2\left((x_{max} - x_{min})^2 + (y_{max}-y_{min})^2 + (z_{max}-z_{min})^2\right)$ and the scale factor $s$ multiplies $\tau$ at each iteration, i.e. $\tau^{k+1} = s \tau^k$ (used to speed up the evolution).  All simulations performed on a 2.7 GHz Intel Core i5 with 8GB of RAM.}
\resizebox{.99\textwidth}{!}{%
\begin{tabular}{|l|l|l|l|l|l|l|l|l|l|l|l|}
\hline
Figure & Geometry & Elements & \begin{tabular}[c]{@{}l@{}}Characteristic\\ Length $l$\end{tabular} & \begin{tabular}[c]{@{}l@{}}p-Willmore\\ DOFs\end{tabular} & \begin{tabular}[c]{@{}l@{}}Conformal\\ DOFs\end{tabular} & Initial $\tau$ & \begin{tabular}[c]{@{}l@{}}Scale Factor $s$\\ \end{tabular} & \begin{tabular}[c]{@{}l@{}}Upper Bound\\ for $\tau$\end{tabular} & \begin{tabular}[c]{@{}l@{}}Final \\ Iteration\end{tabular} & \begin{tabular}[c]{@{}l@{}}Solver Time \\ (1 step)\\ \end{tabular} & Penalty $\varepsilon$ \\ \hline
1 & tri & 23.42k & 2.40 & 105.43k & 58.566k & 3.2e-8 & 1.10 & 5.0e-3 & 130 & 30.1s & \multirow{7}{*}{1.0e-5} \\ \cline{1-11}
\multirow{3}{*}{2} & \multirow{3}{*}{quad} & \multirow{3}{*}{5.344k} & \multirow{3}{*}{3.37} & 16.034k & \multirow{3}{*}{21.382k} & 6.0e-3 & \multirow{3}{*}{1.02} & \multirow{2}{*}{N/A} & 24 & 4.3s &  \\ \cline{5-5} \cline{7-7} \cline{10-11}
 &  &  &  & \multirow{2}{*}{48.114k} &  & 3.2e-4 &  &  & 60 & \multirow{2}{*}{11.5s} &  \\ \cline{7-7} \cline{9-10}
 &  &  &  &  &  & 3.2e-6 &  & 5.0e-3 & 400 &  &  \\ \cline{1-11}
\multirow{2}{*}{3} & \multirow{2}{*}{tri} & \multirow{2}{*}{34.50k} & \multirow{2}{*}{2.05} & \multirow{2}{*}{N/A} & \multirow{2}{*}{86.266k} & \multirow{2}{*}{N/A} & \multirow{2}{*}{N/A} & \multirow{2}{*}{N/A} & \multirow{2}{*}{N/A} & 2.1s &  \\ \cline{11-11}
 &  &  &  &  &  &  &  &  &  & 4.3s &  \\ \cline{1-11}
4 & tri & 9.216k & 6.19 & 41.472k & 23.040k & 5.0e-4 & 1.10 & 5.0e-3 & 29 & 16.6s &  \\ \cline{1-8} \cline{9-12} 
5 & quad & 10.24k & 45.4 & 92.162k & 40.960k & 1.0e3 & \multirow{2}{*}{1.20}  & 5.0e5 & 50 & 22.9s & 1.0e-2 \\ \cline{1-7} \cline{9-12} 
6 & tri & 17.24k & 4.65 & 77.600k & 43.106k & 5.0e-3 &  & \multirow{2}{*}{N/A} & 32 & 27.5s & \multirow{4}{*}{1.0e-5} \\ \cline{1-8} \cline{10-11}
7 & tri & 483.2k & \multirow{2}{*}{257} & N/A & 1.2080m & N/A & N/A &  & N/A & 88.9s &  \\ \cline{1-3} \cline{5-11}
8 & tri & 4.610k &  & 20.691k & N/A & 1.0e3 & 1.20 & 1.0e5 & 32 & 3.8s &  \\ \cline{1-11}
\multirow{2}{*}{9} & \multirow{2}{*}{quad} & \multirow{2}{*}{3.072k} & \multirow{2}{*}{6.55} & \multirow{2}{*}{27.650k} & 12.228k & \multirow{2}{*}{5.0} & \multirow{2}{*}{1.20} & \multirow{2}{*}{N/A} & 30 & 4.79s &  \\ \cline{6-6} \cline{10-12} 
 &  &  &  &  & N/A &  &  &  & 19 & 3.72s & N/A \\ \hline
\multirow{3}{*}{10} & \multirow{3}{*}{quad} & 1.336k, $\times 4$, $\times 16$ & \multirow{3}{*}{3.37} & \multirow{3}{*}{variable} & \multirow{3}{*}{variable} & 5.0e-4 & \multirow{3}{*}{1.02} & \multirow{2}{*}{N/A} & 80 & 4.3s & \multirow{2}{*}{1.0e-5} \\ \cline{3-3} \cline{7-7} \cline{10-11}
 &  & 1.336k, $\times 4$, $\times 16$ &  &  &  & 3.2e-4 &  &  & 100 & \multirow{2}{*}{11.5s} &  \\ \cline{3-3} \cline{7-7} \cline{9-10} \cline{12-12} 
 &  & 0.334k, $\times 4$, $\times 16$ &  &  &  & 3.2e-6 &  & 5.0e-3 & 400 &  & 1.0e-3 \\ \hline
11 & tri & 345.9k & 229 & 1.5568m & 864.87k & 5.0e-4 & N/A & N/A & 1 & 707s & \multirow{3}{*}{1.0e-5} \\ \cline{1-11}
\multirow{3}{*}{12} & \multirow{3}{*}{tri} & \multirow{3}{*}{17.24k} & \multirow{3}{*}{4.65} & 51.733k & \multirow{3}{*}{43.106k} & 5.0e-2 & \multirow{2}{*}{1.10} & \multirow{2}{*}{N/A} & 14 & 9.17s &  \\ \cline{5-5} \cline{7-7} \cline{10-11}
 &  &  &  & \multirow{2}{*}{77.599k} &  & 1.0e-4 &  &  & 36 & \multirow{2}{*}{26.4s} &  \\ \cline{7-10} \cline{12-12} 
 &  &  &  &  &  & 1.0e-7 & 1.05 & 2.0e-4 & 200 &  & 1.0e-4 \\ \hline
\end{tabular}%
}
\label{tab:imp}
\end{table*}

\begin{acks}
The authors would like to acknowledge Dr.~Magdalena Toda and Dr.~Hung Tran for their suggestions regarding numerical models for the p-Willmore flow, and Dr.~Giorgio Bornia who assisted with the implementation of .stl and .ply files in the FEMuS library.  Moreover, we would like to thank the anonymous reviewers for their careful attention and helpful feedback. 

The original meshes in Figures~\ref{fig:funcow}, \ref{fig:moocompare}, \ref{fig:microsoftknot}, \ref{fig:meshcomp}, \ref{fig:statueflow}, \ref{fig:dillo} are courtesy of Keenan Crane (Figure~\ref{fig:funcow}), Microsoft (Figure~\ref{fig:microsoftknot}), and the AIM@SHAPE repository (rest).  The texture in Figure~\ref{fig:statueflow} is courtesy of \\ \url{www.myfreetextures.com}.  The research of the second author was partially supported by the NSF grant DMS-1912902.
\end{acks}

\bibliographystyle{ACM-Reference-Format}
\bibliography{modeling}

\appendix
\green{
\section{Theorem~\ref{thm:conformal} implies Cauchy-Riemann}\label{app:CR}
We show that when $x,y$ are Cartesian coordinates on $\mathbb{R}^2$ and $u:\mathbb{R}^2 \to \mathbb{R}^3$ is an immersion of the $(x,y)$-coordinate plane, the equation $du\, J - N\times du = 0$ expresses the traditional Cauchy-Riemann equations on $\mathbb{C} \cong (T\mathbb{R}^2,J)$.  To see this, let 
\[u(x,y) = \begin{pmatrix} u^1(x,y) \\ u^2(x,y) \\ 0 \end{pmatrix},\]
and consider any constant vector field on $T\mathbb{R}^2$, say $e_1 = (1,0)^T$.  Clearly $J(e_1) = e_2$, so
\begin{equation*}
    du\,J(e_1) = du(e_2) = 
    \begin{pmatrix}
    u_x^1 & u_y^1 \\
    u_x^2 & u_y^2 \\
    0 & 0
    \end{pmatrix}
    \begin{pmatrix}
    0 \\ 1
    \end{pmatrix}= 
    \begin{pmatrix}
    u_y^1 \\ u_y^2 \\ 0
    \end{pmatrix}.
\end{equation*}
Since $N = \begin{pmatrix}0 & 0 & 1 \end{pmatrix}^T$ is normal to the immersion at each point, it follows that
\begin{align*}
  &du\,J(e_1) - N \times du(e_1) \\
  &= 
  \begin{pmatrix}
    u^1_y \\ 
    u^2_y \\
    0
  \end{pmatrix}
  - 
  \begin{pmatrix}
  0 \\ 0 \\ 1
  \end{pmatrix}
  \times
  \begin{pmatrix}
    u^1_x \\
    u^2_x \\
    0
  \end{pmatrix}
  =
  \begin{pmatrix}
  u^1_y + u^2_x \\
  u^2_y - u^1_x \\
  0
  \end{pmatrix}
  =0, \\
\end{align*}
This expression implies the classical Cauchy-Riemann equations,
\begin{align*}
    u^1_x &= u^2_y, \\
    u^1_y &= -u^2_x,
\end{align*}
and it is evident that the expression
\[ \left| du\,J(e_1) - N \times du(e_1) \right|^2 = \left( u^1_y + u^2_x \right)^2 + \left( u^2_y - u^1_x \right)^2, \]
measures the failure of these equations to hold. \greener{ This reflects the fact that, in general, $N \times (\cdot)$ is an ``almost-complex structure'' on $u(\mathbb{R}^2) \subset \mathbb{R}^3$, and an immersion which satisfies the above is ``almost holomorphic''.}
}

\end{document}